\documentclass[invmat,draft,numbook]{svjour}
\usepackage{amsmath,amsfonts,amssymb}

\spnewtheorem*{Theorem}{Theorem}{\bf}{\it}
\spnewtheorem*{Proposition}{Proposition}{\bf}{\it}
\spnewtheorem*{Corollary}{Corollary}{\bf}{\it}
\spnewtheorem*{Lemma}{Lemma}{\bf}{\it}

\def\al{\alpha}
\def\AMN{{\rm A}_N(M)}

\def\be{\beta}
\def\bi{\bar\imath}
\def\bj{\bar\jmath}
\def\BMN{{\rm B}_N(M)}
\def\bx{{\boxed{\phantom{\square}}\kern-.4pt}}

\def\CC{{\mathbb C}}
\def\Cel{{\rm C}_{\ts l}(L)}
\def\Cela{{\rm C}_{\ts\la}(L)}
\def\Cem{{\rm C}_{\ts m}(M)}
\def\Cen{{\rm C}_{\ts n}(N)}
\def\CLl{(\CC^L)^{\ts\ot\ts l}}
\def\CMm{(\CC^M)^{\ts\ot\ts m}}
\def\CNn{(\CC^N)^{\ts\ot\ts n}}
\def\com{\ts,\hskip-.5pt}

\def\de{\delta}
\def\De{\Delta}
\def\dim{\operatorname{dim}\ts}

\def\End{\operatorname{End}\ts}

\def\Fom{F_\Om\ts(M)}

\def\g{{\mathfrak{g}}}
\def\ga{\gamma}
\def\gap{\gamma^{\,\prime}}
\def\ge{\geqslant}
\def\glMN{{\mathfrak{gl}_{N+M}}}
\def\glL{{\mathfrak{gl}_L}}
\def\glM{{\mathfrak{gl}_M}}
\def\glN{{\mathfrak{gl}_N}}

\def\id{{\rm id}}
\def\Idl{{\rm I}_{\ts l}(L)}
\def\Idm{{\rm I}_{\ts m}(M)}

\def\io{\iota}

\def\la{\lambda}
\def\La{\Lambda}
\def\Lac{\La^{\hskip-.5pt\raise-.5pt\hbox{$\circ$}}}
\def\Lan{\La^o}
\def\lap{\la^{\ts\prime}}
\def\Lap{\La^{\ts\prime}}
\def\Lapp{\La^{\ts\prime\prime}}
\def\lcd{\ts,\,\ldots,}
\def\le{\leqslant}
\def\lm{{\la/\ns\mu}}

\def\mi{-}
\def\mup{\mu^{\ts\prime}}
\def\mv{\kern127pt}
\def\mw{\kern-81pt}

\def\ns{\hskip-1pt}
\def\nup{\nu^{\ts\prime}}

\def\Om{\Omega}
\def\op{\oplus}
\def\ot{\otimes}

\def\ph{\varphi}
\def\Pp{P^{\vee}}

\def\Qp{Q^{\,\vee}}

\def\Rb{\,\overline{\hskip-2.5pt R\hskip-4.5pt\phantom{\bar{t}}}\ts}
\def\Rbp{\Rb^{\,\vee}}
\def\Rp{R^{\,\vee}}
\def\RR{{\mathbb R}}
\def\Rt{\ts\widetilde{\ns R\ts}}

\def\sgn{\operatorname{sgn}\ts}
\def\si{\sigma}
\def\so{\mathfrak{so}}
\def\sp{\mathfrak{sp}}

\def\T{{\cal T}}
\def\th{\theta}
\def\Th{\Theta}
\def\ts{\hskip1pt}
\def\Tt{\ts\widetilde{T\ts}}

\def\UgL{\operatorname{U}(\g_L)}

\def\UgMN{\operatorname{U}(\g_{N+M})}
\def\UgN{\operatorname{U}(\g_N)}
\def\UL{\operatorname{U}(\glL)}
\def\UM{\operatorname{U}(\glM)}
\def\UMN{\operatorname{U}(\glMN)}
\def\UN{\operatorname{U}(\glN)}
\def\Up{\Upsilon}
\def\US{\operatorname{U}(\g_N)}
\def\UsoM{\operatorname{U}(\so_M)}
\def\UsoMN{\operatorname{U}(\so_{N+M})}
\def\UspM{\operatorname{U}(\sp_M)}
\def\UspMN{\operatorname{U}(\sp_{N+M})}

\def\Val{V^\ast_{\ns l}}
\def\Vamn{V^{\ts\ast}_{\ns mn}}
\def\Vlm{V_{\la}(\mu)}
\def\Vom{V_\Om\ts(M)}
\def\Vt{\ts\widetilde{\ts V\ts}}

\def\Wal{W^\ast_{\ns l}}
\def\Wamn{W^{\ts\ast}_{\ns mn}}
\def\Wlm{W_{\la}(\mu)}
\def\Wom{W_\Om\ts(M)}

\def\XL{\operatorname{X}\ts(\glL,\si)}
\def\XMN{\operatorname{X}\ts(\glMN,\si)}
\def\XN{\operatorname{X}\ts(\glN,\si)}
\def\xp{x^{\ts\prime}}	

\def\YL{\operatorname{Y}(\glL)}
\def\YMN{\operatorname{Y}(\glMN)}
\def\YN{\operatorname{Y}(\glN)}
\def\YS{\operatorname{Y}(\glN,\si)}
\def\YSL{\operatorname{Y}(\glL,\si)}

\def\ze{\omega}
\def\ZZ{{\mathbb Z}} 

\title{Representations of Twisted Yangians\\
Associated with skew Young Diagrams\\ }
\author{Maxim Nazarov}
\institute{Department of Mathematics, University of York,
York YO10 5DD, England\\
\email{mln1@york.ac.uk}}
\titlerunning{Representations of Yangians}
\authorrunning{Maxim Nazarov}
\date{}
\dedication{\it To Professor I.\,M.\,Gelfand on his 90-th birthday}
\begin{document}
\maketitle

\textwidth=125mm
\textheight=185mm
\parindent=8mm
\evensidemargin=0pt
\oddsidemargin=0pt
\frenchspacing 

\small

\noindent{\bf Abstract.}
Let $G_M$ be either the orthogonal group $O_M$ or
the symplectic group $Sp_M$ over the complex field; 
in the latter case the non-negative integer $M$
has to be even. Classically, the irreducible polynomial
representations of the group $G_M$ are
labeled by partitions $\mu=(\mu_{\ts1},\mu_{\ts2}\ts,\,\ldots)$
such that $\mup_1+\mup_2\le M$ in the case $G_M=O_M$,
or $2\mup_1\le M$ in the case $G_M=Sp_M$.
Here $\mup=(\mup_{\ts1},\mup_{\ts2}\ts,\,\ldots)$ is the partition
conjugate to $\mu$. Let $W_\mu$ be the irreducible polynomial
representation
of the group $G_M$ corresponding to $\mu$.

Regard $G_N\times G_M$ as a subgroup of $G_{N+M}$.
Then take any irreducible polynomial representation
$W_\la$ of the group $G_{N+M}$. 
The vector space
$\Wlm={\rm Hom}_{\,G_M}(\ts W_\mu\ts\com W_\la\ts)$
comes with a natural action of the group $G_N$.
Put $n=\la_1-\mu_1+\la_2-\mu_2+\ldots\,\,$.
In this article, for any standard Young tableau $\Om$ of
skew shape $\lm$ we give a realization of $\Wlm$
as a subspace in the $n\ts$-fold tensor product
$\CNn$, compatible with the action of the group $G_N$.
This subspace is determined as the image of a certain linear operator
$\Fom$ on $\CNn$, given by an explicit formula.

When $M=0$ and $\Wlm=W_\la$ is an irreducible representation of
the group $G_N$, we recover the classical realization of $W_\la$
as a subspace in the space of all traceless tensors in $\CNn$.
Then the operator $F_\Om\ts(0)$ may be regarded as the analogue
for $G_N$ of the Young symmetrizer, corresponding to the
standard tableau $\Om$ of shape $\la\ts$.
This symmetrizer is a certain linear operator on
$\CNn$ with the image equivalent to the irreducible
polynomial representation of the complex general linear group
$GL_N$, corresponding to the partition $\la\ts$. Even in the case 
$M=0$, our formula for the operator $\Fom$ is~new.

Our results are applications of the representation
theory of the twisted Yangian, corresponding to the
subgroup $G_N$ of $GL_N$. This twisted Yangian
is a certain one-sided coideal subalgebra of the Yangian corresponding
to $GL_N$. In particular, $\Fom$ is an intertwining
operator between certain representations of the twisted Yangian in $\CNn$. 


\vskip8pt

\noindent {\bf Mathematics Subject Classification (2000).}
17B35, 17B37, 20C30, 22E46, 81R50.

\vskip8pt

\noindent {\bf Key words.} 
Brauer algebra,
Yangians,
Young symmetrizers.

\normalsize\newpage

\noindent{\bf 0. Brief introduction}

\bigskip\noindent
This article draws on the ideas beginning with the Schur--Weyl duality
\cite{W}, passing through Gelfand--Zetlin bases \cite{GZ1,GZ2} and
continuing through Yangians \cite{D}, with the aim of realizing
explicitly irreducible representations
of the orthogonal or symplectic classical group $G_N$. 
This group corresponds to a symmetric or alternating, 
non-degenerate bilinear form $\langle\ ,\,\rangle$ on 
the $N$-dimensional complex vector space $\CC^N$.
The irreducuble representations of the group
$G_N$ considered in this article
are \textit{polynomial\/}. By definition, they are subrepresentations
of tensor powers of the defining representation acting on $\CC^N$.
This work provides new explicit realization of the
representation of $G_N$ on the
vector space (\ref{1.4}). This vector space
describes the multiplicities in the restriction 
of an irreducible representation of the group $G_{N+M}$ to
the subgroup $G_M$. 
Results follow for the branching rules for restricting
irreducible representations from $G_{N+M}$ to the
subgroup $G_N\times G_M$, see \cite{P}.

Here is an overview of this article.
Section 1 gives an exposition of the principal
results, with detailed references given throughout. Section 2 recalls 
the classical realization \cite{W} of any irreducible polynomial  
representation of the general linear group $GL_N$.
This realization involves elements of the
symmetric group rings, 
known as Young symmetrizers \cite{Y1}.
We also recall the approach to Young symmetrizers
due to Cherednik \cite{C2}. 
Following this approach,
in Section 3 we construct analogues of the
Young symmetrizers for the group $G_N$. 
This construction provides a realization of any
irreducible polynomial representation of the group $G_N$,
more explicit than in \cite{W}.
It is motivated by the
representation theory of Yangians and of their twisted 
analogues \cite{O2}.
The main results concerning branching
rules for the groups $GL_N$ and $G_N$ are stated as 
Theorems 1.6 and 1.8, respectively. 
Theorem 1.6 belongs to Cherednik \cite{C3},
its proof given in Section 4 is new.
Theorem~1.8 is new, its proof is given in
Section 5. 
It is hoped that this work will further
motivate the interest of readers in Yangians,
see \cite{MN,MO,NO}.

A word of explanation is necessary in regard to our scheme of 
referring to theorems, propositions, lemmas and corollaries.
When referring to them,
we indicate the subsections where they respectively appear. 
There will be no more than one of each of them in every subsection,
so our scheme should cause no confusion. For example, Theorems 1.6
and 1.8 mentioned above~are \textit{\it the} theorems
appearing in Subsections 1.6 and 1.8, respectively. 


\section{Main results}\label{S1}

\textbf{1.1.}
Let $\nu=(\nu_{\ts1},\nu_{\ts2}\ts,\,\ldots)$ be any partition of 
a non-negative integer $n$. The parts of $\nu$  are arranged in
the non-increasing order\ts: $\nu_{\ts1}\ge\nu_{\ts2}\ge\ldots\ge0$.
As usual, denote by
$\nup=(\nup_{\ts1},\nup_{\ts2}\ts,\,\ldots)$
the conjugate partition. In particular,
$\nup_{\ts1}$ is the number of non-zero parts of $\nu$.
Take any positive integer $N\ge\nup_{\ts1}$. Let $V_\nu\subset\CNn$
be the irreducible polynomial representation of the complex
general linear group $GL_N$ corresponding to the partition $\nu$. 
We will also regard representations of the group $GL_N$ as modules
over the general linear Lie algebra $\glN$. Then $V_\nu$ is an
irreducible $\glN$-module of the highest weight $(\nu_{\ts1}\lcd\nu_N)$.
Here we choose the  Borel subalgebra in $\glN$ consisting of the
upper triangular matrices, and fix the basis of the diagonal matrix
units $E_{11}\lcd E_{NN}$ in the corresponding Cartan subalgebra of $\glN$.

Now let $\la=(\la_1,\la_2\ts,\,\ldots)$ and $\mu=(\mu_1,\mu_2\ts,\,\ldots)$
be any two partitions. Take any non-negative integer $M$ such that
$\lap_1\le N+M$ and $\mup_1\le M$. Consider the irreducible
representations $V_\la$ and $V_\mu$ of the groups $GL_{N+M}$ and $GL_M\ts$,
respectively. The decomposition $\CC^{\ts N+M}=\CC^N\!\op\CC^M$ provides
an embedding of the direct product $GL_N\times GL_M$ into
$GL_{N+M}$\ts. Consider the vector space

\vskip-16pt
\begin{equation}\label{1.0}
\Vlm={\rm Hom}_{\,GL_M}(\ts V_\mu\ts\com V_\la\ts)\,;
\end{equation}

\vskip4pt\noindent
it comes with a natural action of the group $GL_N$\ts.
This action of $GL_N$ may be reducible.
The vector space $\Vlm$ is non-zero if and only if $\la_i\ge\mu_i$
and $\la'_i-\mu'_i\le N$ for each $i=1\com2\com\ts\ldots$\,; see
\cite[Section I.5]{M}\ts.
In this article, we will consider certain embeddings of $\Vlm$
into the $n\ts$-fold tensor product $\CNn$ where

\vskip-16pt
$$
n=\la_1-\mu_1+\la_2-\mu_2+\ldots\,.
$$

\vskip4pt\noindent
These embeddings will be compatible with the action of the group $GL_N$\ts.

Suppose that $\la_i\ge\mu_i$ for each $i=1\com2\com\ts\ldots$\,\ts.
Consider the {\it skew Young diagram}

\vskip-12pt
$$
\lm=\{\,(i\com j)\in\ZZ^2\ |\ i\ge1,\ \la_i\ge j>\mu_i\,\}\,.
$$

\vskip4pt\noindent
When $\mu=(0\com0\ts,\ts\ldots\ts)$, this is the usual Young diagram
of the partition $\la$. We will employ the standard
graphic representation \cite{M}
of Young diagrams on the plane $\RR^2$ with two matrix style coordinates.
Here the first coordinate increases from top to bottom, while the second
coordinate increases from left to right. The element $(i\com j)\in\lm$
is represented by the unit box with the bottom right corner
at the point $(i\com j)\in\RR^2$. 

The set $\lm$ consists of $n$ elements.
A \textit{standard tableau} of shape $\lm$ is any bijection
$\Om:\lm\to\{1\lcd n\}$ such that $\Om(i\com j)<\Om(i+1\com j)$ and
$\Om(i\com j)<\Om(i\com j+1)$ for all possible $i\com j$\ts.
Graphically, the tableau $\Om$ is represented by placing the
numbers $\Om(i\com j)$ into the corresponding boxes of $\lm$
on the plane $\RR^2$. 
By filling the boxes with the numbers
$1\lcd n$ by rows downwards, from left to right in every row, we get
the \textit{row tableau} $\Om^r$ of shape $\lm$. The
\textit{column tableau} $\Om^c$ of shape $\lm$ is also defined in the
obvious way. Both $\Om^r$ and $\Om^c$ are standard tableaux.
Below we represent $\Om^r$ and $\Om^c$ for the partitions
$\la=(5,\ns3,\ns3,\ns3,\ns3,\ns0,\ns0,\ts\ldots)$ and
$\mu=(3,\ns3,\ns2,\ns0,\ns0,\ts\ldots)$\ts:\!
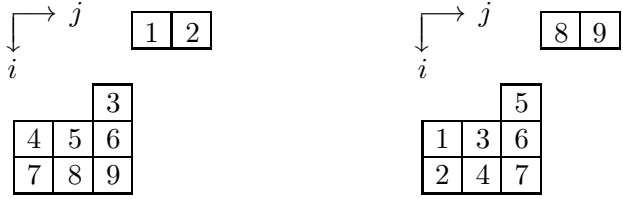
\bigskip

\vbox{
$$
\longrightarrow\,j\mw\kern205pt\longrightarrow\,j\mw\kern130pt
$$
\vglue-24.5pt
$$
\vert\mw\kern232.5pt\vert\mw\kern155.5pt
$$
\vglue-27pt
$$
\bigr\downarrow\mw\kern228.2pt\bigr\downarrow\mw\kern155.6pt
$$
\vglue-16pt
$$
i\mw\kern232.5pt i\mw\kern155.5pt
$$
\vglue-44pt
$$
\phantom{\bx}
\phantom{\bx}
\phantom{\bx}
{\bx}
{\bx}
\kern80pt
\phantom{\bx}
\phantom{\bx}
\phantom{\bx}
{\bx}
{\bx}
$$
\vglue-16pt
$$
\phantom{\bx}
\phantom{\bx}
\phantom{\bx}
\phantom{\bx}
\phantom{\bx}
\kern80pt
\phantom{\bx}
\phantom{\bx}
\phantom{\bx}
\phantom{\bx}
\phantom{\bx}
$$
\vglue-16pt
$$
\phantom{\bx}
\phantom{\bx}
{\bx}
\phantom{\bx}
\phantom{\bx}
\kern80pt
\phantom{\bx}
\phantom{\bx}
{\bx}
\phantom{\bx}
\phantom{\bx}
$$
\vglue-16pt
$$
{\bx}
{\bx}
{\bx}
\phantom{\bx}
\phantom{\bx}
\kern80pt
{\bx}
{\bx}
{\bx}
\phantom{\bx}
\phantom{\bx}
$$
\vglue-16pt
$$
{\bx}
{\bx}
{\bx}
\phantom{\bx}
\phantom{\bx}
\kern80pt
{\bx}
{\bx}
{\bx}
\phantom{\bx}
\phantom{\bx}
$$
\vglue-84pt
$$
\kern45pt{1}\kern10pt{2}\kern134pt{8}\kern9pt{9}
$$
\vglue-3pt
$$
\kern1pt{3}\kern149pt{5}
$$
\vglue-16pt
$$
{4}\kern10pt{5}\kern9pt{6}\kern119pt{1}\kern10pt{3}\kern9pt{6}\kern29pt
$$
\vglue-16pt
$$
{7}\kern10pt{8}\kern9pt{9}\kern119pt{2}\kern10pt{4}\kern9pt{7}\kern29pt
$$
}

\bigskip\smallskip
In Section 2,
for every standard tableau $\Om$ of shape $\lm\ts$ we will
define a certain element $e_\Om$ of the symmetric group ring $\CC S_n$.
If $\mu=(0\com0\ts,\ts\ldots\ts)\ts$, $e_\Om$ is a diagonal
matrix element of the irreducible representation of the group $S_n$  
labeled by the partition $\la\ts$; see (\ref{2.0}).
If $\Om$ is the row or column tableau of shape $\la\ts$,
this matrix element is given explicitly by (\ref{2.1}) and (\ref{2.2}).
When $\mu\neq(0\com0\ts,\ts\ldots\ts)$,
the element $e_\Om\in\CC S_n$ is defined by (\ref{2.85}) and (\ref{2.9}).

We denote by $E_\Om$ the operator on the $n\ts$-fold tensor product $\CNn$,
corresponding to the element $e_\Om\in\CC S_n$ under the action of the
group $S_n$ by permutations of tensor factors. The vector space $\Vlm$
will be realized as the image of the operator $E_\Om$.
Denote this image by $V_\Om\ts$. The subspace $V_\Om$ in $\CNn$
is preserved by the natural action of the group $GL_N$.
Thus $V_\Om$ can be regarded as a representation of $GL_N$. 

\begin{Proposition}
The representations $\Vlm$ and $V_\Om$ of\/ $GL_N$ are equivalent. 
\end{Proposition}

We prove this proposition in Subsection 4.6.
If $\mu=(0\com0\ts,\ts\ldots\ts)$, then
the representation $V_\Om$ of $GL_N$
is equivalent to $V_\la$ by the classical
duality theorem of Schur \cite[Section IV.4]{W}.
Note that for any partitions $\la$ and $\mu$, the
operator $E_\Om$ on the space $\CNn$
does not depend on the integer $M$. 


\medskip\noindent\textbf{1.2.}
There is a description of the operator $E_\Om$
on $\CNn$ of another kind. This description is
obtained by the \textit{fusion procedure\/}, due to Cherednik.
For every $k=1\lcd n$ put $c_k(\Om)=j-i$ if $k=\Om(i\com j)$.
The difference $c_k(\Om)$ is the \textit{content} of the box
occupied by the number $k$ in the tableau $\Om$. In the above example
$n=9$\ts, and the sequences of contents $c_1(\Om^r)\lcd c_9(\Om^r)$ 
and $c_1(\Om^c)\lcd c_9(\Om^c)$ are respectively 
$$
3\com4\com0\com\mi3\com\mi2\com\mi1\com\mi4\com\mi3\com\mi2
\ts\ \ {\rm and}\ \ 
\mi3\com\mi4\com\mi2\com\mi3\com0\com\mi1\com\mi2\com3\com4\ts.
$$

Introduce $n$ complex
variables $t_1(\Om)\lcd t_n(\Om)$ with the constraints
\begin{equation}\label{1.1}
t_k(\Om)=t_l(\Om)
\text{\ \ if $k$ and $l$ occur in the same row of $\Om$\ts.}
\kern10pt
\end{equation}
Alternatively to (\ref{1.1}), as in \cite[Section 2]{NT2},
we can impose the constraints
\begin{equation}\label{1.2}
t_k(\Om)=t_l(\Om)
\text{\ \ if $k$ and $l$ occur in the same column of $\Om$\ts.}
\kern-8pt
\end{equation}
Thus the number of independent variables among
$t_1(\Om)\lcd t_n(\Om)$ equals the number of non-empty rows of the
diagram $\lm$ in the case (\ref{1.1}), or the number of non-empty
columns of $\lm$ in the case (\ref{1.2}).

Order lexicographically
the set of all pairs $(k\com l)$ with $1\le k<l\le n$.
Take the ordered product over this set,
\begin{equation}\label{1.3}
\prod_{1\le k<l\le n}^{\longrightarrow}\ 
\left(1-\frac{P_{kl}}{\ts c_k(\Om)-c_l(\Om)+t_k(\Om)-t_l(\Om)}\ts\right)
\end{equation}

\vskip-5pt
\noindent
where $P_{kl}$ denotes the operator on the vector space $\CNn$
exchanging the $k$th and $l$th tensor factors. Consider
(\ref{1.3})
as a rational function of the constrained variables $t_1(\Om)\lcd t_n(\Om)$.
The next result goes back to \cite{C3}.

\begin{Theorem}
The product {\rm(\ref{1.3})} is regular at
$t_1(\Om)=\ldots=t_n(\Om)$. The value of 
{\rm(\ref{1.3})} at $t_1(\Om)=\ldots=t_n(\Om)$
coincides with the operator $E_\Om$\ts.
\end{Theorem}

Most of the results of \cite{C3} were given without proofs.
In Section 2 of the present article we give all necessary details
of the proof of this theorem. 


\smallskip\medskip\noindent\textbf{1.3.}
The principal aim of this article
is to give analogues of the operator $E_\Om$
for the classical complex Lie
groups $O_N$ and $Sp_N$. Let $G_N$ be one of these two Lie groups.
We will regard $G_N$ as the subgroup in $GL_N$, preserving a non-degenerate
bilinear form $\langle\ ,\,\rangle$ on $\CC^N$, symmetric in the case
$G_N=O_N$ or alternating in the case $G_N=Sp_N$. In the latter case
$N$ has to be even. Throughout this article, we always assume that the
integer $N$ is positive.

The irreducible polynomial representations of $G_N$ are
labeled by the partitions $\nu$ of $n=0\com1\com2\com\,\ldots$ 
such that  $\nup_1+\nup_2\le N$ if $G_N=O_N$, and $2\ts\nup_1\le N$
if $G_N=Sp_N$; see \cite[Sections V.7 and VI.3]{W}.
Denote by $W_\nu$ the irreducible representation of $G_N$ 
corresponding to $\nu$. Take any two distinct numbers $k,l\in\{1\lcd n\}$.
By applying the bilinear
form $\langle\ ,\,\rangle$ to a tensor $w\in\CNn$
in the $k$th and $l$th tensor factors, we obtain a certain tensor
$\widehat{w}\in(\CC^N)^{\ot\ts(n-2)}$. The tensor $w$ is called
\textit{traceless} if $\widehat{w}=0$ for all distinct $k$ and $l$.
Denote by $\CNn_{\,\ts0}$ the subspace in $\CNn$ consisting of 
traceless tensors, this subspace is $G_N\ts$-invariant.
Choose any embedding of the irreducible representation
$V_\nu$ of the group $GL_N$, into the space $\CNn$. Then
$W_\nu$ can be embedded into $\CNn$ as 
$V_\nu\cap\CNn_{\,\ts0}$.

Denote by $\g_N$ the Lie algebra of $G_N$,
so that $\g_N=\so_N$ or $\g_N=\sp_N$.
We will regard $\g_N$ as a Lie subalgebra in $\glN$. 
We will also consider representations of the group $G_N$ as $\g_N$-modules.
The $\g_N$-module $W_\nu$ is irreducible unless $\g_N=\so_N$ and
$2\ts\nup_1=N$, in which case $W_\nu$ is a direct sum of two
irreducible $\so_N$-modules; see Subsection 3.6.
The $\so_N\ts$-module $W_\nu$ may be reducible because the group $O_N$
has two connected components.

Take a non-negative integer $M$ and choose
a non-degenerate bilinear form on the space $\CC^M$,
symmetric in the case $G_N=O_N$ or alternating in the case $G_N=Sp_N$.
In the latter case the integer $M$ has to be even.
Consider the corresponding subgroup $G_M\subset GL_M$. 
If $M=0$, the group $G_M$ consists only of the unit element.
Equip the direct sum 
$
\CC^N\oplus\CC^M=\CC^{\ts N+M}
$
with the bilinear form, which is the sum of the forms
on the direct summands. We get
an embedding of the direct product $G_N\times G_M$ into $G_{N+M}$. 
Consider the irreducible
representations $W_\la$ and $W_\mu$ of the groups $G_{N+M}$ and $G_M\ts$,
respectively. Here we assume that the partitions $\la$ and $\mu$ satisfy
the conditions from \cite[Sections V.7 and VI.3]{W} for the groups
$G_{N+M}$ and $G_M\ts$, respectively, as
described above. Introduce the vector space
\begin{equation}\label{1.4}
\Wlm={\rm Hom}_{\,G_M}(\ts W_\mu\ts\com W_\la\ts)\ts,
\end{equation}
it comes with a natural action of the group $G_N$\ts.
This action of the group $G_N$ may be reducible.
The vector space $\Wlm$ is non-zero if and only if $\la_i\ge\mu_i$ 
and $\la'_i-\mu'_i\le N$ for each $i=1\com2\com\ts\ldots\,\ts$;
see \cite[Proposition 10.1]{P} in the case $G_N=O_N$ and
\cite[Proposition 10.3]{P} in the case $G_N=Sp_N\ts$.
Thus for any given $N$ we have $\Wlm\neq\{0\}$ if and only if 
$\Vlm\neq\{0\}$, see Subsection 1.1 above. Further, for any given $N$
we have the inequality
\begin{equation}\label{1.333}
\dim\Wlm\le\dim\Vlm\,.
\end{equation}
For every standard tableau $\Om$ of shape $\lm\ts$, the results of
the present article
provide a distinguished embedding of the vector space $\Wlm$
into $\Vlm$, compatible with the action of the group $G_N\subset GL_N\ts$;
see \cite{KS}.


\medskip\noindent\textbf{1.4.}
As in the case of the general linear group $GL_N$, suppose
that $\la_i\ge\mu_i$ for all $i=1\com2\com\ts\ldots\,\ts$.
Take the skew Young diagram $\la/\mu$. For any standard tableau $\Om$
of shape $\la/\mu$, we will now construct an embedding of $\Wlm$ into
the tensor product $\CNn$, where $n$ is the number of
elements in the set $\la/\mu$.
This embedding will be compatible with the action of the group $G_N$.

Take any basis $u_1\lcd u_N$ in the vector space $\CC^N$.
Let $v_1\lcd v_N$ be the dual basis in $\CC^N$, so that
$\langle\ts u_i\com v_j\ts\rangle=\de_{ij}$ for $i\com j=1\lcd N$.
The vector
\begin{equation}\label{1.44444444}
w\ts(N)\,=\,\sum_{i=1}^N\,u_i\ot v_i\,\in\,\CC^N\ns\ot\ts\CC^N
\end{equation}
does not depend on the choice
of the basis $u_1\lcd u_N$ and is
invariant under the action of the group $G_N$ on $\CC^N\ns\ot\,\CC^N$.
Introduce the linear operator
\begin{equation}\label{1.45}
Q\ts(N):\ u\ot v\,\mapsto\,\langle\ts u\com v\ts\rangle\cdot w\ts(N)
\end{equation}
in $\CC^N\ns\ot\,\CC^N\ts$; it commutes with the action of the group $G_N$.
The tableau $\Om$ defines the sequence of contents
$c_1(\Om)\lcd c_n(\Om)$. In the case $G_N=O_N$, we will use
the variables $t_1(\Om)\lcd t_n(\Om)$ with the constraints (\ref{1.2}).
In the case $G_N=Sp_N$ it is more convenient to use
the variables $t_1(\Om)\lcd t_n(\Om)$ with
the constraints (\ref{1.1}).
Take the ordered product over the pairs $(k\com l)$
\begin{equation}\label{1.5}
\prod_{1\le k<l\le n}^{\longrightarrow}\ 
\left(1-\frac{Q_{kl}}
{\ts c_k(\Om)+c_l(\Om)+t_k(\Om)+t_l(\Om)+N+M}\ts\right)
\end{equation}

\vskip-4pt
\noindent
where $Q_{kl}$ is the linear operator on $\CNn$,
acting as $Q(N)$ in the $k$th and $l$th tensor factors, and acting
as the identity in the remaining $n-2$ tensor factors.
Here the pairs $(k\com l)$ are
ordered lexicographically, as in (\ref{1.3}).

Throughout this article, we use the following convention. Whenever
the double sign $\pm$ or $\mp$ appears, the upper sign corresponds
to the case of a symmetric form $\langle\ ,\,\rangle$ while the lower
sign corresponds to the case of an alternating form.
Now multiply (\ref{1.5}) by (\ref{1.3}) on
the right, and
consider the result as a rational function of
the constrained variables $t_1(\Om)\lcd t_n(\Om)$.

\begin{Theorem}
a) If $G_N=O_N$ and the variables $t_1(\Om)\lcd t_n(\Om)$ obey
the constraints {\rm(\ref{1.2})}, then
the ordered product of\/ {\rm(\ref{1.5})} and\/ {\rm(\ref{1.3})} is regular
at $t_1(\Om)=\ldots=t_n(\Om)=-\frac12\,$.

\medskip\noindent
b) If $G_N=Sp_N$ and the variables $t_1(\Om)\lcd t_n(\Om)$ obey
{\rm(\ref{1.1})}, then
the ordered product of\/ {\rm(\ref{1.5})} and\/ {\rm(\ref{1.3})} is regular
at $t_1(\Om)=\ldots=t_n(\Om)=\frac12\,$.

\medskip\noindent
c) The operator value $\Fom$ of the ordered product
of\/ {\rm(\ref{1.5})} and\/ {\rm(\ref{1.3})}
at $t_1(\Om)=\ldots=t_n(\Om)=\mp\ts\frac12$
is divisible on the left and on the right~by~$E_\Om\ts$.
\end{Theorem}

Note that unlike $E_\Om\ts$, the operator $\Fom$ on the vector space $\CNn$
may depend on the non-negative integer $M$. Observe that 
$$
c_k(\Om)+c_l(\Om)\ge3-2\la'_1\quad\textrm{if}\quad k\neq l\,.
$$
If the partition $\la$ satisfies the condition
$2\la'_1\le N+M$, every factor in
(\ref{1.5}) is regular at $t_1(\Om)=\ldots=t_n(\Om)=\mp\ts\frac12\,$.
Then by Theorem 1.2 we have 
$$
\Fom=\,
\prod_{1\le k<l\le n}^{\longrightarrow}\, 
\left(1-\frac{Q_{kl}}
{\ts c_k(\Om)+c_l(\Om)+N+M\mp1}\ts\right)
\cdot E_{\Om}\,.
$$
The condition $2\la'_1\le N+M$ is satisfied when
$G_N=Sp_N$, but may be not satisfied when $G_N=O_N$.
The proof of Theorem 1.4 is given
at the end of Subsection 3.4. This proof also
provides an explicit formula for
$\Fom\ts$ in the case $G_N=O_N$ for~any~$\la\ts$.
 
When $G_N=O_N$ and $\Om=\Om^c$, this explicit formula
for the operator $\Fom$ on the space $\CNn$
is particularly simple. Namely, for $G_N=O_N$
\begin{equation}\label{1.6}
F_{\Om^c}(M)\,=\,
\prod_{(k,\ts l)}^{\longrightarrow}\ 
\left(1-\frac{Q_{kl}}
{\ts c_k(\Om^c)+c_l(\Om^c)+N+M-1}\ts\right)
\cdot E_{\Om^c}
\end{equation}

\vskip-5pt
\noindent
where the ordered product is taken over all pairs $(k\com l)$
such that $k$ and $l$ appear in different columns of the tableau $\Om^c$.
For any such pair we have
$$
c_k(\Om^c)+c_l(\Om^c)\ge3-\la'_1-\la'_2\ge3-N-M\,,
$$
so that each of the denominators in (\ref{1.6}) is non-zero.
The operator $E_{\Om^c}$ on $\CNn$ corresponds to the element
$e_{\Om^c}$ of the symmetric group ring $\CC S_n$. If $M=0$,
the element $e_{\Om^c}\in\CC S_n$
can be written explicitly by using~(\ref{2.2}).

When $G_N=Sp_N$, the definition of the operator $\Fom$ on $\CNn$ can
be simplified for the row tableau $\Om=\Om^r$. Namely, for $G_N=Sp_N$
\begin{equation}\label{1.7}
F_{\Om^r}(M)\,=\,
\prod_{(k,\ts l)}^{\longrightarrow}\ 
\left(1-\frac{Q_{kl}}
{\ts c_k(\Om^r)+c_l(\Om^r)+N+M+1}\ts\right)
\cdot E_{\Om^r}
\end{equation}

\vskip-5pt
\noindent
where the ordered product is taken over all pairs $(k\com l)$
such that $k$ and $l$ appear in different rows of the tableau $\Om^r$.
The operator $E_{\Om^r}$ on $\CNn$ corresponds to the element
$e_{\Om^r}$ of the group ring $\CC S_n$.
If $M=0$,
the element $e_{\Om^r}\in\CC S_n$
can be written explicitly by using (\ref{2.1}). 
The simplified formulas (\ref{1.6}) and (\ref{1.7}) will be 
derived at the end of Subsection 3.4.

The vector space $\Wlm$ will be realized as the image of the operator
$\Fom\ts$. Denote this image by $\Wom$.
The subspace $\Wom$ in $\CNn$ is preserved by the natural action of the group
$G_N$ because the operator $\Fom$ commutes with this action by definition.
Thus $\Wom$ can be regarded as a representation of the group $G_N$.

\begin{Proposition}
Representations $\Wlm$ and $\Wom$ of\/ $G_N$ are equivalent. 
\end{Proposition}

The proof is given in Subsection 5.6.
Due to Theorem~1.4, the image $\Wom$ of the operator $\Fom$
is contained in the subspace $V_\Om\subset\CNn$. If $M=0$,
then we have the equality
\begin{equation}\label{1.4444}
W_\Om\ts(0)=V_\Om\cap\CNn_{\,\ts0}
\quad\text{for}\quad
\mu=(0\com0\ts,\ts\ldots\ts)\,;
\end{equation}
see Proposition 3.3.
For general $M$ and $\mu$,
the image $\Wom$ of the operator $\Fom$ may differ from
the intersection $V_\Om\cap\CNn_{\,\ts0}$.
Still our proof of Proposition 1.4 is based on the equality (\ref{1.4444}).
Even when $M=0$ and $\mu=(0\com0\ts,\ts\ldots\ts)\ts$,
our formulas for the operator $F_\Om\ts(0)$ are new. The operator
$F_\Om\ts(0)$ can be regarded as an analogue of the Young symmetrizer
\cite{Y1} for the classical groups $O_N$~and~$Sp_N$ instead of $GL_N$,
see Subsection~3.3 below. This provides a solution
to a problem formulated by Weyl, see \cite[p.\ 149]{W}.

If $\mu\neq(0\com0\ts,\ts\ldots\ts)$, the operator $\Fom$
can also be defined via (\ref{3.46}) and (\ref{3.4444})\ts;
see (\ref{2.85}) and (\ref{2.9}).
Our definition of the operator $\Fom$ is motivated by the
representation theory of Yangians~\cite{MNO}, see below.


\smallskip\medskip\noindent\textbf{1.5.}
By definition, the vector space $\Vlm$ is irreducible under
the natural action of
the subalgebra of $GL_M$\ts-invariants in the universal enveloping algebra
$\UMN\ts$. We denote this subalgebra by $\AMN\ts$;
it coincides with the centralizer of the subalgebra $\UM\subset\UMN\ts$.
In Section 4 we describe the action of the algebra $\AMN$ on $\Vlm$
explicitly, by
using the {\it Yangian} $\YN$ of the general linear Lie algebra $\glN$.   
The Yangian $\YN$ is a deformation of the universal
enveloping algebra of the polynomial current Lie algebra $\glN[x]$
in the class of Hopf algebras, see \cite{D} for instance.

The unital associative algebra $\YN$ has a 
family of generators $T_{ij}^{(a)}$ where $a=1,2,\ts\ldots\ts$ and
$i\ts,\ns j=1\lcd N$. The defining relations for these generators
can be written in terms of the formal power series
\begin{equation}\label{1.31}
T_{ij}(x)=
\de_{ij}\cdot1+T_{ij}^{(1)}x^{-\ns1}+T_{ij}^{(2)}x^{-\ns2}+\,\ldots
\,\in\,\YN\,[[x^{-1}]]\,.
\end{equation}
Here $x$ is the formal parameter. Let $y$ be another formal parameter.  
Then the defining relations in the associative algebra $\YN$
can be written as
\begin{equation}\label{1.32}
(x-y)\cdot[\ts T_{ij}(x)\ts,T_{kl}(y)\ts]\ts=\;
T_{kj}(x)\ts T_{il}(y)-T_{kj}(y)\ts T_{il}(x)\,,
\end{equation}
where $i\com j\com k\com l=1\lcd N\ts$.
The square brackets in (\ref{1.32}) denote the usual commutator.
If $N=1$, the algebra $\YN$ is commutative.
Using the series (\ref{1.31}),
the coproduct $\De:\YN\to\YN\ot\YN$ is defined by
\begin{equation}\label{1.33}
\De\bigl(T_{ij}(x)\bigr)\ts=\ts\sum_{k=1}^N\ T_{ik}(x)\ot T_{kj}(x)\,;
\end{equation}
the tensor product at the right-hand side of the equality (\ref{1.33})
is taken over the subalgebra $\CC[[x^{-1}]]\subset\YN\,[[x^{-1}]]\ts$.
The counit homomorphism $\varepsilon:\YN\to\CC$ is determined by
the assignment $\,\varepsilon:\,T_{ij}(x)\ts\mapsto\ts\de_{ij}\cdot1$.

The antipode ${\rm S}$ on $\YN$ can be defined by using the element
\begin{equation}\label{1.71}
T(x)=\sum_{i,j=1}^N\, E_{ij}\ot
T_{ij}(x)\in\End(\CC^N)\ot\YN\,[[x^{-1}]]\,,
\end{equation}
where the matrix units $E_{ij}$ are regarded as basis elements of the
algebra $\End(\CC^N)\ts$. The formal power series (\ref{1.71}) in
$x^{-1}$ is invertible, because its leading term is $1$.
The anti-automorphism ${\rm S}$ is defined by the assignment
$$
\,\id\ot{\rm S}\ts:\ts T(x)\mapsto T(x)^{-1}\,.
$$
We also use the involutive automorphism $\xi_N$
of the algebra $\YN$ defined by the assignment

\vskip-16pt
\begin{equation}\label{1.51}
\,\id\ot\xi_N:\ts T(x)\mapsto T(-x)^{-1}\,.
\end{equation}

\vskip4pt\noindent
For references and more details on the definition of the Yangian $\YN$
see \cite[Section 1]{MNO}.
Some of these details are also given in Section 4 below.

The defining relations (\ref{1.32}) show that for any $z\in\CC\,$,
the assignment
\begin{equation}\label{tau}
\tau_z:\,T_{ij}(x)\ts\mapsto\,T_{ij}(x-z)
\quad\textrm{for all}\quad
i\com j=1\lcd N
\end{equation}
defines an automorphism $\tau_z$ of the algebra $\YN$. Here the formal
power series $T_{ij}(x-z)$ in $(x-z)^{-1}$ should be re-expanded in $x^{-1}$.
Regard the matrix units $E_{ij}\in\glN$ as generators of the universal
enveloping algebra $\UN$. The relations (\ref{1.32}) show that the assignment
\begin{equation}\label{1.52}
\al_N:\ts T_{ij}(x)\,\mapsto\,\de_{ij}\cdot1-E_{ji}\,x^{-1}
\end{equation}
defines a homomorphism $\al_N:\YN\to\UN$
of associative algebras. By definition,
the homomorphism $\al_N$ is surjective.

By pulling the standard action of the algebra $\UN$ on the space $\CC^N\ts$
back through the composition of the homomorphisms
\begin{equation}\label{eval}
\al_N\circ\,\tau_z:\YN\to\UN\,,
\end{equation}
we obtain a module over the algebra $\YN$, 
called an \textit{evaluation module\/}. 
To indicate the dependence on the
parameter $z\ts$, let us denote this $\YN\ts$-module by $V(z)$.
The operator $E_\Om$ on the vector space
$\CNn$ admits the following interpretation in terms of the tensor products
of evaluation modules over the Hopf algebra $\YN$.
Let $P_0$ be the linear operator on $\CNn$ reversing the order of
the tensor factors. This operator corresponds to the element of
the maximal length in the symmetric group $S_n$.

\begin{Proposition}
The operator $E_\Om\ts P_0$ is an intertwiner of the $\YN$-modules
$$
V(c_n(\Om))\ot\ldots\ot V(c_1(\Om))
\,\ts\longrightarrow\,
V(c_1(\Om))\ot\ldots\ot V(c_n(\Om))\,.
$$
\end{Proposition}

By Proposition 1.5, the image $V_\Om$ of the operator $E_\Om$
is a submodule in the tensor product of
evaluation $\YN\ts$-modules $V(c_1(\Om))\ot\ldots\ot V(c_n(\Om))$. 
This interpretation of $E_\Om$ and $V_\Om$ is due to Cherednik~\cite{C3}.
We obtain Proposition 1.5 as a particular case of Proposition 4.2.

If $M=0$ and
$\mu=(0\com0\ts,\ts\ldots\ts)\ts$,
the image of the operator $E_\Om$ on $\CNn$ is equivalent to
$V_\la$ as a representation of the group $GL_N$. Proposition 1.5 then
turns $V_\la$ into $\YN\ts$-module. The resulting
$\YN\ts$-module can also be obtained from the irreducible
$\glN\ts$-module $V_\la$ by pulling back through the
homomorphism $\al_N\ts$; this is a particular case of Theorem 1.6 below.


\smallskip\medskip\noindent\textbf{1.6.}
Olshanski \cite{O1} has defined a homomorphism from the algebra $\YN$ 
to the subalgebra $\AMN$ of $GL_M$-invariants in $\UMN$,
for each non-negative integer $M$. Along with the centre
of the algebra $\UMN$, the image of this homomorphism generates
the algebra $\AMN$. We will use the following version of this homomorphism,
to be denoted by $\al_{NM}\ts$.

Let the indices $i\com j$ range over the set $\{1\lcd N+M\}$.
Consider the basis of the matrix units $E_{ij}$
in the Lie algebra $\glMN\ts$. We assume that the subalgebras $\glN$ and
$\glM$ in $\glMN$ are spanned by elements $E_{ij}$ where
\begin{equation}\label{ind}
1\le i\com j\le N
\quad\text{and}\quad
N+1\le i\com j\le N+M\ts,
\end{equation}

\vskip-3pt
\noindent
respectively. The subalgebra in the Yangian $\YMN$
generated by $T_{ij}^{(a)}$ where $1\le i\com j\le N\,$,
by definition coincides with the Yangian $\YN$.
Let us denote by $\ph_M$ this natural embedding $\YN\to\YMN$. 
Consider also the involutive automorphism $\xi_{N+M}$
of the algebra $\YMN\ts$, see (\ref{1.51}).
The image of the homomorphism
$$
\al_{N+M}\circ\ts\xi_{N+M}\circ\ts\ph_M:\ts\YN\to\UMN
$$
belongs to the subalgebra $\AMN\subset\UMN\ts$.
Moreover, this image along with the centre
of the algebra $\UMN\ts$, generates the subalgebra $\AMN\ts$.
For the detailed proofs of these two assertions, see \cite[Section 2]{MO}.
In the present article, we use the homomorphism $\YN\to\UMN$
\begin{equation}\label{1.69}
\al_{NM}=\,\al_{N+M}\circ\ts\xi_{N+M}\circ\ts\ph_M\circ\ts\xi_N\ts.
\end{equation}
Note that when $M=0$, the homomorphism (\ref{1.69}) coincides with $\al_N\ts$. 

Take any formal power series $g(x)\in\CC[[x^{-1}]]$ with the leading
term $1$. The assignment

\vskip-14pt
\begin{equation}\label{1.61}
T_{ij}(x)\mapsto\,g(x)\cdot T_{ij}(x)
\end{equation}

\vskip2pt\noindent
defines an automorphism of the algebra $\YN\ts$,
see (\ref{1.31}) and (\ref{1.32}). Put
\begin{equation}\label{1.62}
g_\mu(x)\ =\ \prod_{k\ge1}\ 
\frac{(x-\mu_k+k)(x+k-1)}{(x-\mu_k+k-1)(x+k)}\ .
\end{equation}

\vskip-5pt
\noindent
In the product (\ref{1.62}) over $k$, only finitely many
factors differ from $1$. So $g_\mu(x)$ is a rational function 
of $x$. We have $g_\mu(\infty)=1\ts$; therefore $g_\mu(x)$
expands as a power series in $x^{-1}$ with the leading term $1$. 

As in Subsection 1.1, suppose that
$\lap_1\le N+M$ and $\mup_1\le M$. The space $\Vlm$ comes
with a natural action of the algebra $\AMN\ts$. Regard $\Vlm$ as a module
over the algebra $\YN\ts$, by using the composition of the homomorphism
$\al_{NM}:\YN\to\AMN$ with the automorphism of $\YN$ defined by
(\ref{1.61}), where $g(x)=g_\mu(x)$.
By Proposition~1.5, the image $V_\Om$ of the operator
$E_\Om$ can also be regarded as an $\YN\ts$-module.

\begin{Theorem}
The $\YN$-modules $\Vlm$ and $V_\Om$ are equivalent.
\end{Theorem}

This result goes back to \cite[Theorem 2.6]{C3}.
The proof of Theorem 1.6 is given in Subsections 4.4 to 4.6
of the present article.
The algebra $\AMN$ acts on $\Vlm$ irreducibly;
the central elements of $\UMN$ act on $\Vlm$
as scalar operators. So Theorem 1.6 has a corollary,
see \cite[Section~4]{NT2}.

\begin{Corollary}
The $\YN$-module $V_\Om$ is irreducible.
\end{Corollary}

The Yangian $\YN$ contains the universal enveloping
algebra $\UN$ as a subalgebra. The embedding $\UN\to\YN$ can be defined
by the assignment
\vspace{-8pt}
\begin{equation}\label{4.4}
E_{ij}\mapsto-T_{ji}^{(1)}.
\end{equation}
This embedding provides an action of $\UN$ on the $\YN\ts$-module $\Vlm\ts$.
On the other hand, the vector space $\Vlm$ comes with a natural action
of $\UN$ as a subalgebra of $\UMN\ts$. This natural action of $\UN$ on
$\Vlm$ coincides with its action as a subalgebra in $\YN$,
see Subsection~4.3.

The subspace $V_\Om\subset\CNn$ is preserved by the standard action
of the Lie algebra $\glN$ on $\CNn$, because this subspace is the
image of the operator $E_\Om$ corresponding to
an element of the symmetric group ring $\CC S_n\ts$.
Hence $\UN$ acts naturally on the vector space $V_\Om$ as well.
This natural action of $\UN$ on $V_\Om$ coincides with its action
as a subalgebra in $\YN$, see again Subsection 4.3.

Note that the natural action of the Lie algebra $\glN$ on 
the vector space $\Vlm$ may be reducible.
Using Theorem 1.6 and its Corollary 1.6,
we can identify the vector space $\Vlm$
with the subspace $V_\Om$ in $\CNn$ uniquely, up to multiplication
in $V_\Om$ by a non-zero complex number.
Theorem~1.6 can be regarded as sharpening of Proposition 1.1.
Moreover, we will obtain Proposition 1.1 in the course of
the proof of Theorem 1.6. In the proof of Theorem 1.6,
we will use Proposition 2.4.


\smallskip\medskip\noindent\textbf{1.7.}
Let us now consider the universal enveloping algebra $\UgMN$.
Denote by $\BMN$ the subalgebra of $G_M$\ts-invariants
in $\UgMN$. Then $\BMN$ contains the subalgebra $\UgN\subset\UgMN$.
In the case $\g_N=\sp_N$, $\BMN$ coincides with the
centralizer of the subalgebra $\UspM\subset\UspMN\ts$. 
In the case $\g_N=\so_N$, $\BMN$ is contained in the
centralizer of the subalgebra $\UsoM\subset\UsoMN\ts$,
but may not coincide with the centralizer.
 
In the case $G_N=Sp_N$,
the vector space $\Wlm$ is irreducible under the action of the
algebra $\BMN\ts$.  In the case $G_N=O_N$, the $\BMN\ts$-module $\Wlm$
is either irreducible, or splits as a direct
sum of two irreducible $\BMN\ts$-modules.
It is irreducible if $W_\la$ is irreducible
as a $\so_{N+M}\ts$-module, that is, if $2\la'_1\neq N+M$.
But the condition $2\la'_1\neq N+M$
is not necessary for the irreducibility of the $\BMN\ts$-module $\Wlm$
in the case $G_N=O_N$; see \cite[Proposition 10.1]{P}.
In any case, $\Wlm$ is irreducible under the joint action 
of the algebra $\BMN$ and the subgroup $G_N\subset G_{N+M}$.

In Section~5 we explicitly describe the action of
$\BMN$ in $\Wlm\ts$, by using the {\it twisted Yangian}
$\YS$. Here $\si$ is the involutive automorphism of the Lie algebra $\glN$,
such that $-\si$ is the operator conjugation with respect to the
bilinear form $\langle\ ,\,\rangle$ on $\CC^N$. Then
$\g_N\subset\glN$ is the subalgebra of $\si$-fixed points.
The associative algebra $\YS$ is a deformation of the universal
enveloping algebra of the \textit{twisted polynomial current Lie algebra}
$$
\{A(x)\in\glN[x]:\sigma(A(x))=A(-x)\}\,.
$$
The deformation $\YS$ is not a Hopf algebra, but a coideal
subalgebra in the Hopf algebra $\YN$.
The definition of the twisted Yangian $\YS$
was motivated by the work of
Sklyanin [S] on quantum integrable systems with boundary conditions.
This definition was given by Olshanski \cite{O2} with an assistance
from the author of the present article, see \cite[pp.\ 273--274]{MO}. 

As in Subsection 1.5, let the indices
$i\com j$ range over the set $\{1\lcd N\}$. 
The subalgebra $\YS$ of the associative algebra $\YN$ is defined
in terms of the generating series (\ref{1.31}) as follows.
Let $\Tt(x)$ be the element of the algebra
$\End(\CC^N)\ot\YN\,[[x^{-1}]]$, obtained by applying to $T(x)$
the conjugation with respect to $\langle\ ,\,\rangle$
in the first tensor factor, and by changing $x$ to $-x$.
Then consider the element 
\begin{equation}\label{1.72}
\Tt(x)\,T(x)\in\End(\CC^N)\ot\YN\,[[x^{-1}]]\,. 
\end{equation}
The subalgebra $\YS$ in $\YN$ is generated by the coefficients
of all the formal power series from $\YN\,[[x^{-1}]]$, appearing in
the expansion of the element (\ref{1.72}) relative to the basis
of matrix units $E_{ij}$ in $\End(\CC^N)\ts$.

To give the defining relations for these generators of of
$\YS\ts$, let us introduce the \textit{extended twisted Yangian} $\XN\ts$.
The unital associative algebra $\XN$ has a 
family of generators $S_{ij}^{\ts(a)}$ where $a=1,2,\ts\ldots\,\ts$. Put
\begin{equation}\label{1.771}
S_{ij}(x)=
\de_{ij}\cdot1+S_{ij}^{\ts(1)}x^{-\ns1}+S_{ij}^{\ts(2)}x^{-\ns2}+\,\ldots
\,\in\,\XN\,[[x^{-1}]]\,.
\end{equation}
The defining relations for these generators are given in Subsection 5.1,
using
\begin{equation}\label{1.772}
S(x)\ts=\sum_{i,j=1}^N\, E_{ij}\ot S_{ij}(x)
\,\in\,\End(\CC^N)\ot\XN\,[[x^{-1}]]\,.
\end{equation}
One can define a homomorphism $\pi_N:\XN\to\YS$ by assigning
\begin{equation}\label{piN}
\id\ot\pi_N:S(x)\mapsto\Tt(x)\,T(x)\ts.
\end{equation}
By definition, the homomorphism $\pi_N$ is surjective. 
Further, the algebra $\XN$ has a distinguished family of central elements
$D^{(1)}\ns\com D^{(2)}\ns\com\,\ldots\ts\,$.
These elements of $\XN$ are defined in Subsection 5.1, using the series
\begin{equation}\label{1.773}
D(x)\,=\,1+D^{(1)}x^{-\ns1}+D^{(2)}x^{-\ns2}+\,\ldots
\,\in\,\XN\,[[x^{-1}]]\,.
\end{equation}
By \cite[Theorem 6.4]{MNO} the kernel of the homomorphism $\pi_N$
coincides with the (two-sided) ideal generated by the central elements
$D^{(1)}\ns\com D^{(2)}\ns\com\,\ldots\ts\,$. 

Thus the algebra $\YS$ is defined by the generators $S_{ij}^{\ts(a)}$
satisfying the relation $D(x)=1$ and the 
\textit{reflection equation\/} (\ref{5.1}).
This terminology has been used by physicists;
see \cite[Section 3]{MNO} for the references,
and for more details on the definition of the algebra $\YS$.
In the present article we need the algebra $\XN\ts$, which is
determined by (\ref{5.1}) alone, because this algebra admits
an analogue of the automorphism $\xi_N$ of $\YN\ts$.
By \cite[Proposition 6.5]{MNO} the assignment
\begin{equation}\label{1.751}
\textstyle
\,\id\ot\eta_N:\ts S(x)\mapsto{S(-x-\frac{N}2)}^{-1}
\end{equation}
defines an involutive automorphism $\eta_N$ of the algebra $\XN\ts$.
However, $\eta_N$ does not
determine an automorphism of the algebra $\YS\ts$, because 
the map $\eta_N$ does
not preserve the ideal of $\XN$ generated by the central elements
$D^{(1)}\ns\com D^{(2)}\ns\com\,\ldots\ts\,$; see \cite[Subsection 6.6]{MNO}.

Note that when $z\neq0$, the automorphism $\tau_z$ of the algebra
$\YN$ does not preserve the subalgebra $\YS\subset\YN\ts$; see (\ref{tau}). 
There is no analogue of the automorphism $\tau_z$ for the algebra $\XN$.

Let us now regard $E_{ij}$ as basis vectors of the Lie algebra $\glN$.
The defining relations (\ref{5.1}) of
the algebra $\XN$ imply that the assignment
\begin{equation}\label{1.752}
\be_N:\,S_{ij}(x)\,\mapsto\,\de_{ij}\cdot1-\frac{E_{ji}+\si(E_{ji})}
{\textstyle x\pm\frac12}
\end{equation}

\vskip-2pt\noindent
defines a homomorphism of associative algebras
$
\be_N:\XN\to\US\,;
$
see \cite[Proposition 3.11]{MNO}. According to our general convention,
here the upper sign in $\pm$ corresponds
to the case $G_N=O_N$ while the lower sign corresponds to $G_N=Sp_N$.
The homomorphism $\be_N$ is surjective.
Moreover, $\be_N$ factors through $\pi_N\ts$. Note that the homomorphism
$\YS\to\US$ corresponding to $\be_N\ts$ cannot be obtained from
$\al_N:\YN\to\UN$ by restriction to the subalgebra $\YS$,
because the image of $\YS$ relative to $\al_N$ is not contained in 
the subalgebra $\US\subset\UN$. The link between the homomorphisms
$\al_N$ and $\be_N$ was given by \cite[Proposition 2.4]{N3},
see also \cite[Lemma 3.8]{MN} and Lemma 5.4 of the present article.

The formulas (\ref{1.33}),(\ref{1.71}) and the definition
(\ref{1.72}) imply that for any choice of the symmetric or
alternating non-degenerate form $\langle\ ,\,\rangle$ on $\CC^N$,
the subalgebra $\YS$ in $\YN$ is also a right coideal: 
\begin{equation}\label{1.73}
\De\ts(\ts\YS)\subset\YS\ot\YN\,.
\end{equation}
Although this fact is underlying for our results, it is not directly
used in the present article. Using (\ref{1.73}), for any $\YS\ts$-module 
$W$ and  any $\YN\ts$-module $V$, one
turns the vector space $W\ot V$ into a $\YS\ts$-module again. 
In our case $W$ is going to be the trivial $\YS\ts$-module $\CC$
defined via the restriction to $\YS$ of the counit homomorphism
$\ts\varepsilon:\YN\to\ts\CC\ts$.
 
Let us extend $\si$ to an automorphism of the associative algebra $\UN$.
For any $z\in\CC$, define the \textit{twisted evaluation module\/}
$\Vt(z)$ over the algebra $\YN$
by pulling the standard action of the algebra $\UN$
on the vector space $\CC^N\ts$ back through the 
composition of homomorphisms
\begin{equation}\label{teval}
\si\ts\circ\,\al_N\circ\,\tau_{-z}\ts:\,\YN\to\,\UN\,,
\end{equation}
see (\ref{eval}).
It follows from the definition (\ref{1.72}) that the evaluation
module $V(z)$~and~the twisted evaluation module $\Vt(z)$ over 
$\YN$ have the same restriction to 
$\YS\subset\YN\ts$; see Subsection 5.2 for the explanation.
The linear operator $\Fom$ on
$\CNn$ has the following interpretation,
in terms of the restrictions to $\YS$ of tensor products
of evaluation 
modules over the Hopf algebra $\YN\ts$; see Proposition 1.5. Put
\begin{equation}\label{1.777}
\textstyle
d_k(\Om)=c_k(\Om)+\frac{M}2\mp\frac12
\quad\text{for each}\quad
k=1\lcd n\,.
\end{equation}

\begin{Proposition}
The operator $\Fom$ is an intertwiner of $\YS$-modules
$$
\Vt(d_1(\Om))\ot\ldots\ot\Vt(d_n(\Om))
\,\longrightarrow\ts
V(d_1(\Om))\ot\ldots\ot V(d_n(\Om))\,.
$$
\end{Proposition}

Therefore the image $\Wom$ of the operator $\Fom$
is a submodule in the restriction of the
tensor product of evaluation $\YN\ts$-modules
\begin{equation}\label{1.7777777}
\textstyle
V(d_1(\Om))\ot\ldots\ot V(d_n(\Om))
\end{equation}
to $\YS\subset\YN\ts$. This interpretation is the source
of the definition of $\Fom$ and $\Wom$.
Proof of Proposition 1.7 is given in Subsection~5.3.

If $\mu=(0\com0\ts,\ts\ldots\ts)$ and $M=0$,
the image of the operator $\Fom$ on $\CNn$ is equivalent to
$W_\la$ as a representation of the 
group $G_N$; see (\ref{1.4444}).
Proposition 1.7 turns $W_\la$ into a $\YS\ts$-module.
This $\YS\ts$-module can also be
obtained from the $\g_N\ts$-module $W_\la$ by pulling back through the
homomorphism $\be_N$; this fact is a particular case of Theorem 1.8 below.


\smallskip\medskip\noindent\textbf{1.8.}
Extending the results of \cite{O1} from $\glN$ to
other classical Lie algebras $\g_N=\so_N$ and $\g_N=\sp_N\ts$,
Olshanski \cite{O2} defined
a homomorphism from $\YS$ to the subalgebra $\BMN$
of $G_M\ts$-invariants in $\UgMN\ts$,
for each non-negative integer $M$.
Along with the subalgebra of $G_{N+M}\ts$-invariants in
$\UgMN$, the image of this homomorphism generates
the algebra $\BMN$. We will use the following version of this
homomorphism for the algebra $\XN\ts$,
to be denoted by $\be_{NM}\ts$.

Let the indices $i\com j$ range over $\{1\lcd N+M\}$.
In Subsection~1.6 we chose the basis of the matrix units $E_{ij}$
in the Lie algebra $\glMN\ts$ so that the subalgebras $\glN$ and
$\glM$ in $\glMN$ are spanned by elements $E_{ij}\ts$, where
the indices $i\com j$ satisfy (\ref{ind}).  
Now assume that $\g_N\subset\glN$ and $\g_M\subset\glM\ts$.

Consider the extended twisted Yangian $\XMN\ts$, where $-\si$ is the
conjugation with respect to the form $\langle\ ,\,\rangle$ on
$\CC^{N+M}$. By definition, the subalgebra in $\XMN$
generated by those $S_{ij}^{(a)}$ where $1\le i\com j\le N\,$,
coincides with $\XN$.
Let us denote by $\psi_M$ this natural embedding of the algebra
$\XN$ into $\XMN\ts$.

Consider also the involutive automorphism $\eta_{N+M}$
of $\XMN\ts$, see (\ref{1.751}).
The image of the homomorphism
\begin{equation}\label{beNM}
\be_{N+M}\circ\ts\eta_{N+M}\circ\ts\psi_M:\ts\XN\to\UgMN
\end{equation}
belongs to the subalgebra $\BMN\subset\UgMN\ts$.
This image along with the subalgebra of $G_{N+M}\ts$-invariants in
$\UgMN\ts$, generates the subalgebra $\BMN\ts$.
The proofs of these two assertions are contained in
\cite[Section~4]{MO}. In the present article,
we use the homomorphism $\XN\to\UgMN$
\begin{equation}\label{1.769}
\be_{NM}=\,\be_{N+M}\circ\ts\eta_{N+M}\circ\ts\psi_M\circ\ts\eta_N\ts.
\end{equation}
This is an analogue of the homomorphism (\ref{1.69}). Note that when
$M=0$, the homomorphism (\ref{1.769}) coincides with $\be_N\ts$. 

For any formal power series $g(x)\in\CC[[x^{-1}]]$ with the leading
term $1$, the assignment
\vspace{-6pt}
\begin{equation}\label{1.861}
S_{ij}(x)\mapsto\,g(x)\cdot S_{ij}(x)
\end{equation}

\vskip4pt\noindent
defines an automorphism of the algebra $\XN\ts$,
this follows from (\ref{1.771}) and (\ref{5.1}).
Note that (\ref{1.861}) determines an automorphism of the
quotient $\YS$ of $\XN$ if and only if $g(x)=g(-x)\ts$;
see Subsection~5.1.

As in Subsection 1.3, suppose that the partitions $\la$
and $\mu$ satisfy the conditions from \cite{W} for the groups $G_{N+M}$
and $G_{M}\ts$, respectively. Consider the vector space $\Wlm$ which comes
with a natural action of the algebra $\BMN\ts$. Regard $\Wlm$ as a module
over the algebra $\XN\ts$, using the composition of the homomorphism
$\be_{NM}:\XN\to\BMN$ with the automorphism of the algebra $\XN\ts$,
defined by (\ref{1.861}) where 
\begin{equation}\label{gmM}
\textstyle
g(x)=g_\mu(x-\frac{M}2\pm\frac12)\ts;
\end{equation}
see (\ref{1.62}).
By Proposition~1.7, the image $\Wom$ of the operator
$\Fom$ can also be regarded as an $\YS\ts$-module.

\begin{Theorem}
The action of the algebra $\XN$ on $\Wlm$ factors through the homomorphism
$\pi_N:\XN\to\YS\ts$.
The $\YS$-modules $\Wlm$ and $\Wom$ are equivalent.
\end{Theorem}

Together with the explicit construction of the subspace $\Wom$ in $\CNn$,
this analogue of Theorem 1.6 is the main result of the present article.
The proof of Theorem 1.8 is given in Subsections 5.4 to 5.6.


\smallskip\medskip\noindent\textbf{1.9.}
The twisted Yangian $\YS$ contains the enveloping
algebra $\UgN$ as a subalgebra. The embedding $\UgN\to\YS$ can be defined
by
\begin{equation}\label{1.844}
E_{ij}+\si(E_{ij})\ts\mapsto\,-\ts\pi_N\ts(S_{ji}^{(1)})\,,
\end{equation}
see \cite[Proposition 3.12]{MNO}.
This embedding yields an action of $\UgN$ on the $\YS\ts$-module
$\Wlm\ts$, see the first statement of Theorem 1.8.
On the other hand, the vector space $\Wlm$ comes with a natural action
of the subgroup $G_N\subset G_{N+M}$.
The corresponding action of $\UgN$ on
$\Wlm$ coincides with its action as a subalgebra in $\YS$, see
Subsection~5.7.

The subspace $\Wom\subset\CNn$ is preserved by the standard action
of the group $G_N$ on $\CNn$.
So $\UgN$ acts naturally on the vector space $\Wom$ as well.
This natural action of $\UgN$ on $\Wom$ coincides with its action
as a subalgebra in $\YS$, see again Subsection~5.7.
Theorem~1.8 thus agrees with Proposition 1.4.

The $G_M\ts$-invariant elements of $\UgMN$
act on the space $\Wlm$ as scalar operators. 
In the case $G_N=Sp_N$, the $\YS\ts$-module $\Wlm$ is therefore irreducible. 
In the case $G_N=O_N$,  $\Wlm$ is irreducible under the joint action of
the algebra $\YS$ and of the group $G_N$. In our proof
of Theorem 1.8 we construct a surjective linear operator
$\Wlm\to\Wom\ts$, which intertwines the actions of both $\YS$ and $G_N$.
Hence our proof of Theorem~1.8 has a corollary, cf.\ Corollary~1.6.

\begin{Corollary}
a) If $G_N=Sp_N$, the $\YS$-module $\Wom$ is irreducible.

\smallskip\noindent
b) If $G_N=O_N$, $\Wom$ is irreducible under action of\/ $\YS$
and $O_N$.
\end{Corollary}

In the case when $G_N=O_N$,
the $\YS\ts$-module $\Wom$ is either irreducible, or splits
as a direct sum of two irreducible submodules; see the beginning of 
Subsection~1.7 for more details.

Using Corollary 1.9, we can identify $\Wlm$
with the subspace $\Wom$ in $\CNn$ uniquely, up to multiplication
in $\Wom$ by a non-zero complex number. This identification is
compatible with the action of $G_N$. Thus we sharpen Proposition 1.4. 
Moreover, we will obtain Proposition 1.4 in the course of
the proof of Theorem 1.8. In this proof, we will use Proposition~3.5.

By Theorem 1.4, $\Wom$ is a vector subspace in the image $V_\Om$
of the operator $E_\Om$. As explained in Subsection 1.6,
we can identify the vector space $\Vlm$ with $V_\Om$ uniquely,
up to multiplication in $V_\Om$ by a non-zero complex number.
Using the identification of $\Wlm$ with $\Wom$ as above,
we obtain a distinguished embedding of the vector space $\Wlm$
into $\Vlm$; see the inequality (\ref{1.333}). 
This embedding is compatible with the
action of the subgroup $G_N\subset GL_N$, and depends on the choice of
standard tableau $\Om$ of shape~$\lm$. 
This result supports a thesis of Cherednik \cite{C3}, that
Yangians are ``hidden symmetries'' of the classical
representation theory. 

The subspace $V_\Om$ in $\CNn$ can also be regarded as a submodule
in the tensor product of evaluation $\YN\ts$-modules (\ref{1.7777777})\ts;
this follows by setting $z=\frac{M}2\mp\frac12$ in Proposition 4.2.
Denote this $\YN\ts$-submodule by $\Vom\ts$.
By Proposition 1.7, then $\Wom$ is a submodule
in the restriction of $\Vom$ to the subalgebra $\YS\subset\YN\ts$.

Let us now consider the case when $N=2$ and $G_N=Sp_{\ts2}\ts$. Then the
equality of dimensions in (\ref{1.333}) is attained,
see \cite[Proposition 10.3]{P}. In this case, the restriction of the
$\operatorname{Y}(\mathfrak{gl}_{\ts2})\ts$-module $\Vom$ to the subalgebra 
$
\operatorname{Y}(\mathfrak{gl}_{\ts2},\si)
\subset
\operatorname{Y}(\mathfrak{gl}_{\ts2})\ts
$
is irreducible, and coincides with $\Wom\ts$.
So $\Wom$ is an irreducible 
$\operatorname{Y}(\mathfrak{gl}_{\ts2})\ts$-module in this
case. By using Theorem~1.8 along with \cite[Corollary 2.13]{NT1}, 
one then derives \cite[Theorem~5.2]{M1}
which describes $\Wlm$ as
$\operatorname{Y}(\mathfrak{gl}_{\ts2})\ts$-module,
and which has been pivotal for this work of Molev.

In another special case when $N=2$ and $G_N=O_2\ts$, our results on the 
vector space $\Wlm$ are different from those implied by
\cite[Theorem~3.2]{M2} and \cite[Theorem~2.3]{M3}\ts:
unlike Molev, in this case we work with the non-connected Lie group
$O_{\ts2+M}$ rather than with its Lie algebra $\so_{2+M}\ts$. In any case, 
the methods of this article are different from those of \cite{M1,M2,M3}.


\section{Young symmetrizers}\label{S2}

\textbf{2.1.}
We begin this section with recalling several classical facts \cite{Y1,Y2}
about  the irreducible representations of the symmetric group  $S_l$ over
the complex field $\CC$.  These representations are labeled by
partitions of $l$. We will identify partitions with their Young diagrams.
Denote by $U_\la$ the irreducible representation of $S_l$
corresponding to the partition $\la$. We will also regard
representations of the group $S_l$ as modules over the group ring $\CC S_l$.
Fix the chain $S_1\subset S_2\subset\ldots\subset S_l$
of subgroups with the standard embeddings.

There is a decomposition of the space $U_\la$ into a direct sum of
one-dimensional subspaces, labeled by the {standard tableaux} of shape $\la$.
The one-dimensional subspace $U_\La\subset U_\la$ corresponding to
a standard tableau $\La$ is defined as follows. For any
$m\in\{1\lcd l-1\}$ take the tableau obtained from $\La$ by
removing the numbers $m+1\lcd l$. Let the Young diagram $\mu$
be the shape of the resulting tableau.
Then the subspace $U_\La$ is contained in an irreducible
$\CC S_m$-submodule of $U_\la$ corresponding to $\mu$. Any basis of
$U_\la$ formed by vectors $u_\La\in U_\La$ is called a 
\textit{Young basis}. Fix an $S_l$-invariant inner product
$(\,\,,\,)$ on $U_\la$. All the subspaces
$U_\La\subset U_\la$ are then pairwise orthogonal. We choose the vectors
$u_\La\in U_\La$ so that $(\ts u_\La,u_\La\ts)=1$.

For any standard tableau $\La$, denote by 
$S_\La$ (respectively by $S^{\ts\prime}_\La$) the subgroup
in $S_l$ preserving the collections of numbers appearing
in every row (every column) of the tableau $\La$.
Then introduce the elements of $\CC S_l$
$$
p_\La=\sum_{s\in S_\La}s
\ \quad\text{and}\ \quad
q_\La=\sum_{s\in S^{\ts\prime}_\La}s\cdot\sgn s\,.
$$

\vskip-5pt
\noindent
The product $p_\La\ts q_\La\in\CC S_l$ is the \textit{Young symmetrizer}
corresponding to $\La$. By \cite{Y1} the normalized product
$p_\La\ts q_\La\cdot\dim U_\la\ts/\ts l\ts!$ is a minimal idempotent in the
group ring $\CC S_l$. The left ideal in $\CC S_l$ generated by the element
$p_\La\ts q_\La$ is equivalent to the representation $U_\la\ts$, under the
action of the group $S_l$ on this ideal via left multiplication. 

For any standard tableau $\La$ consider the diagonal matrix element of the
representation $U_\la$ corresponding to the vector $u_\La$\,,
\begin{equation}\label{2.0}
e_\La=
\sum_{s\in S_l}\,
(\,u_\La\com\ts s\ns\cdot\ns u_\La\,)
\,s\,\in\,\CC S_l\,.
\end{equation}

\vskip-6pt
\noindent
We have the equality
\begin{equation}\label{2.01}
e_\La^{\ts2}=e_\La\cdot l\ts!\ts/\dim U_\la\ts.
\end{equation}

There is an explicit formula for the element $e_\La$
of the group ring $\CC S_l\ts$.  
This formula is particularly simple when $\La$ is the {row tableau}
$\La^r$, or the {column tableau} $\La^c\ts$. 
Using the lemma from \cite[Section IV.2]{W}, one~can~obtain

\begin{Proposition}
a) There are equalities of elements of\/ $\CC S_l$
\begin{eqnarray}
e_{\La^r}&=&\,p_{\La^r}\ts q_{\La^r}\ts p_{\La^r}/
{\la_1!\ts\la_2!\ts\ldots}\ ,
\label{2.1}
\\
e_{\La^c}&=&\,q_{\La^c}\ts p_{\La^c}\ts q_{\La^c}/
{\la'_1!\ts\la'_2\ts!\ts\ldots}\ .
\label{2.2}
\end{eqnarray} 
b) There exist invertible elements $p\com q\in\CC S_l$ such that 
$$
p_{\La^r}\ts q_{\La^r}\ts p_{\La^r}=p_{\La^r}\ts q_{\La^r}\ts p
\ \quad\text{and}\ \quad
q_{\La^c}\ts p_{\La^c}\ts q_{\La^c}=q_{\La^c}\ts p_{\La^c}\ts q\,.
$$
\end{Proposition}

\smallskip
For any $k=1\lcd l-1$ let $s_k\in S_l$ be the transposition of $k$ and $k+1$.
An expression for the matrix element $e_\La$ with arbitrary $\La$
can be obtained from either of (\ref{2.1}) and (\ref{2.2}) by using the
formulas \cite{Y2} for the action of the generators $s_1\lcd s_{l-1}$ 
of the group $S_l$ on
the vectors of the Young basis. Fix any standard tableau $\La$. For every
$k=1\lcd l$ let $c_k=c_k(\La)$ be the {content}
of the box occupied by $k$ in $\La$.
Consider the tableau $s_k\La$ obtained from $\La$ by exchanging the numbers
$k$ and $k+1$. The resulting tableau may be non-standard. This
happens exactly when $k$ and $k+1$ stand (next to each other) in the same 
row or column of $\La$. But $c_k\neq c_{k+1}$ always; put
$h=(c_{k+1}-c_k)^{-1}$. If $s_k\La$ is non-standard, we have $h=1$
or $h=-1$.

So far the vector $u_\La$ has been determined
up to a multiplier $z\in\CC$ with $|z|=1$. Due to \cite{Y2} all the 
vectors of the Young basis can be further normalized so that
for any standard tableau $\La$ and $k=1\lcd l-1$
\begin{equation}\label{2.3}
s_k\cdot u_\La=
\left\{
\begin{array}{ll}
h\ts u_\La+\sqrt{1-h^2}\, u_{s_k\La}
&\ \textrm{if}\ s_k\La\ \textrm{is standard;}\\[2pt]
h\ts u_\La
&\ \textrm{otherwise.}
\end{array}
\right.
\end{equation} 
This normalization determines all the vectors of the Young basis
up to a common multiplier $z\in\CC$ with $|z|=1$. If the tableau $s_k\La$
is standard, then
\begin{equation}\label{2.55555}
\sqrt{1-h^2}\ts\,u_{s_k\La}=(s_k-h)\,u_\La
\end{equation}
due to (\ref{2.3}).
Then by the definition (\ref{2.0}) we have the relation
\begin{equation}\label{2.4}
(1-h^2)\,e_{s_k\La}=(s_k-h)\,e_\La\ts(s_k-h).
\end{equation} 
For any standard tableau $\La$ there is a sequence
of transpositions $s_{k_1}\lcd s_{k_b}\ns$ such that 
$\La=s_{k_b}\ldots\ts s_{k_1}\La^r$ and the tableaux
$s_{k_a}\ldots\ts s_{k_1}\La^r$ are standard for all $a=1\lcd b$\,;
see for instance \cite[Section 2]{N1}. Using (\ref{2.4}) repeatedly,
one can derive from (\ref{2.1}) an explicit formula for the element
$e_\La\in S_l$. Another explicit formula for $e_\La$
can be derived in a similar way from (\ref{2.2}).


\medskip\noindent\textbf{2.2.}
There is an expression for the element $e_\La\in\CC S_l$
of another kind. This expression is obtained by
the so-called fusion procedure \cite{C2}.
For every two distinct indices
$i\com j\in\{1\lcd l\ts\}$ introduce the rational
function of $x\com y\in\CC$
\begin{equation}\label{2.45}
f_{ij}(x\com y)=1-\frac{(\ts i\ts j\ts)}{x-y}\ ,
\end{equation}
valued in $\CC S_l$\,; 
here $(\ts i\ts j\ts)\in S_l$ is the transposition of $i$ and $j$.
As direct calculation shows, these rational functions satisfy the relations
\begin{equation}\label{2.5}
f_{ij}(x\com y)\,f_{ik}(x\com z)\,f_{jk}(y\com z)=
f_{jk}(y\com z)\,f_{ik}(x\com z)\,f_{ij}(x\com y)
\end{equation}
for pairwise distinct indices $i\com j\com k$. 
For pairwise distinct $i\com j\com m\com n$ we also have

\vspace{-16pt}
\begin{equation}\label{2.55}
f_{ij}(x\com y)\,f_{mn}(z\com w)=
f_{mn}(z\com w)\,f_{ij}(x\com y)\,.
\end{equation}

Now take $l$ complex variables $t_1\lcd t_l$. Order lexicographically
the set of all pairs $(i\com j)$ with $1\le i<j\le l$. The
ordered product over this set,
\begin{equation}\label{2.6}
\prod_{1\le i<j\le l}^{\longrightarrow}\ 
f_{ij}(\ts c_i+t_i\com c_j+t_j\ts)
\end{equation}
is a rational function of $t_1\lcd t_l$ with values in $\CC S_l$.
This function depends only on the differences $t_i-t_j$.
Denote by $\T_\La$ the set of tuples $(t_1\lcd t_l)$ such that
$t_i=t_j$ whenever the numbers $i$ and $j$ appear in the same row of
the standard tableau $\La$. Alternatively, we can choose $\T_\La$
to be the set of all tuples $(t_1\lcd t_l)$ such that
$t_i=t_j$ whenever the numbers $i$ and $j$ appear in the same column of $\La$.
The following proposition goes back to \cite{C2}.

\begin{Proposition}
Restriction to $\T_\La$ of the rational function {\rm(\ref{2.6})}
is regular at $t_1=\ldots=t_l$. The value of this restriction at
$t_1=\ldots=t_l$ coincides with the element $e_\La\in\CC S_l$.
\end{Proposition}

In its present form, Proposition 2.2 has been proved in
\cite[Section~2]{N2}. The proof actually provides an explicit
formula for the element $e_\La\in\CC S_l$\,, different
from those obtained by using (\ref{2.1}) or (\ref{2.2}).


\medskip\noindent\textbf{2.3.}
We need a generalization of Proposition 2.2
to standard tableaux of skew shapes \cite{C2}.
Take any $m\in\{0\lcd l-1\}$\ts. 
Let $\Up$ be standard tableau obtained from $\La$ by removing
the boxes with the numbers $m+1\lcd l$. Let $\mu$ be the shape of
the tableau $\Up$. Define a standard tableau $\Om$ of the skew shape
$\lm$ by setting $\Om(i\com j)=\La(i\com\ns j)-m$ for all
$(i\com j)\in\lm$\ts. Every standard tableau $\Om$ of shape $\lm$
is obtained from a certain $\La$ in this way.

Put $n=l-m$. Denote by $\io_m$ the embedding of the symmetric group $S_n$
into $S_l$ as a subgroup preserving the subset $\{m+1\lcd l\}$\ts;
we extend the map $\io_m$ to $\CC S_n$ by linearity.
Denote by $S_{mn}$ the subgroup $S_m\times\io_m(S_n)$ in $S_l$.
Introduce the linear map
\begin{equation}\label{2.85}
\th_m:\CC S_l\to\CC S_{mn}:s\mapsto
\left\{
\begin{array}{ll}
s&\ \textrm{if}\ s\in S_{mn}\ts;\\[2pt]
0&\ \textrm{otherwise.}
\end{array}
\right.
\end{equation}
By definition, the element $e_\La\in\CC S_l$ is divisible on the
left and on the right by $e_{\ts\Up}\in\CC S_m$. Hence there exists an
element $e_\Om\in\CC S_n$ such that
\begin{equation}\label{2.9}
\th_m(e_\La)=e_{\ts\Up}\cdot\io_m(e_\Om).
\end{equation}
The element $e_\Om\in\CC S_n$ does not depend on the choice of
standard tableau $\Up$ of the shape $\mu$\ts, because the boxes with
the numbers $1\lcd m$ have in $\La$ and $\Up$ the same contents\ts;
see (\ref{2.4}). The generalization of Proposition 2.2 from $\La$
to $\Om$ is based on the following simple observation.

\begin{Lemma}
The image under the map\/ $\th_m$ of the product\/
 {\rm(\ref{2.6})} equals
\begin{equation}\label{2.7}
\prod_{1\le i<j\le m}^{\longrightarrow}
f_{ij}(\ts c_i+t_i\com c_j+t_j\ts)
\hskip5pt\cdot\hskip-6pt
\prod_{m<i<j\le l}^{\longrightarrow} 
f_{ij}(\ts c_i+t_i\com c_j+t_j\ts)\,.
\end{equation}
\end{Lemma}

\begin{proof}
The pairs $(i\com j)$ in the product (\ref{2.6}) 
are ordered lexicographically:
\begin{equation}\label{2.8}
(1\com2)\lcd(1\com l)\,,\ts(2\com3)\lcd(2\com l)\lcd\ldots\lcd(l-1\com l)\,.
\end{equation}
Now expand the product (\ref{2.6}) as a sum of the products of transpositions
$(i_1\ts j_1)\ts\ldots(i_d\ts j_d)=s$
with coefficients from the field $\CC(t_1\lcd t_l)$\ts;
the sum is taken over subsequences
$(i_1\com j_1)\lcd(i_d\com j_d)$ in the sequence (\ref{2.8}).
Let $(i_b\com j_b)$ be the first pair in the subsequence such that
$i_b\le m<j_b$\ts, we suppose the pair exists. Let $(i_c\com j_c)$
be the last pair in the subsequence such that $i_c=i_b$. Note that
$j_c\ge j_b>m$, while for any pair $(i_a\com j_a)$ with $a<b$ we have
$i_a<j_a\le m$. Hence $s(i_b)=j_c>m$ for $i_b\le m$\ts,
and $\th_m(s)=0$.\qed
\end{proof}

Now consider the ordered product on the right-hand side of (\ref{2.7}),
\begin{equation}\label{2.10}
\prod_{m<i<j\le l}^{\longrightarrow} 
f_{ij}(\ts c_i+t_i\com c_j+t_j\ts)\,.
\end{equation}

\begin{Corollary}
Restriction of\/ {\rm(\ref{2.10})} to $\T_\La$ 
is regular at $t_{m+1}=\ldots=t_l$. The value of this restriction at
$t_{m+1}=\ldots=t_l$ coincides with $\io_m(e_\Om)\in\CC S_l$.
\end{Corollary}

\begin{proof}
Consider the restriction of the rational function (\ref{2.7}) to
$\T_\La$. By Proposition 2.2 and Lemma 2.3, this
restriction is regular at $t_1=\ldots=t_l\ts$; the value of this
restriction at $t_1=\ldots=t_l$ is $\th_m(e_\La)$. By Proposition 2.2
applied to the tableau $\Up$ instead of to $\La$, 
the value at $t_1=\ldots=t_m$ of the restriction of the product
over $1\le i<j\le m$ in (\ref{2.7}) is $e_{\ts\Up}$. Therefore the
value at $t_{m+1}=\ldots=t_l$ of restriction to $\T_\La$ of the
product (\ref{2.10}) is the element $\io_m(e_\Om)$,
determined by the relation (\ref{2.9}).
\qed
\end{proof}

The proof of Proposition~2.2 from \cite[Section 2]{N2} implies
the following.

\begin{Proposition}
The element $e_\La$ is divisible on the left and right by $\io_m(e_\Om)$.
\end{Proposition}


\noindent\textbf{2.4.}
For any two standard tableaux $\La$ and $\Lap$ of the same shape $\la\ts$,
consider the matrix element of the representation $U_\la$ corresponding
to the pair of vectors $u_\La$ and $u_{\La'}$ of the Young basis,
$$
e_{\La\ts\La'}=
\sum_{s\in S_l}\,
(\,u_\La\com\ts s\cdot\ns u_{\La'}\,)
\,s\,\in\,\CC S_l\,.
$$
If $\La=\Lap\ts$, then we have the equality $e_{\La\ts\La'}=e_\La$
by the definition (\ref{2.0}).
Now take any $m\in\{0\lcd l-1\}$. Put $n=l-m$ as in Subsection 2.3.
Consider the element $\th_m(\ts e_{\La\ts\La'})\in S_{mn}\ts$,
see the definition (\ref{2.85}).
Let $\Up$ and $\Up^{\ts\prime}$ be the standard tableaux obtained by removing
the boxes with the numbers $m+1\lcd l$ from the tableaux $\La$ and
$\Lap$, respectively.
 
\begin{Lemma}
We have\/ $\th_m(\ts e_{\La\ts\La'})=0$ unless\/ 
$\Up$ and\/ $\Up^{\ts\prime}$ are of the same shape.
\end{Lemma}

\begin{proof}
By definition, the element $e_{\La\ts\La'}\in\CC S_l$ is divisible
by $e_{\ts\Up}\in\CC S_m$ on the left, and by $e_{\ts\Up'}\in\CC S_m$ on
the right. The image $\th_m(\ts e_{\La\ts\La'})\in\CC S_{mn}$ inherits
these two divisibility properties.
If the tableaux $\Up$ and\/ $\Up^{\ts\prime}$ are not of the same shape,
the irreducible representations $U_\Up$ and $U_{\Up'}$ of $S_m$ are not 
equivalent, and $e_{\ts\Up}\ts s\,e_{\ts\Up'}=0$ for any $s\in S_m\,$.
\qed
\end{proof}

Consider the permutational action
of the symmetric group $S_n$ on the tensor product $\CNn$.
We denote by $E_\Om$ the linear operator on $\CNn$ 
corresponding of the element $e_\Om\in\CC S_n$.
In the notation of Subsection 1.2, we have $c_k(\Om)=c_{m+k}$ for
$k=1\lcd n$\ts. Moreover, when $(t_1\lcd t_l)\in\T_\La$ we can put
$t_k(\Om)=t_{m+k}$ for $k=1\lcd n$\ts. Then we obtain Theorem 1.2
as a reformulation of Corollary 2.3, see the definition (\ref{2.45}). 

Consider the image $V_\Om$ of the operator $E_\Om\ts$. 
Recall that the standard tableau $\Om$ is of shape $\la/\mu\ts$.
Take any non-negative integer $M$ such that $\lap_1\le N+M$
and $\mup_1\le M$. Consider the vector space $\Vlm$ defined by~(\ref{1.0}). 

\begin{Proposition}
If\/ $\Vlm\neq\{0\}$, then $V_\Om\neq\{0\}$.
\end{Proposition}

\begin{proof}
Take the operator $E_\La$ on the vector space
$(\CC^{N+M})^{\ot\,l}$, corresponding
to the element $e_\La\in\CC S_l\ts$.
Realize the irreducible
representation $V_\la$ of the group $GL_{N+M}$
from (\ref{1.0}) as the image $V_\La$ of $E_\La\,$.
Split $(\CC^{N+M})^{\ot\,l}$
into the direct sum of subspaces, obtained
from the subspaces $(\CC^M)^{\ot\ts k}\ot\ts(\CC^N)^{\ot\ts(l-k)}$ by some
permutations of the $l$ tensor factors; here $k=0\lcd l$.
Each of these subspaces is preserved by the action of the subgroup
$GL_M\subset GL_{N+M}\ts$.
The vector space $V_\La$ is the sum of the images of all
these subspaces with respect to $E_\La$, but
only the images with $k=m$ may contribute to
$$
\Vlm=\ts{\rm Hom}_{\,GL_M}(\ts V_\mu\ts\com V_\La\ts)\ts.
$$

Further, consider the projections of 
$V_\La$ onto the direct summands of 
$(\CC^{N+M})^{\ot\,l}$ obtained
from the subspaces $(\CC^M)^{\ot\ts k}\ot\ts(\CC^N)^{\ot\ts(l-k)}$ by some
permutations of the $l$ tensor factors.
Again, irreducible representations of $GL_M$ equivalent to $V_\mu$
may occur only in the projections with $k=m$.
Let us denote by $I_m$ the projector onto the direct summand 
\begin{equation}\label{4.7}
\CMm\ot\CNn\subset(\CC^{N+M})^{\ot\,l}\,.
\end{equation}
Let $P_s$ be the operator on $(\CC^{N+M})^{\ot\,l}$
corresponding to permutation $s\in S_l\ts$.

Now suppose that $\Vlm\neq\{0\}$. By the above argument, then
there exist permutations $s^{\ts\prime}$ and $s^{\ts\prime\prime}$ 
in $S_l$ such that
\begin{equation}\label{2.40}
{\rm Hom}_{\,GL_M}(\, V_\mu\,\com\ts 
I_m\,P_{s'}\,E_\La\,P_{s''}\cdot\ts
\CMm\ot\CNn\ts)\neq\{0\}\,.
\end{equation}
The product $P_{s'}\,E_\La\,P_{s''}$ in (\ref{2.40}) can be 
written as a linear combination of the operators $E_{\La'\La''}$
on $(\CC^{N+M})^{\ot\,l}$,
corresponding to the matrix elements $e_{\La'\La''}\in\CC S_l\,$.
In this linear combination the coefficients are taken from $\CC\ts$,
while $\Lap$ and $\Lapp$ range over all standard tableaux of shape $\la\ts$. 
By (\ref{2.40}), there exists at least
one pair of tableaux $\Lap$ and $\Lapp$ such that
\begin{equation}\label{2.4040}
{\rm Hom}_{\,GL_M}(\, V_\mu\,\com\ts 
I_m\,E_{\La'\La''}\cdot\ts
\CMm\ot\CNn\ts)\neq\{0\}\,.
\end{equation}

Restriction of the operator $I_m\,E_{\La'\La''}$
to the subspace (\ref{4.7}) coincides with the restriction
of the operator on $(\CC^{N+M})^{\ot\,l}$, corresponding
to the element $\th_m(\ts e_{\La'\La''})\in\CC S_{mn}\ts$.
Consider the tableaux $\Up^{\ts\prime}$ and $\Up^{\ts\prime\prime}$,
obtained by removing the boxes with the numbers $m+1\lcd l$ from the 
tableaux $\Lap$ and $\Lapp$, respectively.
Due to Lemma 2.4, the inequality (\ref{2.4040}) implies that
$\Up^{\ts\prime}$ and $\Up^{\ts\prime\prime}$ are of the same shape.
But then $e_{\La'\La''}=e_{\La'}\ts e$ for some invertible element
$e\in\CC S_{mn}\ts$, see (\ref{2.55555}). Then
$$
E_{\La'\La''}\cdot\ts\CMm\ot\CNn\ts=\,
E_{\La'}\cdot\ts\CMm\ot\CNn\,.
$$
By applying the relation (\ref{2.9}) to the tableau $\Lap$ instead
of $\La\ts$, the inequality (\ref{2.4040})
now implies that the tableau $\Up^{\ts\prime}$ is of shape $\mu\ts$.
Moreover, the left-hand side of (\ref{2.4040}) then
equals $V_{\Om'}$ for some standard tableau $\Om^{\ts\prime}$
of skew shape $\lm\ts$. By (\ref{2.4}), the
inequality $V_{\Om'}\neq\{0\}$ implies that $V_\Om\neq\{0\}\ts$.
\qed
\end{proof}

The space $\Vlm$
comes with an action of the subgroup $GL_N\subset GL_L\ts$.
The subspace $V_\Om\subset\CNn$ is preserved by the action of 
the group $GL_N\ts$. Let us now consider $\Vlm$ and $V_\Om$ as
representations of the group $GL_N$. In Subsection 4.6
we will prove that
these representations are equivalent, as stated in Proposition~1.1.
Note that the operator $E_\Om$ on $\CNn$ does not depend on $M$.
It is well known that the dimension of the vector space $\Vlm$
is the same~for all integers $M$ such that
$\lap_1\le N+M$ and $\mup_1\le M$; see for instance \cite[Section I.5]{M}.


\smallskip\medskip\noindent\textbf{2.5.}
In this subsection, we collect a few results which
we need for the proof of Theorem 1.6.
Let us keep fixed a standard tableau $\La$ of non-skew shape 
$\la\ts$. Here $\la$ is a partition of $l$. 
For every $k=1\lcd l$ we denote $c_k=c_k(\La)\ts$.
Consider the rational functions (\ref{2.45}) with
pairwise distinct indices $i\com j\in\{1\lcd l+1\ts\}$;
these functions take values in the group ring $\CC S_{\ts l+1}$.
Take the element $e_\La\in\CC S_l$ defined by (\ref{2.0}).
Consider the image of $e_\La$ under the
embedding $\io_{\ts1}:\ts\CC S_l\to\CC S_{\ts l+1}$, see the
beginning of Subsection~2.3. 
For the proof of the following proposition, see \cite[Section 2]{N2}.

\begin{Proposition}
We have equality of rational functions in\/ $x\ts$,
valued in\/ $\CC S_{\ts l+1}$
$$
\ f_{12}(x\com c_1)\,\ldots\ts f_{1,l+1}(x\com c_l)\cdot\io_{\ts1}(e_\La)
\,=\,
\biggl(1\ts-\ts\sum_{k=1}^l\,\frac{(1\,k+1)}x\ts\biggr)
\cdot\io_{\ts1}(e_\La)\,.
$$
\end{Proposition}

Now for each $m=0\lcd l-1$ define a linear map
$\ga_m:\CC S_{l+1}\to\CC S_{l+1}$ as follows.
By definition, for any group element $s\in S_{l+1}$
\begin{equation}\label{gam}
\ga_m(s)\,=\,
\left\{
\begin{array}{ll}
s&\ \textrm{if}\ \ s(1)\neq2\lcd m+1\ts;\\[2pt]
0&\ \textrm{otherwise.}
\end{array}
\right.
\end{equation}

\begin{Lemma}
For any $z_1\lcd z_l\in\CC$ we have equality
of rational fuctions in $x$
$$
\ga_m\ts(\ts f_{1,l+1}(x\com z_l)\,\ldots\ts f_{12}(x\com z_1))=
f_{1,l+1}(x\com z_l)\,\ldots\ts f_{1,m+2}(x\com z_{m+1})\,.
$$
\end{Lemma}

\begin{proof}
Let us expand the product 
$f_{1,l+1}(x\com z_l)\,\ldots\ts f_{12}(x\com z_1)$
as a sum of the products of transpositions 
$(1\ts i_a)\ldots(1\ts i_1)$ with coefficients from $\CC(x)\ts$;
here the sum is taken over all subsequences $i_1\lcd i_a$ in the 
sequence $2\lcd l+1\ts$. By the definition (\ref{gam}) we have $\ga_m(1)=1$,
while for $a\ge1$
$$
\ga_m\ts((1\ts i_a)\ldots(1\ts i_1))\,=\,
\left\{
\begin{array}{ll}
(1\ts i_a)\ldots(1\ts i_1)&\ \textrm{if}\ \ i_1>m+1\ts;\\[2pt]
\,0&\ \textrm{otherwise.\qed}
\end{array}
\right.
$$
\end{proof}

We will need a reformulation of this lemma.
For each $m=0\lcd l-1$ define a linear map
$\gap_m:\CC S_{l+1}\to\CC S_{l+1}$ as follows:
for $s\in S_{l+1}\ts$,
\begin{equation}\label{gamp}
\gap_m(s)\,=\,
\left\{
\begin{array}{ll}
s&\ \textrm{if}\ \ s^{-1}(1)\neq2\lcd m+1\ts;\\[2pt]
0&\ \textrm{otherwise.}
\end{array}
\right.
\end{equation}

\begin{Corollary}
For $z_1\lcd z_l\in\CC$ we have equality
of rational fuctions in $x$
$$
\gap_m\ts(\ts f_{12}(x\com z_1)\,\ldots\ts f_{1,l+1}(x\com z_l))=
 f_{1,m+2}(x\com z_{m+1})\,\ldots\ts f_{1,l+1}(x\com z_l)\,.
$$
\end{Corollary}

This result is derived from Lemma 2.5 by using the
anti-automor\-phism of the group ring
$\CC S_{l+1}$, such that $s\mapsto s^{-1}$ for every group element $s\ts$.


\section{Traceless tensors}\label{S3}

\textbf{3.1.}
For any positive integer $L$ choose a non-degenerate bilinear
form $\langle\ ,\,\rangle$ on the vector space $\CC^L$, symmetric or
alternating. In the latter case $L$ has to be even. 
For any positive integer $l$ consider 
the commutant of the image of the classical group $G_L=O_L$ or
$G_L=Sp_L$ in the operator algebra $\End(\CLl)$. 
This commutant is called the \textit{Brauer centralizer algebra\/},
see \cite[Section V.2]{W}.
We will denote this algebra by $\Cel$.

Consider the linear operator $Q(L)$ on $\CC^L\ot\CC^L$,
see (\ref{1.45}).
This operator commutes with the action of the group $G_L$ in $\CC^L\ot\CC^L$.
For any distinct indices $i\com j\in\{1\lcd l\ts\}$ let $Q_{ij}$ be
the operator on $\CLl$ acting as $Q(L)$ in the $i$th and $j$th
tensor factors, and acting as the identity
in the remaining $l-2$ tensor factors.
Note that $Q_{ij}=Q_{ji}\ts$. Here we also have
\begin{equation}\label{3.99}
Q_{ij}^{\,2}=L\cdot Q_{ij}
\quad\textrm{for any distinct}\quad
i\com j\in\{1\lcd l\ts\}\,.
\end{equation}
Consider the image
of the symmetric group ring $\CC S_l$ in $\End(\CLl)$.
The algebra $\Cel$ is generated by this image and all the
operators $Q_{ij}$ on $\CLl$ with $i<j\ts$.
Let $P_{ij}$ be the operator on $\CLl$ corresponding to
transposition $(i\,j)\in S_l\ts$.
We have the equalities 
\begin{equation}\label{3.9}
Q_{ij}\,(1\mp P_{ij})=0
\quad\textrm{for any distinct}\quad
i\com j\in\{1\lcd l\ts\}\,.
\end{equation}
According to our general convention,
the upper sign in $\mp$ corresponds
to the case $G_L=O_L\ts$, while the lower sign corresponds to the case 
$G_L=Sp_L\ts$.

For any two distinct indices 
$i\com j\in\{1\lcd l\ts\}$ consider the rational function of $x\com y\in\CC$

\vskip-16pt
\begin{equation}\label{3.45}
R_{\ts ij}(x\com y)=1-\frac{P_{ij}}{x-y}
\end{equation}

\vskip-2pt\noindent
which takes values in $\End(\CLl)\ts$;
this is the \textit{Yang rational R-matrix\/}. 
The function (\ref{3.45}) corresponds to (\ref{2.45}) under the
action of the symmetric group $S_l$ on $\CLl$. Due to (\ref{2.5}) we have
the relation
\begin{equation}\label{3.5}
R_{\ts ij}(x\com y)\,R_{\ts ik}(x\com z)\,R_{\ts jk}(y\com z)=
R_{\ts jk}(y\com z)\,R_{\ts ik}(x\com z)\,R_{\ts ij}(x\com y)
\end{equation}
for any pairwise distinct indices $i\com j \com k\,$;
it is called the \textit{Yang--Baxter relation\/}.
Note that
%
\begin{equation}\label{3.55}
R_{\ts ij}(x\com y)\,R_{\ts ji}(y\com x)=1-\frac1{(x-y)^{\ts2}}\,.
\end{equation}

The operator $Q(L)$ on $\CC^L\ot\CC^L$
can be obtained from the permutation operator on $\CC^L\ot\CC^L$
by conjugation in any one of the two tensor factors
relative to the bilinear form $\langle\ ,\,\rangle\ts$.
For any $i\neq j$ introduce the functions
\begin{equation}\label{3.555}
\Rt_{ij}(x\com y)=1+\frac{Q_{ij}}{x+y}
\textrm{\hskip10pt and\hskip10pt}
\Rb_{ij}(x\com y)=1-\frac{Q_{ij}}{x+y+L}\,\ts;
\hskip-10pt
\end{equation}
then we have
\begin{equation}\label{3.6}
\Rt_{ij}(x\com y)\ \Rb_{ij}(x\com y)=1\,.
\end{equation}
As $Q_{ij}=Q_{ji}\ts$, we also have the relations
\begin{equation}\label{3.65}
\Rt_{ij}(x\com y)=\Rt_{ji}(y\com x)
\qquad\textrm{and}\qquad
\Rb_{ij}(x\com y)=\Rb_{ji}(y\com x)\,.
\end{equation}
By applying to the relation (\ref{3.5}) the conjugation relative to the form
$\langle\ ,\,\rangle\ts$ in the $i$th tensor factor, and
by changing $x$ to $-x$, we get the relation
\begin{equation}\label{3.7}
\ \qquad
\Rt_{\ts ik}(x\com z)\,\Rt_{\ts ij}(x\com y)\,R_{\ts jk}(y\com z)=
R_{\ts jk}(y\com z)\,\Rt_{\ts ij}(x\com y)\,\Rt_{\ts ik}(x\com z)\,.
\end{equation}
By using (\ref{3.6}), we obtain from (\ref{3.7}) the relation
\begin{equation}\label{3.8}
\ \,\qquad
\Rb_{\ts ij}(x\com y)\,\Rb_{\ts ik}(x\com z)\,R_{\ts jk}(y\com z)=
R_{\ts jk}(y\com z)\,\Rb_{\ts ik}(x\com z)\,\Rb_{\ts ij}(x\com y)\,.
\end{equation}
Using (\ref{3.6}) and (\ref{3.65}), we also obtain from (\ref{3.7}) 
the relation
\begin{equation}\label{3.85}
\ \,\qquad
\Rt_{\ts ij}(x\com y)\,R_{\ts ik}(x\com z)\,\Rb_{\ts jk}(y\com z)=
\Rb_{\ts jk}(y\com z)\,R_{\ts ik}(x\com z)\,\Rt_{\ts ij}(x\com y)\,.
\end{equation}
The relations (\ref{3.5})\com\ts(\ref{3.55}) and
(\ref{3.6}) to (\ref{3.85}) are equalities
of functions that take their values in the Brauer centralizer algebra $\Cel$.


\medskip\noindent\textbf{3.2.}
The irreducible modules over the algebra $\Cel$, or equivalently,
the irreducible representations of the group
$G_L$ appearing in the tensor product 
$\CLl$, are labeled by the partitions of
$l\com l-2\com\,\ldots$ satisfying the conditions
from \cite{W} for $G_L$ as described in 
Subsection~1.3.
In the present article, we consider only the irreducible
$\Cel$-modules which are labeled by partitions of $l$.
The corresponding irreducible representations of the group $G_L$
appear in the subspace

\vskip-12pt
\begin{equation}\label{3.0}
\CLl_{\,\ts0}\subset\CLl
\end{equation}
of traceless tensors. The images of the groups $S_l$ and $G_L$ in
$\End(\CLl_{\,\ts0})$ span the commutants of each other,
see \cite[Sections V.7 and VI.3]{W} again.

By definition, all
the operators $Q_{ij}$ on $\CLl$ vanish on the subspace (\ref{3.0}).
Denote by $\Idl$ the two-sided ideal in $\Cel$ generated by all the operators
$Q_{ij}\,$. The quotient algebra $\Cel/\,\Idl$ can be identified 
with the image of the group ring
$\CC S_l$ in the operator algebra $\End(\CLl_{\,\ts0})\ts$.
The image of the element $P_{ij}\in\Cel$ in the quotient algebra 
is then identified with the operator on the subspace (\ref{3.0})
corresponding to $(i\ts j)\in S_l\ts$. 

Take any partition $\la$ of $l\ts$ satisfying the above mentioned
conditions. The irreducible representation
$W_\la\subset\CLl_{\,\ts0}$ of the group $G_L$ corresponds to the
irreducible representation $U_\la$ of the group $S_l\ts$.
We will regard $U_\la$ as a $\Cel$-module
using the canonical homomorphism
\begin{equation}\label{3.22}
\Cel\,\longrightarrow\,\Cel/\,\Idl\,.
\end{equation}
Then $U_\la$ is the irreducible $\Cel$-module
corresponding to the partition $\la\ts$.

Fix any standard tableau $\La$ of shape $\la$. As in Section 2,
let $c_k=c_k(\La)$ be the content of the box
occupied by $k$ in $\La$.  Consider the product 
\begin{equation}\label{3.1}
\prod_{1\le i<j\le l}^{\longrightarrow}\,
\Rb_{\ts ij}(\ts c_i+t_i\com c_j+t_j\ts)
\ \cdot
\prod_{1\le i<j\le l}^{\longrightarrow}\,
R_{\ts ij}(\ts c_i+t_i\com c_j+t_j\ts)
\end{equation}
as a rational function of the complex variables
$t_1\lcd t_l$ which takes values in the algebra $\Cel$.
Here the pairs $(i\com j)$ with $1\le i<j\le l$ are ordered
lexicographically, as usual.
Denote by $\T_\La$ the set of tuples $(t_1\lcd t_l)$ such that
$t_i=t_j$ whenever the numbers $i$ and $j$ appear in $\La$
the same column for $G_L=O_L\ts$, or in the same row for $G_L=Sp_L\ts$.
Denote by $E_\La$ the linear operator on $\CLl$
corresponding to the element $e_\La\in\CC S_l\ts$, see (\ref{2.0}).

\begin{Proposition}
Restriction to $\T_\La$ of the rational function 
{\rm(\ref{3.1})} is regular at $t_1=\ldots=t_l=\mp\ts\frac12$.
The value $F_\La$ of this restriction at $t_1=\ldots=t_l=\mp\ts\frac12$ 
is divisible on the left and on the right
by the operator $E_\La$ on $\CLl$.
\end{Proposition}

\begin{proof}
Denote by $\ze$ the involutive
anti-automorphism of the algebra $\Cel$, such that
all the elements $P_{ij}$ and $Q_{ij}$ are $\ze\ts$-invariant.
The element $E_{\La}$ of $\Cel$ corresponds to (\ref{2.0}).
Therefore $E_\La$ is also $\ze\ts$-invariant.
Applying $\ze$ to the product (\ref{3.1}) just
reverses the ordering of the factors. The initial ordering can
then be restored by using the relations (\ref{3.5}) and (\ref{3.8}).
So any value of the function (\ref{3.1}) is $\ze\ts$-invariant.
Thanks to this observation, it suffices to prove the divisibility
of $F_\La$ by $E_\La$ only on the right.

Take $l$ complex variables $x_1\lcd x_l\ts.$
Choose any $k\in\{1\lcd l-1\}$ and write $\xp_1\lcd\xp_l$ for
the sequence of variables obtained  by exchanging the terms 
$x_k$ and $x_{k+1}$ in  the sequence $x_1\lcd x_l$. Using
(\ref{3.5}), we obtain
$$
P_{k,k+1}\ts R_{\ts k+1,k}\ts(\ts x_{k+1}\com x_k)\,\cdot
\prod_{1\le i<j\le l}^{\longrightarrow}\,
R_{\ts ij}(\ts x_i\com x_j\ts)\,
$$
\vskip-6pt
\begin{equation}\label{3.2}
=\,\prod_{1\le i<j\le l}^{\longrightarrow}\,
R_{\ts ij}(\ts\xp_i\com\xp_j\ts)\,\cdot\,
R_{\ts k,k+1}\ts(x_k\com x_{k+1})\,P_{k,k+1}\,.
\end{equation}
Using (\ref{3.8}) along with the second relation in (\ref{3.65}),
we obtain the equality
$$
P_{k,k+1}\ts R_{\ts k+1,k}\ts(x_{k+1}\com x_k)\,\cdot
\prod_{1\le i<j\le l}^{\longrightarrow}\,
\Rb_{\ts ij}(\ts x_i\com x_j\ts)
$$
\vskip-6pt
\begin{equation}\label{3.3}
=\,\prod_{1\le i<j\le l}^{\longrightarrow}\,
\Rb_{\ts ij}(\ts\xp_i\com\xp_j\ts)
\cdot R_{\ts k+1,k}\ts(x_{k+1}\com x_k)\,P_{k,k+1}\,.
\end{equation}

Now suppose that the index $k\in\{1\lcd l-1\}$ is chosen so that
the tableau $s_k\La$ is standard. Then $|\,c_k-c_{k+1}|>1$.
Set $x_i=c_i+t_i$ for each $i=1\lcd l$ in the two equalities
(\ref{3.2}) and (\ref{3.3}). Then these two equalities show that
when multiplying (\ref{3.1}) on the left by
\begin{equation}\label{3.35}
P_{k,k+1}\ts R_{\ts k+1,k}\ts(c_{k+1}\ns+t_{k+1}\com c_k+t_k)\,,
\end{equation}
and dividing on the right by
\begin{equation}\label{3.3535}
R_{\ts k,k+1}\ts(c_k+t_k\com c_{k+1}\ns+t_{k+1})\,P_{k,k+1}\,,
\end{equation}
we get the analogue of (\ref{3.1}) for the standard tableau $s_k\La$
instead~of~$\La\ts$.
As $|\,c_k-c_{k+1}|>1\ts$, the functions (\ref{3.35}) 
are regular at $t_k=t_{k+1}\ts$. Moreover, the values of the functions
(\ref{3.35}) at $t_k=t_{k+1}\ts$ are invertible.
By (\ref{2.4})~and~(\ref{3.55}),
$$
{\textstyle
P_{k,k+1}\ts R_{\ts k+1,k}\ts(c_{k+1}\mp\frac12\,\com\ts c_k\mp\frac12\ts)
\ts\,E_\La}
$$
$$
=\ {\textstyle
E_{s_k\La}\,
R_{\ts k,k+1}\ts(c_k\mp\frac12\,\com\ts c_{k+1}\mp\frac12\ts)\,P_{k,k+1}\,.}
$$

\vskip4pt
\noindent
Hence it suffices to prove Proposition 3.2 for the
tableau $s_k\La$ instead of $\La$.

First consider the case $G_L=O_L$. Here
we select the upper sign in $\mp\,$.
There is a sequence
of transpositions $s_{k_1}\lcd s_{k_b}$ such that 
$s_{k_b}\ldots\ts s_{k_1}\La=\La^c$ and the tableaux
$s_{k_a}\ldots\ts s_{k_1}\La^r$ are standard for all $a=1\lcd b$\,;
see Subsection~2.1. Therefore it suffices to prove Proposition~3.2 only for
$\La=\La^c$.
Then the set $\T_{\La^c}$ consists of all tuples $(t_1\lcd t_l)$
such that $t_i=t_j$ whenever $i$ and $j$ appear in
the same column of $\La^c\ts$.

Suppose that the numbers $k$ and $k+1$ appear in the
same column of $\La^c$.  Using the relations (\ref{3.5}) we can demonstrate
that in (\ref{3.1}), the second ordered product over $1\le i<j\le l$
is divisible on the left by the function 
$R_{\ts k,k+1}\ts(t_k\ns+c_k\com t_{k+1}\ns+c_{k+1})$.
Restriction of this function to $\T_{\La^c}$ equals $1-P_{k,k+1}$.
Now reorder the pairs $(i\com j)$ in the first ordered
product over $1\le i<j\le l$ in (\ref{3.1}) as follows. 
For every $a=1\lcd\la_1\ts$ first come all the pairs $(i\com j)$ where
both $i$ and $j$ occur in the $a$th column of the tableau $\La^c$.
These pairs are ordered between themselves lexicographically.
After them come all the pairs $(i\com j)$ where $i$ occurs in the $a$th column,
while $j$ occurs in the $b$th column, for all $b>a$.
These pairs are again ordered lexicographically.
After these come the pairs $(i\com j)$ where $i$ occurs in the
$b$th column, for some $b>a$.
This reordering involves only transpositions of commuting factors,
so the value of the product (\ref{3.1}) will not change.

Now take any $a\in\{1\lcd\la_1\}\ts$. Let $c\com c+1\lcd d$
be the consecutive numbers appearing in the $a$th column of $\La^c$. 
In the restriction to $\T_{\La^c}$ of the reordered product (\ref{3.1}),
the product of all the factors succeeding
\begin{equation}\label{3.10}
\prod_{c\le i<j\le d}^{\longrightarrow}\,
\Rb_{\ts ij}(\ts c_i+t_i\com c_j+t_j\ts)\,,
\end{equation}
is divisible on the left by $1-P_{k,k+1}$ whenever $c\le k<d$.
This follows from the relations (\ref{3.8}).
Moreover, then the product of the succeeding factors is divisible on
the left by $1-P_{ij}$ for $c\le i<j\le d$. Using the relations
(\ref{3.9}) with the upper sign in $\mp\,$, we now eliminate
all the products (\ref{3.10}) from the restriction
of (\ref{3.1}) to $\T_{\La^c}$, consecutively for
$a=1\lcd \la_1$. The restriction to $\T_{\La^c}$ 
of the product of the remaining factors in (\ref{3.1}) is
regular at $t_1=\ldots=t_l=-\frac12\,$. Thus we get
the first statement of Proposition~3.2. Moreover, we get 
the explicit formula
\begin{equation}\label{3.11}
F_{\La^c}=\,
\prod_{(i,\ts j)}^{\longrightarrow}\, 
\left(1-\frac{Q_{ij}}
{\ts c_i(\La^c)+c_j(\La^c)+L-1}\ts\right)
\cdot E_{\La^c}
\ \ \text{for}\ \ 
G_L=O_L
\end{equation}

\vskip-4pt
\noindent
where  the ordered product is taken over all pairs $(i\com j)$
such that $i$ and $j$ appear in different columns of $\La^c$.
For any such pair we have
$$
c_i(\La^c)+c_j(\La^c)\ge3-\la'_1-\la'_2\ge3-L\,,
$$
so that each of the denominators in (\ref{3.11}) is non-zero.
The operator $E_{\La^c}$ on $\CLl$ corresponds to the element
$e_{\La^c}\in\CC S_l\ts$, and we used Proposition~2.2.
The operator $F_{\La^c}$ given by (\ref{3.11})
is evidently divisible
by $E_{\La^c}$ on the right.

The proof of Proposition 3.2
in the case $G_L=Sp_L$ is much simpler.
As we observed in Subsection~1.4,
in this case every factor
$\Rb_{\ts ij}(\ts c_i+t_i\com c_j+t_j\ts)$
in the product (\ref{3.1}) is regular at $t_1=\ldots=t_l=\frac12\,$,
for any choice of the tableau $\La$. Using Proposition~2.2,
we get from (\ref{3.1}) the explicit formula
\begin{equation}\label{3.2222}
F_\La=\,
\prod_{1\le i<j\le l}^{\longrightarrow}\, 
\left(1-\frac{Q_{ij}}
{\ts c_i+c_j+L+1}\ts\right)
\cdot E_{\La}
\quad\text{for}\quad
G_L=Sp_L\,.
\end{equation}

\vskip-4pt
\noindent
Thus $F_\La$ is divisible by $E_\La$ on the right.
Thanks to the observation made in
the beginning of the proof, $F_\La$ is
also divisible by $E_\La$ on the left.
\qed
\end{proof}

Observe that in the case $G_L=O_L$,
for any tableau $\La$ and any two distinct
numbers $i\com j\in\{1\lcd l\}$ we have $c_i+c_j\ge3-L$
unless both $i$ and $j$
appear in the first column of $\La$. Therefore the factor
$\Rb_{\ts ij}(\ts c_i+t_i\com c_j+t_j\ts)$
in the product (\ref{3.1}) may have a pole at $t_1=\ldots=t_l=-\ts\frac12\,$,
only if both $i$ and $j$ appear in the first column of $\La$.
However, this observation does not facilitate significantly the proof
of Proposition~3.2. Our proof also has a

\begin{Corollary}
If  the tableau $s_k\La$ is standard for some $k\in\{1\lcd l-1\}$, then
$$
P_{\ts k,k+1}\ts R_{\ts k+1,k}\ts(\ts c_{k+1}\com c_k)\ts\,F_\La\,=\,
F_{s_k\La}\,R_{\ts k,k+1}\ts(\ts c_k\com c_{k+1})\,P_{\ts k,k+1}\,.
$$
\end{Corollary}

In the case $G_L=O_L$, our proof of Proposition 3.2 also provides the
formula (\ref{3.11}) for the operator $F_{\La^c}$.
In the case $G_L=Sp_L\ts$, there is a simplified formula 
for the operator $F_{\La^r}$, similar to (\ref{3.11}). We have
\begin{equation}\label{3.12}
F_{\La^r}=\,
\prod_{(i,\ts j)}^{\longrightarrow}\, 
\left(1-\frac{Q_{ij}}
{\ts c_i(\La^r)+c_j(\La^r)+L+1}\ts\right)
\cdot E_{\La^c}
\ \ \text{for}\ \ 
G_L=Sp_L
\end{equation}

\vskip-4pt
\noindent
where the ordered product is taken over all pairs $(i\com j)$
such that $i$ and $j$ appear in different rows of $\La^r$.
The proof of the formula (\ref{3.12})
is similar to that of (\ref{3.11}). Here we work with the rows of 
the tableau $\La^r$. But here from the very beginning
of the proof we can set $t_1=\ldots=t_l=\frac12$ in the first product
over $1\le i<j\le l$ in (\ref{3.1}) where $\La=\La^r$,
at the same time replacing the second product in (\ref{3.1}) by $E_{\La^r}$.
If $k$ and $k+1$ appear in the same row of
$\La^r$, due to (\ref{2.3}) the operator $E_{\La^r}$ is divisible on the
left by
$$
1+P_{k,k+1}=
R_{\ts k,k+1}(\ts c_k+\textstyle\frac12\ts\,\com\ts c_{k+1}+\frac12\ts)\,.
$$
We omit other details of the proof of the formula (\ref{3.12}).


\medskip\noindent\textbf{3.3.}
For any standard tableau $\La\ts$, denote by
$V_\La$ and $W_\La$ the images in the space $\CLl$ of the operators
$E_\La$ and $F_\La\ts$, respectively. We have $W_\La\subset V_\La\ts$
by the second statement of Proposition 3.2. Consider the subspace
(\ref{3.0}) of the traceless tensors. The next proposition is pivotal
for the present article.

\begin{Proposition}
We have the equality of vector spaces\/ $W_\La=V_\La\cap\CLl_{\,\ts0}$\,.
\end{Proposition}

\begin{proof}
For any distinct indices $i\com j\in\{1\lcd l\}$ the operator $P_{ij}$
preserves the subspace (\ref{3.0}) while the operator $Q_{ij}$
vanishes on this subspace. So the action of (\ref{3.1}) on 
this subspace coincides with the action of the second product 
over $1\le i<j\le l$ in (\ref{3.1}). Hence by the definition of the
operator $F_\La$ we get the equality $F_\La\cdot w=E_\La\cdot w$
for any vector $w\in\CLl_{\,\ts0}\,$. Therefore
$$
W_\La\ts\supset\ts E_\La\cdot\CLl_{\,\ts0}=V_\La\cap\CLl_{\,\ts0}\,.
$$
We already noted that $W_\La\subset V_\La\,$. It remains to prove
that $W_\La\subset\CLl_{\,\ts0}$.
Equivalently, we have to prove that $Q_{ij}\,F_\La=0$
for any distinct $i$ and $j\ts$.

Firstly let us prove the equality $Q_{12}\,F_\La=0$ for every $\La$.
Consider the case $G_L=O_L$.
In every standard tableau $\La$, the
number $1$ always occupies the upper left corner. The number $2$
appears in $\La$ either next to the right of $1$, or next down from $1$.
In the latter case, the operator $E_\La$ is divisible on the left
by $1-P_{12}\ts$. By Proposition~3.2,
the operator $F_\La$ 
is then also divisible on the left
by $1-P_{12}\ts$.
But $Q_{12}(1-P_{12})=0$ for $G_L=O_L$ by (\ref{3.9}). 

Suppose that $2$ appears in $\La$ next to the right of $1$. 
Then for any $j\neq1$
$$
c_1+c_j\ge1-\la'_1\ge1+\la'_2-L\ge2-L\,.
$$
So none of the factors $\Rb_{\ts 1j}(\ts c_1+t_1\com c_j+t_j\ts)$
in the product (\ref{3.1}) has a pole at $t_1=\ldots=t_l=-\ts\frac12\,$.
Some of the factors $\Rb_{\ts ij}(\ts c_i+t_i\com c_j+t_j\ts)$ 
with $i\neq1$ may have poles at $t_1=\ldots=t_l=-\ts\frac12\,$.
Then consider the standard tableau
of shape $\la\ts$,
which is obtained by placing the numbers $1$ and $2$ in the first row
as in $\La$, and by filling the remaining boxes of the Young diagram
$\la$ by columns. Let us denote this new tableau by $\Lan$.
Arguing as in the proof of Proposition 3.2,
we show that the equality $Q_{12}\,F_\La=0$ is sufficient to
prove for $\La=\Lan$. This argument involves multiplying the product
(\ref{3.1}) on the left by (\ref{3.35}) and dividing on the right
by (\ref{3.3535}),
for some $k\in\{3\lcd l-1\}\ts$ such that the tableau $s_k\La$ is standard.
Further, arguing as in the proof of Proposition~3.2, we show that 
$F_{\Lan}$ is divisible on the left by
$$
{\textstyle
\Rb_{12}(\ts c_1-\frac12\,\com\ts c_2-\frac12\ts)
\displaystyle
\ts=\ts1-Q_{12}\ts/\ts L\,\ts.}
$$
The equality $Q_{12}\,F_{\Lan}=0$ now follows from
(\ref{3.99}).

In the case $G_L=Sp_L$, the proof of the equality
$Q_{12}\,F_\La=0$ for any standard tableau $\La$ is much simpler.
In this case for any $\La$, every factor
\hbox{\ts}$\Rb_{\ts ij}(\ts c_i+t_i\com c_j+t_j\ts)$
in the product (\ref{3.1}) is regular at $t_1=\ldots=t_l=\frac12\,$.
If the numbers $1$ and $2$ appear in the first row of $\La$, then
the operator $E_\La$ is divisible on the left
by $1+P_{12}\ts$. The operator $F_\La$ is then also
divisible on the left by $1+P_{12}\ts$. But
for $G_L=Sp_L$ we have $Q_{12}(1+P_{12})=0$.
If $1$ and $2$ appear in the first column of $\La$, then
the operator $F_\La$ is divisible on the left by

\vskip-16pt
$$
\textstyle
\Rb_{12}(\ts c_1+\frac12\,\com\ts c_2+\frac12\ts)
\displaystyle
\ts=\ts1-Q_{12}\ts/\ts L\,\ts.
$$

\vskip4pt\noindent
The equality $Q_{12}\,F_{\La}=0$ now follows from (\ref{3.99}).

Let us now prove the equality $Q_{ij}\,F_\La=0$ for any pair
of distinct indices $i$ and $j$. If for some $k\in\{1\lcd l-1\}\ts$
the tableau $s_k\La$ is standard, then by Corollary 3.2
$$
P_{\ts k,k+1}\ts F_\La\,=\,
\frac{F_\La}{c_{k+1}-c_k}\,+\,
F_{s_k\La}\,R_{\ts k,k+1}\ts(\ts c_k\com c_{k+1})\,P_{\ts k,k+1}\,.
$$
If the tableau $s_k\La$ is not standard,
then $P_{k,k+1}\ts F_\La$ equals $F_\La$
or $-F_\La\,$; this equality follows from
the second statement of Proposition~3.2 by (\ref{2.3}).

For any permutation $s\in S_l\ts$, let $P_s$ be
the corresponding operator on the space $\CLl$.
For some elements $R_{\La^\prime}(s)\in\Cel$ that may also 
depend on $\La\ts$, we have
$$
P_s\,F_\La=\ts
\sum_{\La^\prime}\,\ts F_{\La^\prime}\ts R_{\La^\prime}(s)\,,
$$
where $\La^\prime$ ranges over all standard tableaux of shape $\la\,$.
Let us now choose the permutation $s$ so that $s(i)=1$ and $s(j)=2$.
Then we get
$$
Q_{ij}\,F_\La=P_s^{-1}Q_{12}\,P_s\ts F_\La=
\,\sum_{\La^\prime}\,P_s^{-1}\ts Q_{12}\ts F_{\La^\prime}\ts
R_{\La^\prime}(s)=0.\quad\qed
$$
\end{proof}

By definition, the associative algebra $\Cel$ is semisimple.
Denote by $\Cela$ the simple ideal of $\Cel$
corresponding to the irreducible $\Cel\ts$-module~$U_\la$.
Proposition 3.3 gives the following characterization of the
$F_\La\ts$.

\begin{Corollary}
The operator $F_\La\in\Cel$ is the unique element of the simple ideal
$\Cela\ts$,
with the image under\/ {\rm(\ref{3.22})} corresponding to\/
$e_\La\in\CC S_l\ts$.
\end{Corollary}

\begin{proof}
By definition, we have $F_\La\in E_\La+\Idl\ts$;
see Proposition~2.2. Hence the image of $F_\La$ under the 
homomorphism {\rm(\ref{3.22})} corresponds to 
the element $e_\La\in\CC S_l\ts$.
By Proposition 3.3, we also have $F_\La\in\Cela$, see the beginning
of Subsection 3.2. But any element of the ideal $\Cela\subset\Cel$
is uniquely determined by the image of this element under (\ref{3.22}).
\qed
\end{proof}

Note that due to (\ref{2.01}), the above characterization of
the operator $F_\La$ implies the equality

\vskip-14pt
\begin{equation}\label{3.301}
F_\La^{\ts2}=F_\La\cdot l\ts!\ts/\dim U_\la\ts.
\end{equation}


\smallskip\noindent\textbf{3.4.}
We will now extend the results of Subsection 3.2 to standard tableaux
of skew shapes. Take any $m\in\{0\lcd l-1\}$\ts. As in
Subsection~2.3, denote by $\Up$ the standard tableau obtained from $\La$
by removing the boxes with the numbers $m+1\lcd l$. Let $\mu$ be the shape of
the tableau $\Up$. Define a standard tableau $\Om$ of the skew shape
$\lm$ by setting $\Om(i\com j)=\La(i\com\ns j)-m$ for all
$(i\com j)\in\lm$\ts. Put $n=l-m$, the number of elements in the set $\lm$.

Take any $M\in\{0\lcd L-1\}$ and put $N=L-M$. 
Then choose the decomposition $\CC^L=\CC^N\oplus\CC^M$ so that the
subspaces $\CC^N$ and $\CC^M$ in $\CC^L$ are orthogonal 
relative to the form $\langle\ ,\,\rangle$ on $\CC^L$.
In the case $G_L=Sp_L$ both $N$ and $M$ have to be even.
If $M>0$, the restriction of the form $\langle\ ,\,\rangle$ from $\CC^L$
to $\CC^M$ is non-degenerate. Consider the corresponding
subspace of traceless tensors

\vskip-16pt
\begin{equation}\label{3.000}
\CMm_{\,\ts0}\subset\CMm.
\end{equation}
We will assume that the irreducible representation $V_\mu$
of the group $G_M$ appears in the subspace (\ref{3.000}),
so that the partition $\mu$ of $m$ satisfies the conditions from
\cite{W} for $G_M$ as described in Subsection 1.3. In particular,
if $M=1\ts$, then $G_M=O_{\ts1}\ts$ and we have to take $m\in\{0\com1\}$.

Recall that in Subsection 2.4, we denoted by $I_m$ the projector to the 
direct summand (\ref{4.7}). Now, there is a unique $G_M\ts$-invariant
projector  
$$
H_m:\,\CMm\,\to\,\CMm_{\,\ts0}\,.
$$
Let us denote by $J_m$ the composition of the operators $I_m$ and 
$H_m\ot1\ts$. This composition is a projector to the subspace
\begin{equation}\label{3.41}
\CMm_{\,\ts0}\ot\CNn\subset\CLl.
\end{equation}
Regard $J_m$ as an operator on the vector space $\CLl$. Then,
for any linear operator $A$ on $\CLl$,
define an operator on the same space,
$$
A^\vee=\,J_m\,A\,J_m\,.
$$
We may also regard $A^\vee$ as an operator on the subspace (\ref{3.41}).
The projector $J_m$ commutes with the action of subgroup
$G_N\times G_M\subset G_L$ on $\CLl$. So for $A\in\Cel\ts$,

\vskip-16pt
$$
A^\vee\in(\Cem/\,\Idm)\otimes\Cen\,.
$$

\vskip2pt\noindent
Here the quotient algebra $\Cem/\,\Idm$ is identified with 
the image of the symmetric group ring $\CC S_m$ in $\End(\CMm_{\,\ts0})$\ts.
We get a linear map
\begin{equation}\label{3.46}
\Th_m:\Cel\rightarrow(\Cem/\,\Idm)\otimes\Cen:A\mapsto A^\vee\,.
\end{equation}

As the map $\Th_m$ is linear, for any distinct $i\com j\in\{1\lcd l\}$
we have
$$
\Rp_{ij}(x\com y)=1-\frac{\Pp_{ij}}{x-y}
\qquad\textrm{and}\qquad
\Rbp_{ij}(x\com y)=1-\frac{\Qp_{ij}}{x+y+L}\,.
$$
If $i\le m<j$ or $i>m\ge j$, then we have $\Pp_{ij}=0$ and $\Qp_{ij}=0$.
If $i\com j\le m$ then $\Qp_{ij}=0$ also, while
the operator $\Pp_{ij}$ on (\ref{3.41})
acts in the tensor factor $\CMm_{\,\ts0}$ as the permutation corresponding to
$(i\ts j)\in S_m\,$, and it also acts as the identity in
$\CNn$. If $i\com j>m$ then each of the operators 
$\Pp_{ij}$ and $\Qp_{ij}$ acts as the identity on $\CMm_{\,\ts0}$. 
Then the action of $\Pp_{ij}$ on $\CNn$ corresponds
to the transposition $(i-m\com j-m)\in S_n\ts$. Then $\Qp_{ij}$ acts
as $Q(N)$ in the $(i-m)$-th and $(j-m)$-th tensor factors of $\CNn$,
and~acts as the identity in the remaining $n-2$ tensor factors of $\CNn$.
Our extension~of the results of Subsection 3.2 to standard tableaux
of skew shapes is based on the next observation, see Lemma~2.3.

\begin{Proposition}
The image of the product\/
{\rm(\ref{3.1})} under the map\/ $\Th_m$ equals
\begin{equation}\label{3.42}
\prod_{1\le i<j\le m}^{\longrightarrow}
\Rp_{\ts ij}(\ts c_i+t_i\com c_j+t_j\ts)
\end{equation}

\vskip-10pt
\begin{equation}\label{3.43}
\times\,
\prod_{m<i<j\le l}^{\longrightarrow}
\Rbp_{\ts ij}(\ts c_i+t_i\com c_j+t_j\ts)
\ \cdot\!
\prod_{m<i<j\le l}^{\longrightarrow}
\Rp_{\ts ij}(\ts c_i+t_i\com c_j+t_j\ts)\,.
\end{equation}
\end{Proposition}

\begin{proof}
When multiplying the product (\ref{3.1}) on the left by the 
operator $J_m\ts$, the factors 
$\Rb_{\ts ij}(\ts c_i+t_i\com c_j+t_j\ts)$ with
$i<j\le m$ and $i\le m<j$ cancel. Indeed,
since $\CMm_{\,\ts0}\subset(\CC^L)^{\ot m}_{\,\ts0}$,
we have $J_m\,Q_{ij}=0$ for $i<j\le m$. Since 
the subspaces $\CC^M$ and $\CC^N$
are orthogonal with respect to the form $\langle\ ,\,\rangle$ on
$\CC^L$, we also have $J_m\,Q_{ij}=0$ for $i\le m<j$. Thus by multiplying
(\ref{3.1}) by $J_m$ on the left and right, we obtain the product
\begin{equation}\label{3.44}
J_m\,\ts\cdot\hskip-8pt
\prod_{m<i<j\le l}^{\longrightarrow}\hskip-6pt
\Rb_{\ts ij}(\ts c_i+t_i\com c_j+t_j\ts)
\ \cdot\hskip-6pt
\prod_{1\le i<j\le l}^{\longrightarrow}\hskip-6pt
R_{\ts ij}(\ts c_i+t_i\com c_j+t_j\ts)
\ts\,\cdot\ts\,J_m\,.
\end{equation}

Expand the product of all factors 
$R_{\ts ij}(\ts c_i+t_i\com c_j+t_j\ts)$ in (\ref{3.44}) 
as a sum of the operators $P_s$ on $\CLl$
corresponding to permutations $s\in S_l\ts$,
with the coefficients from $\CC(t_1\lcd t_l)$\ts.
The operators $\Rb_{\ts ij}(\ts c_i+t_i\com c_j+t_j\ts)$ with $m<i<j$ 
commute with the operator on $\CLl$, which acts as the projector to
the subspace $\CC^M\subset\CC^L$ along $\CC^N$
in each of the first $m$ tensor
factors, and which acts trivially in the last $n$ tensor factors of $\CLl$.
Hence the summand in (\ref{3.44}) corresponding to $P_s$
vanishes, unless the permutation $s\in S_l$
preserves the subset $\{1\lcd m\}\subset\{1\lcd l\}$.
By Lemma 2.3, the product (\ref{3.44}) then equals
$$
J_m\,\,\cdot\hskip-6pt
\prod_{m<i<j\le l}^{\longrightarrow}
\Rb_{\ts ij}(\ts c_i+t_i\com c_j+t_j\ts)
$$
$$
\times\ 
\prod_{1\le i<j\le m}^{\longrightarrow}
R_{\ts ij}(\ts c_i+t_i\com c_j+t_j\ts)
\ \cdot\!
\prod_{m<i<j\le l}^{\longrightarrow}
R_{\ts ij}(\ts c_i+t_i\com c_j+t_j\ts)
\ts\,\cdot\ts\,J_m
$$
\begin{equation}\label{3.33}
=\ \,
J_m\,\ts\cdot\hskip-6pt
\prod_{m<i<j\le l}^{\longrightarrow}
\Rb_{\ts ij}(\ts c_i+t_i\com c_j+t_j\ts)
\ \cdot\ts\,J_m
\end{equation}
$$
\times\ 
\prod_{1\le i<j\le m}^{\longrightarrow}
\Rp_{\ts ij}(\ts c_i+t_i\com c_j+t_j\ts)
\ \,\cdot\hskip-6pt
\prod_{m<i<j\le l}^{\longrightarrow}
\Rp_{\ts ij}(\ts c_i+t_i\com c_j+t_j\ts)\,.
$$
\smallskip\noindent

Now expand the expression displayed in the line (\ref{3.33}), 
as a sum over all subsequences
\begin{equation}\label{3.47}
(i_1\com j_1)\lcd(i_d\com j_d)
\end{equation}

\vskip4pt\noindent
in the sequence of lexicographically ordered pairs $(i\com j)$
with $m<i<j\le l$. The summand of (\ref{3.33})
corresponding to the subsequence (\ref{3.47}) 
is equal, up to a coefficient from $\CC(t_1\lcd t_l)$\ts, to the product
$
J_m\,Q_{\ts i_1j_1}\ldots Q_{\ts i_dj_d}\,J_m\,.
$
We will prove that for some permutation $s$ of the set
$\{i_1\com j_1\}\cup\ldots\cup\{i_d\com j_d\}\ts$, and for some
pairwise distinct elements $\bi_1\com\bj_1\lcd\bi_c\com\bj_c$
of this set, we have
\begin{equation}\label{3.48}
Q_{\ts i_1j_1}\ldots Q_{\ts i_dj_d}=
P_s\,Q_{\ts\bi_1\bj_1}\ldots Q_{\ts\bi_c\bj_c}\ts.
\end{equation}
Since
$$
J_m\,P_s\,Q_{\ts\bi_1\bj_1}\ldots Q_{\ts\bi_c\bj_c}\,J_m=
\Pp_s\,\Qp_{\ts\bi_1\bj_1}\ldots\Qp_{\ts\bi_c\bj_c}\ts,
$$
we will then also have the equality
$$
J_m\,Q_{\ts i_1j_1}\ldots Q_{\ts i_cj_c}\,J_m=
\Qp_{\ts i_1j_1}\ldots\Qp_{\ts i_cj_c}\ts.
$$
The last equality shows that the expression (\ref{3.33})
equals the first product over $m<i<j\le l$ in (\ref{3.43}).
Hence Proposition 3.4 will follow from (\ref{3.48}).

Since $Q_{ij}=Q_{ji}$ for any $i\neq j$, we can assume that
$\bi_1<\bj_1\lcd\bi_c<\bj_c$ in (\ref{3.48}). We will prove
the equality (\ref{3.48}) under this assumption. We
will also show that the indices $\bi_1\com\bj_1\lcd\bi_c\com\bj_c$
in (\ref{3.48}) can be so chosen that
\begin{equation}\label{3.49}
c\le d\,,\ \,\textrm{and}\ \,
(\bi_b\com\bj_b)\preccurlyeq(i_d\com j_d)
\ \,\textrm{for all}\ \,
b=1\lcd c\,.
\end{equation}
Here $\prec$ indicates the lexicographical ordering.
The operators $Q_{\ts\bi_b\bj_b}$ in (\ref{3.48})
pairwise commute, so their ordering is irrelevant.
We proceed by induction on $d=0\com1\com2\com\,\ldots\,\,$.
If $d=0$, there is nothing to prove. Now suppose that for
some subsequence (\ref{3.47}), the equality (\ref{3.48}) is true
along with (\ref{3.49}).

Suppose that $i<j$ and $(i_d\com j_d)\prec(i,j)\,$.
By the induction assumption, we have

\vskip-16pt
$$
Q_{\ts i_1j_1}\ldots Q_{\ts i_dj_d}\,Q_{ij}=
P_s\,Q_{\ts\bi_1\bj_1}\ldots Q_{\ts\bi_c\bj_c}\ts Q_{ij}\,.
$$

\vskip4pt\noindent
This equality makes the induction step, unless $i$ or $j$
coincides with one of the indices $\bi_1\com\bj_1\lcd\bi_c\com\bj_c\ts$.
Suppose there is such a coincidence.
Let $b$ be the maximal number such that
$\{\bi_b\com\bj_b\}\cap\ts\{i\com j\}\neq\varnothing\,$.
Here $\bi_b\neq j$ because $(\bi_b\com\bj_b)\prec(i\com j)\ts$.
Then there are three possibilities:
\begin{equation}\label{3.333}
i=\bi_b<\bj_b<j\,;
\qquad
\bi_b<\bj_b=i<j\,;
\qquad 
\bi_b<i<j=\bj_b\,.
\end{equation}

Consider, for instance, the first of the possibilities (\ref{3.333}).
Here we have
$$
Q_{\ts\bi_b\bj_b}\ts Q_{ij}=
Q_{\ts\bi_b\bj_b}\ts Q_{\ts\bi_b j}=
P_{\ts\bj_b j}\ts Q_{\ts\bi_b j}\,.
$$
For each $a=1\lcd c$ we have $\bi_a\neq j$. By our choice of $b\ts$,
we have $\bj_a\neq j$ for $a=b\lcd c$.
If $\bj_a\neq j$ also for $a=1\lcd b-1$, then we have the equality
$$
P_s\,Q_{\ts\bi_1\bj_1}\ldots Q_{\ts\bi_c\bj_c}\ts Q_{ij}=
P_s\,P_{\ts\bj_b j}\,
Q_{\ts\bi_1\bj_1}\ldots Q_{\ts\bi_{b-1}\bj_{b-1}}
Q_{\ts\bi_{b+1}\bj_{b+1}}\ldots Q_{\ts\bi_c\bj_c}\ts Q_{\ts\bi_b j}\,,
$$
which makes the induction step. If $\bj_a=j$ for some number $a<b$,
then the number $a$ is unique, and we have the equalities
$$
P_s\,Q_{\ts\bi_1\bj_1}\ldots Q_{\ts\bi_c\bj_c}\ts Q_{ij}=
P_s\,Q_{\ts\bi_1\bj_1}\ldots Q_{\ts\bi_{a-1}\bj_{a-1}}P_{\ts\bj_b\bj_a}
Q_{\ts\bi_{b+1}\bj_{b+1}}\ldots Q_{\ts\bi_c\bj_c}\ts Q_{\ts\bi_b\bj_a}=
$$
$$
P_s\ts P_{\ts\bj_b\bj_a}\ns
Q_{\ts\bi_1\bj_1}\ns\ldots Q_{\ts\bi_{a-1}\bj_{a-1}}\ns Q_{\ts\bi_a\bj_b}
Q_{\ts\bi_{a+1}\bj_{a+1}}\ns\ldots Q_{\ts\bi_{b-1}\bj_{b-1}}\ns
Q_{\ts\bi_{b+1}\bj_{b+1}}\ns\ldots Q_{\ts\bi_c\bj_c}Q_{\ts\bi_b\bj_a}.
$$
Thus we again make the induction step, because
$(\bi_a\com\bj_a)\prec(i\com j)=(\bi_b\com\bj_a)$
implies $\bi_a<\bi_b\,$, so that here we have $\bi_a<\bj_b$ and
$(\bi_a\com\bj_b)\prec(\bi_b\com\bj_a)$.

The second and the third possibilities in (\ref{3.333}) are treated
similarly. We omit the details of this treatment.\qed
\end{proof}

Consider the product (\ref{3.43}) as a rational function
of the variables $t_1\lcd t_l\ts$. This function may actually depend
only on $t_{m+1}\lcd t_l\ts$. The values of this function are linear
operators on the subspace (\ref{3.41}), which act as the identity in
the tensor factor $\CMm_{\,\ts0}$ of (\ref{3.41}).
Regard $E_\Om$ as the operator acting
on the tensor factor $\CNn$ of (\ref{3.41}).

\begin{Corollary}
Restriction of\/ {\rm(\ref{3.43})} to $\T_\La$ 
is regular at $t_{m+1}=\ldots=t_l=\mp\ts\frac12$.
The value of this restriction at
$t_{m+1}=\ldots=t_l=\mp\,\frac12$ is divisible on the left and right
by the operator $\id\ot E_\Om\ts$.
\end{Corollary}

\begin{proof}
Consider the restriction of the rational function (\ref{3.43}) to
$\T_\La$. Due to Propositions 3.2 and 3.4, this
restriction is regular at $t_1=\ldots=t_l\ts$; the value of this
restriction at $t_1=\ldots=t_l$ is $\Th_m(F_\La)\ts$. By Proposition 2.2
applied to the tableau $\Up$ instead of to $\La$, the value
at $t_1=\ldots=t_m$ of the restriction of the function (\ref{3.42})
to $\T_\La$ equals the operator
$
\,(\ts E_{\ts\Up}\,|\,\CMm_{\,\ts0}\ts)\ot\id\,;
$
here the first tensor factor is the restriction
of the operator $E_{\ts\Up}$ in $\CMm$ to the subspace (\ref{3.000}).
Hence the restriction to $\T_\La$ of the
product (\ref{3.43}) is regular at $t_{m+1}=\ldots=t_l$.
The second statement of Corollary 3.4 now follows from Proposition 2.3.
\qed
\end{proof}

In the notation of Subsection 1.2, $c_k(\Om)=c_{m+k}$ for each
$k=1\lcd n$. When $(t_1\lcd t_l)\in\T_\La\ts$, we can put
$t_k(\Om)=t_{m+k}$ for $k=1\lcd n$\ts. Then we get Theorem 1.4 
as a reformulation of Corollary 3.4. Moreover, the value at
$t_{m+1}=\ldots=t_l=\mp\ts\frac12$ of the 
restriction of\/ {\rm(\ref{3.43})} to $\T_\La$, is $\id\ot\Fom$. 

If $M=0$, then $\mu=(0\com0\ts,\ts\ldots\ts)$ and $\Om=\La$,
so that $\Fom=F_\La\ts$. For arbitrary $M$ and $\Om$,
our proof of Corollary~3.4 shows that the operator $\Fom$ on $\CNn$,
as defined by Theorem 1.4, satisfies the relation
\begin{equation}\label{3.4444}
\Th_m(F_\La)=(\ts E_{\ts\Up}\,|\,\CMm_{\,\ts0}\ts)\ot\Fom\,.
\end{equation}
This relation, along with Proposition 3.2,
can be regarded as an alternative definition of the operator $\Fom$.
Using the relation (\ref{3.4444}) along with the proof of
Proposition 3.4, we also obtain the simplified formulas
(\ref{1.6}) and (\ref{1.7}) from the formulas
(\ref{3.11}) and (\ref{3.12}), respectively.


\medskip\noindent\textbf{3.5.}
The image of the operator $\Fom$ acting on the vector space
$\CNn$ is denoted by $\Wom$. If $M=0$, then 
$\mu=(0\com0\ts,\ts\ldots\ts)$ and $\Om=\La\ts$.
In this case $W_\Om\ts(0)=W_\La\ts$, so that
we have the equality (\ref{1.4444}) by Proposition 3.3.
For arbitrary $M\ge0$, we will need the following
analogue of Proposition~2.4. 

\begin{Proposition}
If\/ $\Wlm\neq\{0\}$, then $\Wom\neq\{0\}$.
\end{Proposition}

\begin{proof}
Let us realize the irreducible representation $W_\la$ of the group $G_{N+M}$
from (\ref{1.4}) as the image $W_\La\subset(\CC^{N+M})^{\ts\ot\,l}$
of the operator $F_\La\,$. Define the vectors
$w\ts(M)\in(\CC^M)^{\ts\ot\ts2}$ and
$w\ts(N+M)\in(\CC^{N+M})^{\ts\ot\ts2}$
as in (\ref{1.44444444}). Note~that
\begin{equation}\label{3.51}
w\ts(N+M)=w\ts(N)+w\ts(M)\,.
\end{equation}
The vector space $(\CC^{N+M})^{\ot\,l}$
is the sum of its subspaces, obtained by some permutations
of the $l$ tensor factors from the subspaces of the form
\begin{equation}\label{3.52}
(\CC^M)^{\ts\ot\ts k}_{\,\ts0}\ot
(\CC\ts w\ts (M))^{\ts\ot\ts d}\ot
(\CC^N)^{\ts\ot\ts(l-k-2d)}\,;
\end{equation}
here $k$ and $d$ are non-negative integers such that $k+2d\le l$.
This sum may be not direct, but every summand is preserved by 
the action of the subgroup $G_M\subset G_{N+M}\ts$;
see \cite[Section V.6]{W}. For instance, consider
the subspace (\ref{3.52}) itself. The image of this subspace
under the operator
$F_\La$ coincides with the image of the subspace
\begin{equation}\label{3.53}
(\CC^M)^{\ts\ot\ts k}_{\,\ts0}\ot
(\CC\ts w\ts (N))^{\ts\ot\ts d}\ot
(\CC^N)^{\ts\ot\ts(l-k-2d)}\ts.
\end{equation}
This coincidence follows from (\ref{3.51}) and from the equalities 
$F_\La\,Q_{ij}=0$ for
\begin{equation}\label{3.54}
(i\com j)=(k+1\com k+2)\lcd(k+2d-1\com k+2d)\,;
\end{equation}
here $Q_{ij}$ are operators on $(\CC^{N+M})^{\ts\ot\,l}$.
These operator equalities are implied
by Proposition 3.3\ts; see also the beginning
of the proof of Proposition 3.2.
But the subspace (\ref{3.53}) is contained in
$(\CC^M)^{\ts\ot\ts k}_{\,\ts0}\ot(\CC^N)^{\ts\ot\ts(l-k)}\ts$. Note that
$$
\quad
{\rm Hom}_{\,G_M}(\ts W_\mu\ts\com F_\La\cdot
(\CC^M)^{\ts\ot\ts k}_{\,\ts0}\ot(\CC^N)^{\ts\ot\ts(l-k)})\neq\{0\}
\ \,\Rightarrow\ \,k=m\ts.
$$

Further, define a projector $J_k^{\ts(d)}$ from the vector space 
$(\CC^{N+M})^{\ts\ot\,l}$ to the subspace (\ref{3.52}) as follows.
Let $H_k^{\ts(d)}$ be the linear operator on the tensor product
$(\CC^M)^{\ts\ot\,(k+2d)}\ts$, acting as the (unique) $G_M\ts$-invariant
projector 
$$
H_k:\,
(\CC^M)^{\ts\ot\ts k}
\,\to\,
(\CC^M)^{\ts\ot\ts k}_{\,\ts0}
$$
in the first $k$ tensor factors of $(\CC^M)^{\ts\ot\,(k+2d)}\ts$,
and acting as the operator $(Q(M)/M)^{\ts\ot\ts d}$ in 
the last $2d$ tensor factors of  $(\CC^M)^{\ts\ot\,(k+2d)}\ts$.
Then $J_k^{\ts(d)}$ is the composition of the projector $I_{k+2d}$
to the direct summand 
$$
(\CC^M)^{\ts\ot\ts (k+2d)}\ot
(\CC^N)^{\ts\ot\ts(l-k-2d)}
\subset
(\CC^{N+M})^{\ts\ot\,l}\,,
$$
with the operator $H_k^{\ts(d)}\ot1\ts$.
Note that the projector $J_k^{\ts(d)}$ is $G_M\ts$-equivariant.

In a similar way, by using the operator $(Q(N)/N)^{\ts\ot\ts d}$
instead of the operator $(Q(M)/M)^{\ts\ot\ts d}$, define a projector
$\bar{J}_k^{\ts(d)}$ from $(\CC^{N+M})^{\ts\ot\,l}$ 
to the subspace (\ref{3.53}).
In the notation of Subsection 3.4, we have
$$
J_k^{\ts (0)}=\bar{J}_k^{\ts (0)}=J_k\,.
$$

Let us apply the projector $J_k^{\ts (d)}$ to the subspace
$W_\La\subset(\CC^{N+M})^{\ts\ot\,l}\ts$.
If an irreducible representation
of $G_M$ equivalent to $W_\mu$ occurs in the projection
of $W_\La$ to (\ref{3.52}),
it also occurs in the projection of $W_\La$ to (\ref{3.53}).
This follows from (\ref{3.51}) and the equalities 
$Q_{ij}\ts F_\La=0$ for the pairs of indices $i$ and $j$ given by 
(\ref{3.54}). Either occurrence implies that $k=m$. But again we have
$$
\bar{J}_m^{\ts (d)}\cdot F_\La
\,\subset\, 
\bar{J}_m^{\ts (0)}\cdot F_\La
\,=\,
J_m\cdot F_\La\,.
$$

Now suppose that $\Wlm\neq\{0\}$. 
By the above argument,
there exist permutations $s^{\ts\prime}$ and $s^{\ts\prime\prime}$ 
in $S_l$ such that
\begin{equation}\label{3.56}
{\rm Hom}_{\,G_M}(\, W_\mu\,\com\ts 
J_m\,P_{s'}\,F_\La\,P_{s''}\cdot\ts
\CMm_{\,\ts0}\ot\CNn\ts)\neq\{0\}\,.
\end{equation}
Let $\Lap$ and $\Lapp$ range over the set of
all standard tableaux of shape $\la\ts$. 
Denote by $F_{\La'\La''}$ the unique element of the simple ideal
$\Cela\subset\Cel\ts$, with the image under {\rm(\ref{3.22})} 
corresponding to the matrix element $e_{\La'\La''}\in\CC S_l\ts$.
If $\La=\Lap=\Lapp$, then $F_\La=F_{\La'\La''}$ by Corollary 3.3.
The product $P_{s'}\,F_\La\,P_{s''}$ in (\ref{3.56}) can be 
written as a linear combination of the operators $F_{\La'\La''}$
on $(\CC^{N+M})^{\ot\,l}$, with the coefficients from $\CC\ts$.
Due to (\ref{3.56}), there exists at least
one pair of tableaux $\Lap$ and $\Lapp$ such that
\begin{equation}\label{3.57}
{\rm Hom}_{\,G_M}(\, W_\mu\,\com\ts 
J_m\,F_{\La'\La''}\cdot\ts
\CMm_{\,\ts0}\ot\CNn\ts)\neq\{0\}\,.
\end{equation}

According to Subsection 3.4,
the restriction of the operator $J_m\,F_{\La'\La''}$
to the subspace (\ref{3.41}) is denoted by $\Th_m(F_{\La'\La''})\,$.
Consider the tableaux $\Up^{\ts\prime}$ and $\Up^{\ts\prime\prime}$,
obtained by removing the boxes with the numbers $m+1\lcd l$ from the 
tableaux $\Lap$ and $\Lapp$, respectively.
The operator $\Th_m(F_{\La'\La''})$ on 
the subspace (\ref{3.41}) is divisible on the left and right,
respectively by
$$
(\ts E_{\ts\Up'}\,|\,\CMm_{\,\ts0}\ts)\ot\id
\ \quad\text{and}\ \quad
(\ts E_{\ts\Up''}\,|\,\CMm_{\,\ts0}\ts)\ot\id\,.
$$
These divisibility properties follow from Proposition 3.2\ts;
see the proof of Lemma 2.4. Now the inequality (\ref{3.57})
implies that the tableaux
$\Up^{\ts\prime}$ and $\Up^{\ts\prime\prime}$ are of the same shape.
But then $F_{\La'\La''}=F_{\La'}\ts E$ for the
operator $E$ on the subspace (\ref{3.41}) corresponding to
some invertible element
$e\in\CC S_{mn}\ts$; see the end of the proof of Proposition 2.4. Then
$$
F_{\La'\La''}\cdot\ts\CMm_{\,\ts0}\ot\CNn\ts=\,
F_{\La'}\cdot\ts\CMm_{\,\ts0}\ot\CNn\,.
$$
By applying the relation (\ref{3.4444}) to the tableau $\Lap$ instead
of $\La\ts$, the inequality (\ref{3.57})
now implies that the tableau $\Up^{\ts\prime}$ is of shape $\mu\ts$.
Moreover, the left-hand side of (\ref{3.57})
equals $W_{\Om'}(M)$ for some standard tableau $\Om^{\ts\prime}$
of skew shape $\lm\ts$. The
inequality $W_{\Om'}(M)\neq\{0\}$ implies that $\Wom\neq\{0\}$.
\qed
\end{proof}

The space $\Wlm$
comes with an action of the subgroup $G_N\subset G_L\ts$.
The subspace $\Wom\subset\CNn$ is preserved by the action of 
the group $G_N$
because $\Fom\in\Cen\ts$. Let us consider $\Wlm$ and $\Wom$ as
representations of the group $G_N$. In Subsection 5.6
we will prove that
these representations are equivalent, as stated in Proposition~1.4.


\medskip\noindent\textbf{3.6.}
Let $\g_L$ be the Lie algebra of $G_L$,
so that $\g_L=\so_L$ or $\g_L=\sp_L$.
We regard $\g_L$ as a Lie subalgebra in $\glL$. 
Let us now consider representations of the group $G_L$ as $\g_L$-modules.
If $G_L=Sp_L$, the $\sp_L$-modules $W_\la$ for different
partitions $\la$ of $l=0\com1\com2\com\ts\ldots$ with $2\la'_1\le L$
are irreducible and pairwise non-equivalent. 
One obtains all irreducible finite-dimensional $\sp_L\ts$-modules in this way. 
If we choose the Borel subalgebra, and then fix the basis
in the corresponding Cartan subalgebra of $\sp_L$ as in
\cite[Subsection 1.1]{KT}, then $W_\la$ is
the irreducible $\sp_L\ts$-module of highest weight
$(\la_1\lcd\la_{\ts L/2})$. However, in the present article 
we do not use the highest weight theory of $\sp_L\ts$-modules.

If $G_L=O_L$, we have $\la'_1+\la'_2\le L$. Therefore,
by changing
the length $\la'_1$ of the first column of
the Young diagram of the partition $\la$
to $L-\la'_1$,
we obtain the Young diagram of a certain partition.
Denote this partition by $\la^\ast$. The two representations
$W_\la$ and $W_{\la^\ast}$ are called \textit{associated\/},
and are equivalent as $\so_L$-modules.
The $\so_L$-module $W_\la$ is irreducible unless $\la=\la^\ast$,
that is, unless $2\la'_1=L$. In the last case $W_\la$ splits into
a direct sum of two irreducible $\so_L$-modules.
The irreducible $\so_L$-modules corresponding to different partitions  
$\la$ are pairwise non-equivalent,
except for the pairs of modules $W_\la$ and $W_{\la^\ast}$
with $\la\neq\la^\ast\ts$.
In this way one obtains all irreducible finite-dimensional
non-spinor $\so_L$-modules, see \cite[Section V.9]{W}.
We can choose the Borel subalgebra, and fix the basis
in the corresponding Cartan subalgebra of $\so_L$ as in
\cite[Subsection 1.1]{KT}. If $2\la'_1=L\ts$, then $W_\la$ splits into
a direct sum of two irreducible $\so_L$-modules of highest weights
$$
(\la_1\lcd\la_{\ts L/2-1}\com\la_{\ts L/2})
\quad\textrm{and}\quad 
(\la_1\lcd\la_{\ts L/2-1}\com-\la_{\ts L/2})\,.
$$
If $2\la'_1<L\ts$, then $W_\la$ is
the irreducible $\so_L$-module of the highest weight
$(\la_1\lcd\la_{\ts [L/2]})$.
But again, 
the highest weight theory of $\so_L$-modules is not used
in the present article.

Observe that when $G_L=O_L\ts$, the first column length
of at least one of the two partitions $\la$ and $\la^\ast$ 
does not exceed $L/2$. Therefore, when regarding $W_\la$ as a $\so_L$-module,
we can assume that $2\la'_1\le L$. Under this assumption,
for any two distinct indices 
$i\com j\in\{1\lcd l\}$ we have
$$
c_i+c_j\ge3-2\la'_1\ge3-L\,.
$$
Then every factor
$\Rb_{\ts ij}(\ts c_i+t_i\com c_j+t_j\ts)$
in the product (\ref{3.1}) is regular at $t_1=\ldots=t_l=-\ts\frac12\,$,
for any choice of the standard tableau $\La$ of shape $\la$.
By Proposition 2.2, we then obtain from (\ref{3.1}) the explicit formula 
$$
F_\La\ =
\prod_{1\le i<j\le l}^{\longrightarrow}\, 
\left(1-\frac{Q_{ij}}
{\ts c_i+c_j+L-1}\ts\right)
\cdot E_{\La}
\quad\text{for}\quad
\g_L=\so_L\,,
$$
which is similar to the explicit formula (\ref{3.2222}).


\section{Yangian representations}\label{S4}

\noindent\textbf{4.1.}
Using the formal power series (\ref{1.31}) in $x^{-1}$,
introduce the elements of the algebra
$(\ts\End(\CC^N))^{\ts\ot\ts2}\ot\YN\,[[x^{-1}]]$
$$
T_1(x)=\sum_{i,j=1}^N\, E_{ij}\ot 1\ot T_{ij}(x)
\quad\textrm{and}\quad
T_2(x)=\sum_{i,j=1}^N\, 1\ot E_{ij}\ot T_{ij}(x)\,.
$$
The defining relations (\ref{1.32}) of the algebra $\YN$ are equivalent to
\begin{equation}\label{4.1}
R_{\ts12}(x\com y)\ts\,T_1(x)\ts\,T_2(y)
\ts=\,
T_2(y)\ts\,T_1(x)\ts\,R_{\ts12}(x\com y)\,.
\end{equation}
After multiplying both sides of the equality (\ref{4.1}) by $x-y\ts$,
it becomes an equality of formal Laurent series in $x^{-1}$ and $y^{-1}$
with the coefficients from the algebra
$(\ts\End(\CC^N))^{\ts\ot\ts2}\ot\YN$.
In (\ref{4.1}), we identify the element $R_{12}(x\com y)$ of 
$\End((\CC^N)^{\ts\ot\ts2})(x\com y)$ as defined by (\ref{3.45}),
with the element
$$
R_{12}(x\com y)\ot1\ts\in\ts(\ts\End(\CC^N))^{\ts\ot\ts2}(x\com y)\ot\YN\,.
$$

Using the defining relations of $\YN$ in the form
(\ref{4.1}), together with (\ref{3.55}), one shows that 
the assignment (\ref{1.51}) defines an automorphism of the
algebra $\YN$. Evidently, this automorphism is involutive.

For any $z\in\CC$ consider the evaluation $\YN\ts$-module $V(z)$,
see (\ref{eval}). Let us denote by $\rho_z$ the corresponding homomorphism
$\YN\to\End(\CC^N)\ts$. Then $T(x)\mapsto R_{12}(x\com z)$
under the homomorphism
\begin{equation}\label{4.15}  
\id\ot\rho_z:\ts\End(\CC^N)\ot\YN\to(\ts\End(\CC^N))^{\ts\ot\ts2}\,;
\end{equation}
see the definitions (\ref{1.71}),(\ref{tau}),(\ref{1.52}) and (\ref{3.45}).
More generally, consider the tensor product of evaluation $\YN\ts$-modules
$V(z_1)\ot\ts\ldots\ot V(z_n)\ts$, for any $z_1\lcd z_n\in\CC$.
The corresponding homomorphism
\begin{equation}\label{rhozz}
\rho_{z_1...\ts z_n}:\ts\YN\to(\ts\End(\CC^N))^{\ts\ot\ts n}
\end{equation}
is the composition of the map
$\rho_{z_1}\!\ot\ts\ldots\ot\ts\rho_{z_n}$
with $n\ts$-fold comultiplication map $\YN\to\YN^{\ts\ot\ts n}$.  
By definition (\ref{1.33}), then
\begin{equation}\label{4.2}
T(x) \mapsto R_{12}(x\com z_1)\ts\ldots\ts  R_{1,n+1}(x\com z_n)
\end{equation}
under the homomorphism
\begin{equation}\label{idrhozz} 
\id\ot\rho_{z_1...\ts z_n}:\,
\End(\CC^N)\ot\YN\to(\ts\End(\CC^N))^{\ts\ot\ts(n+1)}\ts.
\end{equation}

Now consider the standard action of the algebra
$\UN$ on $\CNn$. Denote by $\varpi_n$ the homomorphism $\UN\to\End(\CNn)$.
Using the matrix units $E_{ij}\in\End(\CC^N)\ts$,
we can write the operator of transposition of the
tensor factors in $(\CC^N)^{\ts\ot\ts2}$ as the sum
$$
\sum_{i,j=1}^N\ts E_{ij}\ot E_{ji}\,. 
$$
So by the definition (\ref{1.52}) of the homomorphism
$\al_N:\YN\to\UN\ts$,
\begin{equation}\label{4.25}
T(x)\,\mapsto\,
1\,-\,\sum_{k=1}^n\,\ts\frac{P_{1,k+1}}x
\end{equation}
under the homomorphism
\begin{equation}\label{idvarpial}
\id\ot(\varpi_n\circ\ts\al_N):
\End(\CC^N)\ot\YN\to(\ts\End(\CC^N))^{\ts\ot\ts(n+1)}\ts.
\end{equation}


\noindent\textbf{4.2.}
Now take any standard tableau $\Om$ with $n$ boxes.
Consider the operator $E_\Om$ on the vector space $\CNn$,
as defined in Subsection 1.1. 
We denote by $P_0$ the linear operator on $\CNn$ reversing the order of
the tensor factors. Setting $z=0$ in the following proposition,
we obtain Proposition 1.5. 

\begin{Proposition}
Put $z_k=c_k(\Om)+z$ for\/ $k=1\lcd n\ts$. Then
the operator $E_\Om\ts P_0$ is an intertwiner of the $\YN$-modules
\begin{equation}\label{4.3}
V(z_n)\ot\ldots\ot V(z_1)
\,\ts\longrightarrow\,
V(z_1)\ot\ldots\ot V(z_n)\,.
\end{equation}
\end{Proposition}

\begin{proof}
The action of $\YN$ on the module at the 
right-hand side of (\ref{4.3}) can be described
by the assignment (\ref{4.2}). The action of $\YN$ on the module
at the left-hand side of (\ref{4.3}) can be described
by the assignment
\begin{equation}\label{4.44}
T(x)\mapsto R_{12}(x\com z_n)\ts\ldots\ts  R_{1,n+1}(x\com z_1)\,.
\end{equation}

Take $n$ complex variables $x_1\lcd x_n$. Using
(\ref{3.5}) repeatedly, we obtain the equality
of rational functions in $x\com x_1\lcd x_n$
$$
R_{12}(x\com x_1)\ts\ldots\ts  R_{1,n+1}(x\com x_n)\,\cdot
\prod_{1\le i<j\le n}^{\longrightarrow}\,
R_{\ts i+1,j+1}(\ts x_i\com x_j\ts)\,\cdot\,(1\ot P_0)
$$
\vskip-16pt
$$
=\ 
\prod_{1\le i<j\le n}^{\longrightarrow}\,
R_{\ts i+1,j+1}(\ts x_i\com x_j\ts)\,\cdot\,(1\ot P_0)\,\cdot\,
R_{12}(x\com x_n)\ts\ldots\ts R_{1,n+1}(x\com x_1)\,,
$$
with values in the tensor product $(\ts\End(\CC^N))^{\ts\ot\ts(n+1)}$.
Using the constrained variables $t_1(\Om)\lcd t_n(\Om)$ from
Subsection 1.2, put
$$
x_k=c_k(\Om)+t_k(\Om)
\quad\text{for each}\quad
k=1\lcd n\,.
$$
Then at $t_1(\Om)=\ldots=t_n(\Om)=z$ we have $x_k=z_k$ for
$k=1\lcd n$. By Theorem 1.2, the above equality of rational
functions in $x\com x_1\lcd x_n$ yields
$$
R_{12}(x\com z_1)\ts\ldots\ts  R_{1,n+1}(x\com z_n)
\,\cdot\,(1\ot E_\Om\ts P_0)
$$

\vskip-14pt
\begin{equation}\label{4.5}
=\
(1\ot E_\Om\ts P_0)\,\cdot\,
R_{12}(x\com z_n)\ts\ldots\ts R_{1,n+1}(x\com z_1)\,.
\end{equation}

\vskip2pt
\noindent
In view of (\ref{4.2}) and (\ref{4.44}), the equality (\ref{4.5})
proves Proposition 4.2.
\qed
\end{proof}

Due to Proposition 4.2, for any $z\in\CC$ the subspace $V_\Om$ in $\CNn$ 
can be regarded as a submodule
in the tensor product of evaluation $\YN\ts$-modules

\vskip-16pt
$$
V(c_1(\Om)+z)\ot\ts\ldots\ot V(c_n(\Om)+z)\,.
$$
This fact will be also used in Section 5, for the particular value
$z=\frac{M}2\mp\frac12\,$.


\smallskip\medskip\noindent\textbf{4.3.}
Consider the embedding $\UN\to\YN$ as defined by ({\ref{4.4}).
By the definition (\ref{1.52}),
the homomorphism $\al_N:\YN\to\UN$ is identical on the subalgebra 
$\UN\subset\YN\ts$.
The automorphism $\xi_N$ of $\YN$
is also identical on this subalgebra, see the definition (\ref{1.51}).
It follows that the restriction of $\al_{NM}$ to this subalgebra
coincides with the natural embedding $\UN\to\UMN\ts$.
For the particular choice (\ref{1.62}) of the series $g(x)$,
the automorphism (\ref{1.61}) is identical on the subalgebra
$\UN\subset\YN\ts$, because the coefficient of $g_\mu(x)$ at 
$x^{-1}$~is~$0\ts$.
So the action of the subalgebra
$\UN\subset\YN\ts$ on the $\YN\ts$-module $\Vlm$ coincides
with its natural action, corresponding to the natural
action of the group $GL_N$ on $\Vlm$.

Furthermore, the action of the subalgebra
$\UN\subset\YN\ts$ on any evaluation module $V(z)$ over $\YN$
coincides with the natural action of $\UN$
on the vector space $\CC^N$ of $V(z)\ts$.
The assigment (\ref{4.4}) determines a Hopf algebra
embedding, because by (\ref{1.33}) we have
$$
\De\ts\bigl(\ts T_{ij}^{(1)}\bigr)\ts=\,
T_{ij}^{(1)}\ns\ot\ts1+1\ot\ts T_{ij}^{(1)}\,.
$$
Therefore the action of the subalgebra
$\UN\subset\YN\ts$ on the tensor product of evaluation $\YN\ts$-modules 
$V(c_1(\Om))\ot\ldots\ot V(c_n(\Om))$ 
coincides with the natural action of $\UN$
on the vector space $\CNn$. The $\YN\ts$-module $V_\Om$ is
defined as a submodule of this tensor product of evaluation modules.
Hence the action of $\UN\subset\YN\ts$ on this submodule
coincides with the natural action of $\UN$ on the subspace
$V_\Om\subset\CNn$.


\smallskip\medskip\noindent\textbf{4.4.}
In the remainder of this section we prove Theorem 1.6.
Fix a standard tableau $\La$ of non-skew shape 
$\la\ts$, such that the tableau $\Om$ is obtained from $\La$
by removing the boxes with numbers $1\lcd m\ts$. Here $\la$
is a\ partition of $l$, and $m=l-n$. As in Sections 2 and 3,
write $c_k=c_k(\La)$ for $k=1\lcd l$. The standard tableau 
of non-skew shape $\mu\ts$, obtained by removing the boxes
with the numbers $m+1\lcd l$ from the tableau $\La\ts$, is denoted by $\Up$. 

Put $L=N+M$. Take the vector space $\CC^{L}=\CC^N\ns\op\ts\CC^M$.
In the present subsection
the rational functions $R_{12}(x\com y)\ts\lcd\ts R_{1,l+1}(x\com y)$
of the variables $x$ and $y$ will take values in the algebra
$(\ts\End(\CC^L))^{\ts\ot\ts(l+1)}\ts$. 

We always denote by $E_\La$ the operator
on the vector space $\CLl$, corresponding to
$e_\La\in\CC S_l\ts$. By Proposition 2.5, we have the equality

\vskip-4pt
$$
R_{12}(x\com c_1)\ts\ldots\ts  R_{1,l+1}(x\com c_l)
\ts\cdot\ts(\ts1\ot E_\La)
$$

\vskip-10pt
\begin{equation}\label{4.6}
=\ 
\biggl(1\ts-\ts\sum_{k=1}^l\,\frac{P_{\ts1,k+1}}x\ts\biggr)
\cdot\ts(\ts1\ot E_\La)\,\ts,
\end{equation}

\vskip-2pt
\noindent
of rational functions in $x$,
with values in the algebra $(\ts\End(\CC^L))^{\ts\ot\ts(l+1)}\ts$.
The operator $P_{\ts1,k+1}$ on 
$(\CC^L))^{\ts\ot\ts(l+1)}$ corresponds to $(1\,k+1)\in S_{\ts l+1}\ts$,
see (\ref{3.45}). 
The equality (\ref{4.6}) is the starting point for our proof of Theorem 1.6. 

First consider the case, when $M=0$ and
$\mu=(0\com0\ts,\ts\ldots\ts)\ts$. In this case $N=L$, $n=l$
and the standard tableau $\Om=\La$ has a non-skew shape.
In this case, the left-hand side of the equality (\ref{4.6})
describes the action of the Yangian $\YL$ on the submodule $V_\Om=V_\La$
of the tensor product $V(c_1)\ot\ts\ldots\ot V(c_l)$
of evaluation $\YL\ts$-modules; see (\ref{4.2}).
The right-hand side of (\ref{4.6})
describes the action of $\YL$ on the module $\Vlm=V_\la\ts$,
defined in Subsection~1.6\ts; see (\ref{4.25}).
Indeed, here $g_\mu(x)=1$ and $\al_{NM}=\al_N$. The image of the operator
$E_\Om=E_\La\,$, as a $\glL\ts$-submodule in $\CLl$,
is equivalent to the $\glL\ts$-module~$V_\la\ts$. 
Thus the equality (\ref{4.6}) shows that the $\YL\ts$-modules
$V_\La$ and $V_\la$ are equivalent.

Let us now prove Theorem 1.6 for $M\ge1$. As in Subsection 2.4,
split the vector space

\vskip-16pt
$$
\CLl=(\CC^N\ns\op\ts\CC^M)^{\ts\ot\ts l}
$$
into the direct sum of subspaces, obtained
from
$
(\CC^M)^{\ot\ts k}\ot\ts(\CC^N)^{\ot\ts(l-k)}
$
by some permutations of the $l$ tensor factors. Here $k=0\lcd l$. 
Consider the projector onto the direct summand (\ref{4.7}),
\begin{equation}\label{4.65}
I_m:\,\CLl\,\to\,\CMm\ot\CNn\,.
\end{equation}
\noindent
Note that the operator (\ref{4.65}) is $GL_N\times GL_M\ts$-equivariant;
here we use the embedding $GL_N\times GL_M\to GL_L$ chosen in the
beginning of Subsection~1.6.

Consider the $\YL\ts$-module,
obtained by pulling the tensor product of the evaluation $\YL\ts$-modules
corresponding to $c_1\lcd c_l$ back through the automorphism $\xi_L$ of
$\YL\ts$; see (\ref{1.51}). The action of
$\YL$ on this module is described by the assignment 
\begin{equation}\label{vlup}
\sum_{i,j=1}^L\,E_{ij}\ot T_{ij}(x)\,\ts\mapsto\,\,\ts
R_{1,l+1}(-x\com c_l)^{-1}\ldots\ts R_{12}(-x\com c_1)^{-1}
\end{equation}
\begin{equation}\label{vl}
=\ 
f(x)\,\cdot\ts
R_{1,l+1}(\ts x\com-c_l)\ts\ldots\ts R_{12}(\ts x\com-c_1)\,,
\end{equation}
where
\begin{equation}\label{fx}
f(x)\ =\ \prod_{k=1}^l\,\frac{(x+c_k)^2}{\ts(x+c_k)^2-1\ts}\ ;
\end{equation}

\vskip-2pt\noindent
see (\ref{3.55}) and (\ref{4.2}). Denote by $V_l$
the restriction of this $\YL\ts$-module
to the subalgebra $\YN\subset\YL\ts$; here
we use the natural embedding $\ph_M:\YN\to\YL\ts$.
That is, $\ph_M:\ts T_{ij}(x)\ts\mapsto\ts T_{ij}(x)$ 
for $1\le i\com j\le N$ by definition.
Further, denote by $\Val$ the $\YN\ts$-module
obtained by pulling the $\YN\ts$-module $V_l$ back through
the automorphism $\xi_N$ of $\YN\ts$.
Note that the vector space of the $\YN\ts$-modules
$V_l$ and $\Val$ is $\CLl$.

Take the vector space $\CMm\ot\CNn$.
For $k=m+1\lcd l\/$ put
\begin{equation}\label{PV}
R_{1,k+1}^{\,\,\wedge}(x\com y)\,=\,\ts
1-\frac{\,P_{1,k+1}^{\,\wedge}\,}{x-y}
\end{equation}
where $P_{1,k+1}^{\,\wedge}$ denotes the operator on
$\CC^N\ot\CMm\ot\CNn$, acting as transposition in the first
and $(k+1)$-th tensor factors, and acting as the identity
in the remaining $l-1$ tensor factors. One can define
an action of the algebra $\YN$ on the vector space $\CMm\ot\CNn$
by the assignment
$$
\sum_{i,j=1}^N\,E_{ij}\ot T_{ij}(x)\,\ts\mapsto\,\,\ts
f(x)\,\cdot\,
R_{1,l+1}^{\,\,\wedge}(\ts x\com-c_l)\,\ldots\ts 
R_{1,m+2}^{\,\,\wedge}(\ts x\com-c_{m+1})\,.
$$

Denote by $V_{mn}$ the $\YN\ts$-module defined 
by the above displayed assignment.
Further, denote by $\Vamn$ the $\YN\ts$-module
obtained by pulling the $\YN\ts$-module $V_{mn}$ back through
the automorphism $\xi_N$ of $\YN\ts$. 
The action of the algebra $\YN$ on $\Vamn$
is then described by the assignment
$$
\sum_{i,j=1}^N\,E_{ij}\ot T_{ij}(x)\,\ts\mapsto\,\,\ts
h(x)\,\cdot\,
R_{1,m+2}^{\,\,\wedge}(x\com c_{m+1})\,\ldots\ts 
R_{1,l+1}^{\,\,\wedge}(x\com c_l)\,,
$$
where
\begin{equation}\label{gx}
h(x)\ =\  
f(-x)^{-1}\,\cdot\!
\prod_{k=m+1}^l\!\frac{(x-c_k)^2}{\ts(x-c_k)^2-1\ts}
\ = \ 
\prod_{k=1}^m\,\frac{\ts(x-c_k)^2-1\ts}{(x-c_k)^2}\ .
\end{equation}

Note that the $\YN\ts$-module $\Vamn$
can also be obtained as follows. Take the $\YN\ts$-module,
obtained by pulling the tensor product of evaluation $\YN\ts$-modules
with the parameters
$
c_{m+1}=c_1(\Om)\,\lcd\,c_l=c_n(\Om)
$
back through the automorphism
$T_{ij}(x)\mapsto h(x)\,T_{ij}(x)$ of the algebra $\YN$.
The vector space of this $\YN\ts$-module is $\CNn$.
Then by regarding the tensor product $\CMm\ot\CNn$ as
$\YN\ts$-module where every element of the Hopf algebra
$\YN$ acts on $\CMm$ via the counit homomorphism $\varepsilon:\YN\to\CC$,
we obtain the $\YN\ts$-module $\Vamn\ts$.

\begin{Proposition}
The projection {\rm (\ref{4.65})} is
an intertwining operator of\/ $\YN\ts$-modules $V_l\to V_{mn}\ts$.
\end{Proposition}

\begin{proof}
This result follows
by comparing (\ref{vlup}),(\ref{vl})
with the definition of the $\YN\ts$-module
$V_{mn}\ts$, and by using Lemma 2.5 in the case
\begin{equation}\label{zc}
z_1=-c_1\ts\lcd\ts z_l=-c_l\,.\quad\qed
\end{equation}
\end{proof}

\begin{Corollary}
The projection {\rm (\ref{4.65})} is
an intertwining operator of\/ $\YN\ts$-modules $\Val\to\Vamn\ts$.
\end{Corollary}

We end this subsection with the next lemma.
Consider the rational functions $g_\mu(x)$ and $h(x)\ts$, 
defined by (\ref{1.62}) and (\ref{gx}), respectively.

\begin{Lemma}
We have the equality\/ $g_\mu(x)\,h(x)=1$.
\end{Lemma}

\begin{proof}
Consider the product at the right-hand side of the equalities
(\ref{gx}). This product is symmetric in $c_1\lcd c_m$ and
therefore does not depend on the choice of a standard tableau $\Up$
of shape $\mu$. We choose $\Up$ to be the row tableau
of shape $\mu\ts$. For any index $i\ge1$, take the boxes
in the $i\,$th row of the Young diagram $\mu$ in their natural order,
from the leftmost to the rightmost box. The contents of these
boxes form the sequence $1-i\lcd\mu_i-i$ which is increasing by $1$;
this sequence is empty if $\mu_i=0$. Therefore
$$
h(x)\ = \ 
\prod_{k=1}^m\,\left(
\frac{\ts x-c_k+1\ts}{x-c_k}\ts\cdot\ts
\frac{\ts x-c_k-1\ts}{x-c_k}\right)
$$
$$
=\ \,
\prod_{i\ge1}\,\left(
\frac{\ts x+i\ts}{x-\mu_i+i}\ts\cdot\ts
\frac{\ts x-\mu_i+i-1\ts}{x+i-1}\right)\,= \
g_\mu(x)^{-1}.\quad\qed
$$ 
\end{proof}


\noindent\textbf{4.5.}
Let us continue our proof of Theorem 1.6.
Consider the image of the subspace (\ref{4.7})
under the operator $E_\La$ on $\CLl$. Note that
this image is contained in the subspace $V_\La\subset\CLl$.
Consider the $\YN\ts$-module $V_l$ defined in Subsection 4.4\ts;
the vector space of this module is $\CLl$.

\begin{Proposition}
The image of the subspace\/ {\rm (\ref{4.7})} under the operator $E_\La$
is an\/ $\YN$-submodule of $V_l\ts$.
\end{Proposition}

\begin{proof}
The action of the coefficients of the series 
$T_{ij}(x)$ with $1\le i\com j\le N$
on the $\YN\ts$-module $V_l$ is described by the assignment (\ref{vlup}). 
Here we use the natural embedding $\ph_M:\YN\to\YL\ts$. 
Consider the product (\ref{vl}) in the algebra
$(\ts\End(\CC^L))^{\ts\ot\ts(l+1)}(x)\ts$.
We have a relation in this algebra,
$$
R_{1,l+1}(x\com-c_l)\ts\ldots\ts R_{12}(x\com-c_1)
\ts\cdot\ts(\ts1\ot E_\La)\,=\, 
(\ts1\ot E_\La)
$$
\begin{equation}\label{4.8}
\times\ \ts 
R_{12}(x\com-c_1)\ts\ldots\ts R_{1,l+1}(x\com-c_l)\,;
\end{equation}
see (\ref{4.5}). Expand the product in the second line of the display
(\ref{4.8}), as

\vskip-3pt
$$
R_{12}(x\com-c_1)\ts\ldots\ts R_{1,l+1}(x\com-c_l)\,=\,
\sum_{i,j=1}^L\,E_{ij}\ot A_{ij}(x)
$$

\vskip-2pt\noindent
for certain rational functions
$A_{ij}(x)\in(\ts\End(\CC^L))^{\ts\ot\ts l}(x)\ts$.
Then consider the restrictions of the operator values of the
functions $A_{ij}(x)$ with $1\le i\com j\le N$
to the subspace (\ref{4.7}). Using Corollary 2.5 
in the case (\ref{zc}), we prove that

\vskip-4pt
$$
\sum_{i,j=1}^N\,E_{ij}\ot (\ts A_{ij}(x)\,|\,\CMm\ot\CNn)
$$

\vskip-4pt
$$
=\ \, 
R_{1,m+2}^{\,\,\wedge}(x\com-c_{m+1})\,\ldots\ts 
R_{1,l+1}^{\,\,\wedge}(x\com-c_l)\,.
$$

\vskip7pt\noindent
In particular, the operator values of the functions $A_{ij}(x)$
with $1\le i,j\le N$ preserve the subspace (\ref{4.7}).
Now Proposition~4.5 follows from (\ref{4.8}).
\qed
\end{proof}

The $\YN\ts$-module $\Val$ is obtained from $V_l$
by pulling back through an automorphism of $\YN\ts$.
Therefore Proposition 4.5 has a corollary.

\begin{Corollary}
The image of the subspace\/ {\rm (\ref{4.7})} under the operator $E_\La$
is an\/ $\YN$-submodule of $\Val\ts$.
\end{Corollary}


\noindent\textbf{4.6.}
In this subsection we complete the proof of Theorem 1.6.
Consider the image of the subspace (\ref{4.7}) under the linear operator
$$
I_m\ts E_\La:\CLl\to\CMm\ot\CNn\,.
$$ 
This image coincides with the vector subspace
\begin{equation}\label{4.9}
V_{\ts\Up}\ot V_\Om\,\subset\,\CMm\ot\CNn\ts.
\end{equation}
Indeed,

\vskip-12pt
$$
I_m\ts E_\La\ts\,|\,\ts\CMm\ot\CNn\ts=\ts\,E_{\ts\Up}\ot E_\Om\,\ts;
$$

\vskip4pt
\noindent
see (\ref{2.85}) and (\ref{2.9}). It now follows from
Corollaries 4.4 and 4.5 that the vector subspace (\ref{4.9})
is a submodule in the $\YN\ts$-module $\Vamn\ts$. 
Let us denote this submodule of $\Vamn$ by $V$.
Note that $V$ is a subquotient of the $\YN\ts$-module $\Val$
by definition.

The description of the $\YN\ts$-module $\Vamn$ given after (\ref{gx})
yields the following description of the
$\YN\ts$-module $V$. 
Take the $\YN\ts$-module $V_\Om$ as defined in Subsection 1.5. 
Pull $V_\Om$ back through the automorphism
$T_{ij}(x)\mapsto h(x)\,T_{ij}(x)$ of $\YN\ts$.
Extend the resulting action of $\YN$ on the vector
space $V_\Om$ to the vector space $V_{\ts\Up}\ot V_\Om\ts$ so that
every element of $\YN$ acts on $V_{\ts\Up}$
as the identity. Then we obtain the $\YN\ts$-module $V$.

The subspace $V_{\ts\Up}\subset\CMm$ is equivalent to $V_\mu$
as a representation of the group $GL_M\ts$.
Consider the subspace $V_\La\subset\CLl$ as a
representation of $GL_L\ts$,
equivalent to $V_\la$. Then regard $V_\La$ as $\YN\ts$-module
by pulling back through the homomorphism
$\al_{NM}:\YN\to\UL\ts$; see definition (\ref{1.69}).

\begin{Proposition}
\ns\!The $\YN$-module $V\!$ is a subquotient of\/ $\YN$-module~$V_\La$.
\end{Proposition}

\begin{proof}
By the definition (\ref{1.69}), we have
$$
\al_{NM}=\,\al_L\circ\ts\xi_L\circ\ts\ph_M\circ\ts\xi_N\ts.
$$
Consider $V_\La$ as a submodule in the tensor product of
the evaluation $\YL\ts$-modules with the parameters
$c_1\lcd c_l$. We have already shown that the action of
$\YL$ on this submodule factors through the homomorphism
$\al_L:\YL\to\UL\ts$. So the $\YN\ts$-module $V_\La$ as defined above
can also be obtained by pulling the action of $\YL$ on $V_\La$
back through the injective homomorphism
$$
\xi_L\circ\ts\ph_M\circ\ts\xi_N:\YN\to\YL\,.
$$
Thus $V_\La$ is a submodule in the $\YN\ts$-module $\Val$.
But by definition, $V$ is a quotient of a certain $\YN\ts$-submodule of
$\Val$. The latter submodule of $\Val$ is contained in $V_\La\,$.
\qed
\end{proof}

Consider the restriction of the representation $V_\La$ of 
the group $GL_L$
to the subgroup $GL_M\ts$. Realize the vector space (\ref{1.0}) as
\begin{equation}\label{4.10}
{\rm Hom}_{\,GL_M}(\ts V_{\ts\Up}\ts\com V_\La\ts)\ts.
\end{equation}
Since the image of the homomorphism $\al_{NM}$ is contained in
the subalgebra of $GL_M\ts$-invariants $\AMN\subset\UL$, the action
of the algebra $\YN$ on $V_\La$ induces an action of $\YN$ on (\ref{4.10}).
This action of $\YN$ on (\ref{4.10})
is irreducible, see \cite[Section 2]{MO}.

The operator (\ref{4.65}) is $GL_N\times GL_M\ts$-equivariant,
and the vector space $V_\Up\ot V_\Om$ of the $\YN\ts$-module $V$
has a natural action 
of the groups $GL_N$ and $GL_M\ts$. The action of $GL_M$ on $V$
commutes with the action of the algebra $\YN\ts$.
By Proposition 4.6, the $\YN\ts$-module 
\begin{equation}\label{4.11}
{\rm Hom}_{\,GL_M}(\ts V_{\ts\Up}\ts\com V\ts)
\end{equation}
is a subquotient of 
(\ref{4.10}). Since the $\YN\ts$-module (\ref{4.10})
is irreducible, it must be equivalent to the $\YN\ts$-module
(\ref{4.11})\ts; see Proposition 2.4.

The $\YN\ts$-module (\ref{4.11})  
can also be obtained by pulling the $\YN\ts$-module $V_\Om\ts$,
as defined in Subsection 1.5, back through the automorphism
(\ref{1.61}) of $\YN\ts$, where $g(x)=g_\mu(x)^{-1}$.
Here we use Lemma 4.4. The proof of Theorem 1.6
is now complete. 

Note that (\ref{4.11}) is also a subquotient of (\ref{4.10})
as a representation of the group $GL_N\ts$. Thus we obtain
Proposition 1.1 together with Theorem 1.6. Proposition 1.1 could
be proved independently of Theorem 1.6. We chose the present
proofs, because they have analogues for the classical groups
$O_N$ and $Sp_N$ instead of $GL_N$. These analogues will be given
in the next section. 


\section{Twisted Yangians}\label{S5}

\textbf{5.1.}
Using the formal power series (\ref{1.771}) in $x^{-1}$,
introduce the elements of the algebra
$\End(\CC^N)^{\ts\ot\ts2}\ot\XN\,[[x^{-1}]]$
$$
S_1(x)=\sum_{i,j=1}^N\, E_{ij}\ot 1\ot S_{ij}(x)
\quad\textrm{and}\quad
S_2(x)=\sum_{i,j=1}^N\, 1\ot E_{ij}\ot S_{ij}(x)\,.
$$
The defining relations of the algebra $\XN$ can be written as
\begin{equation}\label{5.1}
R_{\ts12}(x\com y)\,S_1(x)\,\Rt_{12}(x\com y)\,S_2(y)
\ts=\ts
S_2(y)\,\Rt_{12}(x\com y)\,S_1(x)\,R_{\ts12}(x\com y)\,.
\end{equation}
After multiplying both sides of the equality (\ref{5.1}) by $x^2-y^2\ts$,
it becomes an equality of formal Laurent series in $x^{-1}$ and $y^{-1}$
with the coefficients from the algebra $\End(\CC^N)^{\ts\ot\ts2}\ot\XN$.
In (\ref{5.1}),
we identify the elements $R_{12}(x\com y)$ and $\Rt_{12}(x\com y)$
of $\End(\CC^N)^{\ts\ot\ts2}(x\com y)$ as defined by 
(\ref{3.45}) and (\ref{3.555}), respectively with the elements 
$R_{12}(x\com y)\ot1$ and $\Rt_{12}(x\com y)\ot1$ of
$$
\End(\CC^N)^{\ts\ot\ts2}(x\com y)\ot\XN\,.
$$

Using (\ref{3.55}) and (\ref{3.6}) where $L=N$,
one derives from the defining relations (\ref{5.1}) that the assignment
(\ref{1.751}) defines an automorphism of the algebra $\XN$.
For details of this argument see \cite[Subsection 6.5]{MNO}.

By dividing each side of the equality (\ref{5.1})
by $S_2(y)$ on the left and right, and then setting $y=-x$, 
we obtain the equality
$$
Q_{12}\,S_1(x)\,R_{\ts12}(x\com -x)\,S_2(-x)^{-1}=\ts
S_2(-x)^{-1}\ts R_{\ts12}(x\com -x)\,S_1(x)\,Q_{12}\,.
$$
As in (\ref{5.1}),
here we identify the element $Q_{12}\in\End(\CC^N)^{\ts\ot\ts2}$ with
$$
Q_{12}\ot1\ts\in\ts\End(\CC^N)^{\ts\ot\ts2}\ot\XN\,.
$$
Since the image of the operator $Q_{12}=Q(N)$ in $(\CC^N)^{\ts\ot\ts2}$
is one-dimensional,
the last displayed equality implies the existence of a formal power series
$D(x)$ in $x^{-1}$ with coefficients in $\XN$ and leading term $1$,
such that
\begin{equation}\label{5.2}
Q_{12}\,S_1(x)\,R_{\ts12}(x\com -x)\,S_2(-x)^{-1}=\ts
\Bigl(1\mp\frac1{2x}\ts\Bigr)\,D(x)\, Q_{12}\,.
\end{equation}
By using (\ref{3.55}) when $y=-x$,
one derives from (\ref{5.2}) that
$D(x)\ts D(-x)=1$.

By \cite[Theorem 6.3]{MNO}, all coefficients of the series $D(x)$
belong to the centre of the algebra $\XN$. By \cite[Theorem 6.4]{MNO},
the kernel of the surjective homomorphism $\pi_N:\XN\to\YS$ is generated
by the coefficients of the series $1-D(x)$.
For any series $g(x)\in\CC[[x^{-1}]]$ with the leading term $1$,
the definition (\ref{5.2}) of the series $D(x)$
shows that the assignment (\ref{1.861}) determines an automorphism
of the quotient algebra $\YS$ of $\XN$, if and only if $g(x)=g(-x)$.


\medskip\noindent\textbf{5.2.}
For any $z\in\CC\ts$, consider the restriction of the
evaluation $\YN\ts$-module $V(z)$ to the subalgebra $\YS\subset\YN$.
By definition, this subalgebra is generated by the coefficients
of all the formal power series from $\YN\,[[x^{-1}]]$, appearing in
the expansion of the element (\ref{1.72}) relative to the basis
of matrix units $E_{ij}$ in $\End(\CC^N)\ts$. Under the 
homomorphism (\ref{4.15}) corresponding to the $\YN\ts$-module $V(z)\ts$,
we have $\Tt(x)\mapsto\Rt_{12}(x\com z)\ts$;
see (\ref{3.555}) and Subsection 4.1. Therefore, 
under the homomorphism (\ref{4.15})
\begin{equation}\label{5.21}
\Tt(x)\,T(x)\ts\mapsto\ts\Rt_{12}(x\com z)\,R_{12}(x\com z)\,.
\end{equation}

Now consider the twisted evaluation $\YN\ts$-module $\Vt(z)$,
see (\ref{teval}). Let us denote by $\tilde{\rho}_z$ 
the corresponding homomorphism $\YN\to\End(\CC^N)\ts$.
Then $T(x)\mapsto\Rt_{12}(x\com z)$
under the homomorphism
$$ 
\id\ot\tilde\rho_z:\ts\End(\CC^N)\ot\YN\to(\ts\End(\CC^N))^{\ts\ot\ts2}\ts.
$$
Therefore under this homomorphism
\begin{equation}\label{5.22}
\Tt(x)\,T(x)\ts\mapsto\ts R_{12}(x\com z)\,\Rt_{12}(x\com z)\,.
\end{equation}
The obvious equality of the right-hand sides of (\ref{5.21}) and 
(\ref{5.22}) explains why the restrictions of the $\YN\ts$-modules
$V(z)$ and $\Vt(z)$ to the subalgebra $\YS\subset\YN$ coincide.

Consider the restriction of
the tensor product of evaluation $\YN\ts$-modules
$V(z_1)\ot\ts\ldots\ot V(z_n)\ts$ to the subalgebra $\YS\subset\YN\ts$,
for any $z_1\lcd z_n\in\CC$.
The action of $\YN$ on this tensor product
defines the homomorphism (\ref{rhozz}). Then
$$
\Tt(x)\,T(x)\mapsto
$$
\begin{equation}\label{5.225}
\Rt_{1,n+1}(x\com z_n)\ts\ldots\ts\Rt_{12}(x\com z_1)\,
R_{12}(x\com z_1)\ts\ldots\ts  R_{1,n+1}(x\com z_n)
\end{equation}
under the homomorphism (\ref{idrhozz}), see the formula (\ref{4.2}).

Furthermore, consider the restriction of
the tensor product of twisted evaluation $\YN\ts$-modules
$\Vt(z_1)\ot\ts\ldots\ot\Vt(z_n)\ts$ to $\YS\subset\YN\ts$.
The action of $\YN$ on this tensor product
defines a homomorphism
$$
\tilde\rho_{z_1...\ts z_n}:\,
\YN\to(\ts\End(\CC^N))^{\ts\ot\ts n}\,.
$$
This homomorphism is the composition of the map
$\tilde\rho_{z_1}\!\ot\ts\ldots\ot\ts\tilde\rho_{z_n}$
with the $n\ts$-fold comultiplication map $\YN\to\YN^{\ts\ot\ts n}$.  
By  the definition (\ref{1.33}),
$$
T(x)\mapsto\Rt_{12}(x\com z_1)\ts\ldots\ts\Rt_{1,n+1}(x\com z_n)
$$
under the homomorphism
\begin{equation}\label{idtilrhozz} 
\id\ot\tilde\rho_{z_1...\ts z_n}:\,
\End(\CC^N)\ot\YN\to(\ts\End(\CC^N))^{\ts\ot\ts(n+1)}\ts.
\end{equation}
Therefore, under (\ref{idtilrhozz})
$$
\Tt(x)\,T(x)\mapsto
$$
\begin{equation}\label{5.226}
R_{1,n+1}(x\com z_n)\ts\ldots\ts R_{12}(x\com z_1)\,
\Rt_{12}(x\com z_1)\ts\ldots\ts\Rt_{1,n+1}(x\com z_n)\,.
\end{equation}
Note that when $n>1$, the images of the product $\Tt(x)\ts T(x)$ under
the homomorphisms (\ref{idrhozz}) and (\ref{idtilrhozz}) may be
non-equal; see Subsection 5.3.

Let us now consider the extended twisted Yangian $\XN\ts$, and the
homomorphism $\be_N:\XN\to\US$ defined by the 
assignment (\ref{1.752}). Consider the standard action of the algebra
$\US$ on the tensor product $\CNn$. The corresponding homomorphism 
$\US\to\End(\CNn)$ is just the restriction of the
homomorphism $\varpi_n:\UN\to\End(\CNn)$ to the subalgebra
$\US\subset\UN\ts$, see the end of Subsection 4.1. 
Using the matrix units $E_{ij}\in\End(\CC^N)\ts$,
we can write the operator $-\ts Q(N)$ on $(\CC^N)^{\ts\ot\ts2}$ as

\vskip-12pt
\begin{equation}\label{sisi}
\sum_{i,j=1}^N\ts E_{ij}\ot\si(E_{ji})\,=\,
\sum_{i,j=1}^N\ts\si(E_{ij})\ot E_{ji}\,. 
\end{equation}

\vskip-4pt\noindent
Hence

\vskip-12pt
\begin{equation}\label{5.23}
S(x)\,\mapsto\,
1\ts-\ts{\textstyle\bigl(\ts x\pm\frac12\ts\bigr)^{-1}}\,
\sum_{k=1}^n\,\ts(\ts P_{\ts1,k+1}-Q_{1,k+1}\ts)
\end{equation}

\vskip-6pt\noindent
under the homomorphism
$$
\id\ot(\varpi_n\circ\ts\be_N):
\End(\CC^N)\ot\XN\to(\ts\End(\CC^N))^{\ts\ot\ts(n+1)}\ts.
$$
Note that under the homomorphism (\ref{idvarpial}),
\begin{equation}\label{5.24}
\hskip25pt
\Tt(x)\,T(x)\,\mapsto\,
\Bigl(1\ts+\ts x^{-1}\,\sum_{k=1}^n\,Q_{\ts1,k+1}\Bigr)
\Bigl(1\ts-\ts x^{-1}\,\sum_{k=1}^n\,P_{\ts1,k+1}\Bigr)\ts;
\end{equation}
see (\ref{4.25}). The elements of the algebra 
$\End(\CC^N))^{\ts\ot\ts(n+1)}(x)$ at the right-hand sides
of the assignments (\ref{5.23}) and (\ref{5.24}) are different.
This difference was explained by \cite[Proposition 2.4]{N3},
see also Lemma 5.4 below. 


\smallskip\medskip\noindent\textbf{5.3.}
Now take any standard tableau $\Om$ of shape $\lm\ts$.
We assume that the partitions $\la$ of $l$ and $\mu$ of $m$
satisfy the conditions from \cite{W} for the groups $G_{M+N}$
and $G_M\ts$, respectively, as described in Subsection 1.3.
Consider the operator $\Fom$ on the vector space $\CNn\ts$; this
operator is defined by Theorem 1.4. Here $n=l-m$.
Let us also keep to the notation (\ref{1.777}).

\smallskip\medskip\noindent{\it Proof of Proposition 1.7.\/}
The action of the algebra $\YS$ on the tensor products
$V(d_1(\Om))\ot\ts\ldots\ot V(d_n(\Om))\ts$ and
$\Vt(d_1(\Om))\ot\ts\ldots\ot\Vt(d_n(\Om))\ts$ is 
explicitly described by the formulas (\ref{5.225}) and (\ref{5.226}),
respectively, where 
\begin{equation}\label{5.333}
z_1=d_1(\Om)\,\lcd\,z_n=d_n(\Om)\,.
\end{equation}

As in Subsection 4.2, take $n$ complex variables $x_1\lcd x_n$. Using
(\ref{3.5}), (\ref{3.7}) and (\ref{3.85}) for $L=N$ repeatedly, 
we obtain the equality
of rational functions in the variables $x\com x_1\lcd x_n$
$$
\Rt_{1,n+1}(x\com x_n)\ts\ldots\ts\Rt_{12}(x\com x_1)\,\,
R_{12}(x\com x_1)\ts\ldots\ts  R_{1,n+1}(x\com x_n)
$$
$$
\times\ 
\prod_{1\le i<j\le n}^{\longrightarrow}\,
\Rb_{\ts i+1,j+1}(\ts x_i\com x_j\ts)
\ \ \cdot\hskip-5pt
\prod_{1\le i<j\le n}^{\longrightarrow}\,
R_{\ts i+1,j+1}(\ts x_i\com x_j\ts)
$$
$$
=\ 
\prod_{1\le i<j\le n}^{\longrightarrow}\,
\Rb_{\ts i+1,j+1}(\ts x_i\com x_j\ts)
\ \ \cdot\hskip-5pt
\prod_{1\le i<j\le n}^{\longrightarrow}\,
R_{\ts i+1,j+1}(\ts x_i\com x_j\ts)
$$
$$
\times\ \ 
R_{1,n+1}(x\com x_n)\ts\ldots\ts R_{12}(x\com x_1)\,\,
\Rt_{12}(x\com x_1)\ts\ldots\ts\Rt_{1,n+1}(x\com x_n)
$$

\vskip4pt\noindent
with values in the tensor product $(\ts\End(\CC^N))^{\ts\ot\ts(n+1)}$.
Using the variables $t_1(\Om)\lcd t_n(\Om)$ constrained as in Theorem 1.4,
put
$$
x_k=c_k(\Om)+t_k(\Om)+{\textstyle\frac{M}2}
\quad\text{for each}\quad
k=1\lcd n\,.
$$
At $t_1(\Om)=\ldots=t_n(\Om)=\mp\frac12\ts$
the above equality of rational
functions in $x\com x_1\lcd x_n$ then yields the equality
of rational functions in $x$,
$$
\Rt_{1,n+1}(x\com d_n(\Om))\ts\ldots\ts\Rt_{12}(x\com d_1(\Om))\,\,
R_{12}(x\com d_1(\Om))\ts\ldots\ts R_{1,n+1}(x\com d_n(\Om))
$$
$$
\times\ \Fom\,=\,\Fom\ \times
$$
$$
R_{1,n+1}(x\com d_n(\Om))\ts\ldots\ts R_{12}(x\com d_1(\Om))\,\,
\Rt_{12}(x\com d_1(\Om))\ts\ldots\ts\Rt_{1,n+1}(x\com d_n(\Om))\,.
$$

\medskip\noindent
In view of (\ref{5.225}),(\ref{5.226}) and (\ref{5.333}),
this equality proves Proposition 1.7.
\qed


\smallskip\medskip\noindent\textbf{5.4.}
In this and the next two subsections, we prove Theorem 1.8. 
We use the same method as in the proof of Theorem~1.6.
As in Subsection 4.4, fix a standard tableau $\La$ of shape 
$\la\ts$, such that the tableau $\Om$ is obtained from $\La$
by removing the boxes with the numbers $1\lcd m\ts$. Here $\la$
is a\ partition of $l$, and $m=l-n$.
For $k=1\lcd l$ write 
$
\textstyle
d_k\ts=\,c_k(\La)\mp\frac12\ .
$
The standard tableau of shape $\mu\ts$, 
obtained by removing the boxes
with numbers $m+1\lcd l$ from the tableau $\La\ts$, is denoted by $\Up$. 

Put $L=N+M$. Consider the vector space $\CC^{L}=\CC^N\ns\op\ts\CC^M$.
In the present subsection, the rational functions
$$
R_{12}(x\com y)\ts\lcd\ts R_{1,l+1}(x\com y)
\quad\text{and}\quad
\Rt_{12}(x\com y)\ts\lcd\ts\Rt_{1,l+1}(x\com y)
$$
will take values in the algebra
$(\ts\End(\CC^L))^{\ts\ot\ts(l+1)}\ts$.
Take the linear operator $F_\La$
on the vector space $\CLl$, this operator is defined by Proposition 3.2.

\begin{Lemma}
We have the equality in $(\ts\End(\CC^L))^{\ts\ot\ts(l+1)}\ts$
$$
\Rt_{1,l+1}(x\com d_l)\ts\ldots\ts\Rt_{12}(x\com d_1)\,\ts
R_{12}(x\com d_1)\ts\ldots\ts R_{1,l+1}(x\com d_l)
\ts\cdot\ts(\ts1\ot F_\La)
$$

\vskip-20pt
\begin{equation}\label{5.51}
=\ 
\biggl(
1\ts-\ts\sum_{k=1}^l\,\frac{P_{\ts1,k+1}-Q_{\ts1,k+1}}
{x\pm\textstyle\frac12}
\ts\biggr)
\cdot\ts(\ts1\ot F_\La)\,.
\end{equation}
\end{Lemma}

\begin{proof}
Due to Proposition 3.2, the operator $F_\La$ is divisible
on the left by the operator $E_\La\ts$. Using Proposition 2.5 twice,
we can rewrite the product at the left hand side of the 
desired equality (\ref{5.51}) as
\begin{equation}\label{5.52}
\biggl(
1\ts+\ts\sum_{k=1}^l\,\frac{Q_{\ts1,k+1}}{x\mp\textstyle\frac12}
\ts\biggr)
\biggl(
1\ts-\ts\sum_{k=1}^l\,\frac{P_{\ts1,k+1}}{x\pm\textstyle\frac12}
\ts\biggr)
\cdot\ts(1\ot F_\La)\,\ts;
\end{equation}
see (\ref{4.6}). Due to Proposition 3.3, for any distinct
indices $i\com j\in\{1\lcd l\}$
$$
Q_{\ts1,i+1}\ts P_{\ts1,j+1}\cdot(1\ot F_\La)=
P_{\ts1,j+1}\ts Q_{\ts j+1,i+1}\cdot(1\ot F_\La)=0\,.
$$
Therefore, by using the relations
$$
Q_{\ts1,k+1}\ts P_{\ts1,k+1}=\pm\,Q_{\ts1,k+1}
\quad\text{for each}\quad
k=1\lcd l
$$
we can rewrite the product (\ref{5.52}) as at the
right-hand side of (\ref{5.51}).
\qed
\end{proof}

The equality (\ref{5.51})
is the starting point for our proof of Theorem~1.8. 
First, consider the case when $M=0$ and
$\mu=(0\com0\ts,\ts\ldots\ts)\ts$. In this case $N=L$, $n=l$
and the standard tableau $\Om=\La$ has a non-skew shape.
In this case, the left-hand side of (\ref{5.51})
describes the action of the twisted
Yangian $\YSL$ on its module $W_\Om\ts(0)=W_\La(0)\ts$;
see (\ref{5.225}) and the end of Subsection 1.7.

The right-hand side of the equality (\ref{5.51})
describes the action of the extended twisted Yangian 
$\XL$ on its module $\Wlm=W_\la\ts$,
defined in Subsection~1.8\ts; see (\ref{5.23}).
Indeed, here we have $g_\mu(x)=1$ and $\be_{NM}=\be_N$. 
The image of the operator
$F_\Om\ts(0)=F_\La\,$, as a $\g_L\ts$-submodule in $\CLl$,
is equivalent to the $\g_L\ts$-module~$W_\la\ts$. 
Thus the equality (\ref{5.51}) shows that in the case $M=0$
the action of the algebra $\XL$ on $W_\la$ factors through the 
homomorphism $\pi_L:\XL\to\YSL\ts$; see (\ref{piN}).
The equality (\ref{5.51}) also shows that
$W_\La$ and $W_\la$ are equivalent as $\YSL\ts$-modules.

Now suppose that $M\ge1$.
In Subsection 3.4, we introduced a projector
\begin{equation}\label{5.53}
J_m:\,\CLl\,\to\,\CMm_{\,\ts0}\ot\CNn\,.
\end{equation}
\noindent
Note that 
the linear operator (\ref{5.53}) is $G_M\times G_N\ts$-equivariant
by definition.

Take the tensor product of the evaluation $\YL\ts$-modules
with the parameters $d_1\lcd d_l\ts$. Consider the restriction
of this tensor product to the subalgebra $\YSL\subset\YL\ts$.
Pull this restriction back through the homomorphism
$\ts\pi_L\circ\ts\eta_L\ts:\ts\XL\to\YSL\ts$,
see (\ref{piN}) and (\ref{1.751}).
The action of $\XL$ on this module is described by the assignment 
$$
\sum_{i,j=1}^L\,E_{ij}\ot S_{ij}(x)\,\ts\mapsto\,\ts
\textstyle
R_{1,l+1}(-x\ns-\ns\frac{L}2\ts\com\ts d_l)^{-1}\ldots\,
R_{12}(-x\ns-\ns\frac{L}2\ts\com\ts d_1)^{-1}
$$

\vskip-16pt
\begin{equation}\label{wlup}
\textstyle
\times\ \,
\Rt_{12}(-x\ns-\ns\frac{L}2\ts\com\ts d_1)^{-1}\ldots\,
\Rt_{1,l+1}(-x\ns-\ns\frac{L}2\ts\com\ts d_l)^{-1}
\end{equation}
$$
\textstyle
=\ \,f(\ts x\ns+\ns\frac{L}2\ns\mp\ns\frac12\ts)\,\cdot\,
R_{1,l+1}(x\ns+\ns\frac{L}2\ts\com\ts-d_l)\,\ldots\,
R_{12}(x\ns+\ns\frac{L}2\ts\com\ts-d_1)
$$
\begin{equation}\label{wl}
\textstyle
\times\ \,
\Rt_{12}(x\ns-\ns\frac{L}2\ts\com\ts-d_1)\,\ldots\,
\Rt_{1,l+1}(x\ns-\ns\frac{L}2\ts\com\ts-d_l)\,,
\end{equation}

\vskip12pt\noindent
see (\ref{5.225}). 
Here the function $f(x)\in\CC(x)$ is defined by (\ref{fx})\ts;
in the notation of Subsection 4.4 we have $d_k=c_k\mp\frac12$
for every $k=1\lcd l$. Denote by $W_l$
the restriction of this $\XL\ts$-module
to the subalgebra $\XN$ of $\XL\ts$; here
we use the natural embedding $\psi_M:\XN\to\XL\ts$.
That is, $\psi_M:\ts S_{ij}(x)\ts\mapsto\ts S_{ij}(x)$ 
for $1\le i\com j\le N$.$\phantom{{}_{()}}$

Further, denote by $\Wal$ the $\XN\ts$-module
obtained by pulling the $\XN\ts$-module $W_l$ back through
the automorphism $\eta_N$ of $\XN\ts$.
Note that the vector space of the $\YN\ts$-modules
$W_l$ and $\Wal$ is $\CLl$.

Take the vector space $\CMm_{\,\ts0}\ot\CNn$.
For $k=m+1\lcd l\/$ we will keep using the notation (\ref{PV}).
But from now on we will regard
$P_{1,k+1}^{\,\wedge}$ as an operator on the subspace
\begin{equation}\label{5.54}
\CC^N\ot\CMm_{\,\ts0}\ot\CNn\subset\CC^N\ot\CMm\ot\CNn\,.
\end{equation}
For $k=m+1\lcd l\/$ consider the operator on
$\CC^N\ot\CMm\ot\CNn$, acting as $Q(N)$ in the first
and $(k+1)\ts$-th tensor factors, and acting as the identity
in the remaining $l-1$ tensor factors. 
The restriction of this operator
to the subspace  (\ref{5.54}) will be denoted by $Q_{1,k+1}^{\,\wedge}$.
Using this notation, put
$$
\Rt_{1,k+1}^{\,\,\wedge}(x\com y)\,=\,\ts
1+\frac{\,Q_{1,k+1}^{\,\wedge}\,}{x+y}\,.
$$

One can define
an action of $\XN$ on the space $\CMm_{\,\ts0}\ot\CNn$ by the assignment

\vskip-16pt
$$
\sum_{i,j=1}^N\,E_{ij}\ot S_{ij}(x)\,\ts\mapsto
$$
$$
\textstyle
f(\ts x\ns+\ns\frac{L}2\ns\mp\ns\frac12\ts)\,\cdot\,
R_{1,l+1}^{\,\,\wedge}(x\ns+\ns\frac{L}2\ts\com\ts-d_l)\,\ldots\,
R_{1,m+2}^{\,\,\wedge}(x\ns+\ns\frac{L}2\ts\com\ts-d_{\ts m+1})
$$
\begin{equation}\label{wmn}
\textstyle
\times\ \,
\Rt_{1,m+2}^{\,\,\wedge}(x\ns-\ns\frac{L}2\ts\com\ts-d_{\ts m+1})\,\ldots\,
\Rt_{1,l+1}^{\,\,\wedge}(x\ns-\ns\frac{L}2\ts\com\ts-d_l)\,;
\end{equation}

\vskip6pt\noindent
see below for the proof of this assertion.
Denote by $W_{mn}$ the $\XN\ts$-module defined 
by the assignment (\ref{wmn}).
Denote by $\Wamn$ the $\XN\ts$-module
obtained by pulling the $\XN\ts$-module $W_{mn}$ back through
the automorphism $\eta_N$ of $\XN\ts$. 
Using the function $h(x)\in\CC(x)$ defined by (\ref{gx}),
the action of $\XN$ on $\Wamn$
is described by
$$
\sum_{i,j=1}^N\,E_{ij}\ot S_{ij}(x)\,\ts\mapsto
$$

\vskip-16pt
$$
\textstyle
h(\ts x\ns-\ns\frac{M}2\ns\pm\ns\frac12\ts)\,\cdot\ts
\Rt_{1,l+1}(x\com d_l\ns+\ns\frac{M}2)\ts\ldots\ts
\Rt_{1,m+2}(x\com d_{\ts m+1}\!+\ns\frac{M}2)\,\ts
$$
$$
\textstyle
\times\ \,
R_{1,m+2}(x\com d_{\ts m+1}\!+\ns\frac{M}2)\ts\ldots\ts 
R_{1,l+1}(x\com d_l\ns+\ns\frac{M}2)\,.
$$

\vskip6pt
Note that the $\XN\ts$-module $\Wamn$
can also be obtained as follows. 
Take the tensor product of the evaluation $\YN\ts$-modules
with parameters
$$
\textstyle
d_{\ts m+1}\ns+\ns\frac{M}2\ts=\ts d_1(\Om)
\,\lcd\, 
d_l\ns+\ns\frac{M}2\ts=\ts d_n(\Om)\,.
$$
Consider the restriction
of this tensor product to $\YS\subset\YN\ts$.
Then regard this restriction as $\XN\ts$-module by using the
homomorphism $\pi_N:\XN\to\YS\ts$. Pull this $\XN\ts$-module
back through the automorphism
\begin{equation}\label{autom}
\textstyle
S_{ij}(x)\,\mapsto\,
h(\ts x\ns-\ns\frac{M}2\ns\pm\ns\frac12\ts)\cdot S_{ij}(x)
\end{equation}

\vskip4pt\noindent
of the algebra $\XN$.
The vector space of the resulting $\XN\ts$-module is $\CNn$.
Then by regarding $\CMm_{\,\ts0}\ot\CNn$ as
$\XN\ts$-module where every element of $\XN$ acts on $\CMm_{\,\ts0}$
trivially, we obtain the $\XN\ts$-module $\Wamn\ts$.

\begin{Proposition}
The projector {\rm (\ref{5.53})} is
an intertwiner of\/ $\XN\ts$-modules $W_l\to W_{mn}\ts$.
\end{Proposition}

\begin{proof}
Expand the product 
at the right hand-side of the equality (\ref{wl}) as the sum

\vskip-16pt
$$
\sum_{i,j=1}^L\,E_{ij}\ot B_{ij}(x)
$$
for certain rational functions
$B_{ij}(x)\in(\ts\End(\CC^L))^{\ts\ot\ts l}(x)\ts$.
It suffices to show that the sum

\vskip-16pt
\begin{equation}\label{EJB}
\sum_{i,j=1}^N\,E_{ij}\ot(\ts J_m\,B_{ij}(x))
\end{equation}
is equal to the product at the right-hand side of the assignment
(\ref{wmn}), multiplied by $J_m$ on the right. To show this, let us
expand the right-hand side of the equality
(\ref{wl}) as the sum of the products
\begin{equation}\label{PPQQ}
P_{\ts 1i_a}\ldots P_{\ts 1i_1}\ts
Q_{\ts 1j_1}\ldots Q_{\ts 1j_b}
\end{equation}
with coefficients from $\CC(x)\ts$;
here the sum is taken over all subsequences $i_1\lcd i_a$
and $j_1\lcd j_b$ in the sequence $2\lcd l+1\ts$. 
Consider four cases.

1) Suppose that $j_1\le m+1$. Also suppose that $a=0$ or $i_1>m+1$.
In this case, the product (\ref{PPQQ}) is divisible on the left 
by $Q_{\ts 1j_1}$ or by $Q_{\ts i_1j_1}$. 
The product (\ref{PPQQ}) does not contribute to the sum (\ref{EJB})
in this case, because the subspaces $\CC^M$ and $\CC^N$ in $\CC^L$
are orthogonal.

2) Suppose that $i_1\le m+1$ and $b=0$.
In this case,
the product (\ref{PPQQ}) does not contribute to the sum (\ref{EJB})\ts; 
see the proof of Lemma 2.5.

3) Suppose that $i_1\le m+1$ and $b>0$.
Then the product (\ref{PPQQ}) is divisible on the left 
by $Q_{\ts i_1k}$ for some index $k\in\{1\lcd l+1\}\ts$.
If $k=1$ or $k>m+1$, then
(\ref{PPQQ}) does not contribute to the sum (\ref{EJB}), 
because the subspaces $\CC^M$ and $\CC^N$ in $\CC^L$
are orthogonal. If $1<k\le m+1$, then 
(\ref{PPQQ}) does not contribute to (\ref{EJB}),
because the subspace $\CMm_{\,\ts0}\subset(\CC^L)^{\ts\ot\ts m}$
consists of traceless tensors with respect to the bilinear form
$\langle\ ,\,\rangle$ on $\CC^L$.

4) It remains to consider the case, when both $i_1\lcd i_a$
and $j_1\lcd j_b$ are subsequences in 
the sequence $m+2\lcd l+1\ts$. 
Suppose this is the case.
Write (\ref{PPQQ}) as the sum 

\vskip-16pt
$$
\sum_{i,j=1}^L\,E_{ij}\ot B_{ij}
$$
for certain elements $B_{ij}\in(\ts\End(\CC^L))^{\ts\ot\ts l}\ts$.
If $b>0$, the product (\ref{PPQQ}) can also be written as

\vskip-16pt
$$
P_{\ts 1i_a}\ldots P_{\ts 1i_1}\ts
P_{\ts j_{b-1}j_b}\ldots P_{\ts j_1j_2}\ts Q_{\ts 1j_b}\ts.
\hskip-20pt
$$

\vskip4pt\noindent
Using this observation and Lemma 2.5, we prove that in our remaining case
$$
\sum_{i,j=1}^N\,E_{ij}\ot J_m\,B_{ij}\,=\,
P_{\ts 1i_a}^{\ts\wedge}\ldots P_{\ts 1i_1}^{\ts\wedge}\ts
Q_{\ts 1j_1}^{\ts\wedge}\ldots Q_{\ts 1j_b}^{\ts\wedge}\,J_m\,.
\quad\qed
$$
\end{proof}

\begin{Corollary}
The projector {\rm (\ref{5.53})} is
an intertwiner of\/ $\XN\ts$-modules $\Wal\to\Wamn\ts$.
\end{Corollary}


\noindent\textbf{5.5.}
Let us continue our proof of Theorem 1.8.
Consider the image of the subspace (\ref{3.41})
under the operator $F_\La$ on $\CLl$. Note that
this image is contained in the subspace $W_\La\subset\CLl$.
Consider the $\XN\ts$-module $W_l$ defined in Subsection 5.4\ts;
the vector space of this module is $\CLl$.

\begin{Proposition}
The image of the subspace\/ {\rm (\ref{3.41})} under the operator $F_\La$
is an\/ $\XN$-submodule of $W_l\ts$.
\end{Proposition}

\begin{proof}
The action of the coefficients of the series 
$S_{ij}(x)$ with $1\le i\com j\le N$
on the $\XN\ts$-module $W_l$ is described by the assignment (\ref{wlup}). 
Here we use the natural embedding $\psi_M:\XN\to\XL\ts$. 
Consider the product in the algebra
$(\ts\End(\CC^L))^{\ts\ot\ts(l+1)}(x)\ts$,
displayed at the right-hand side of the equality (\ref{wl}).
We have a relation in this algebra,
$$
\textstyle
R_{1,l+1}(x\ns+\ns\frac{L}2\ts\com\ts-d_l)\,\ldots\,
R_{12}(x\ns+\ns\frac{L}2\ts\com\ts-d_1)
$$
$$
\textstyle
\times\ \,
\Rt_{12}(x\ns-\ns\frac{L}2\ts\com\ts-d_1)\,\ldots\,
\Rt_{1,l+1}(x\ns-\ns\frac{L}2\ts\com\ts-d_l)\ts\cdot\ts
(1\ot F_\La)\ =\ (1\ot F_\La)
$$
$$
\textstyle
\times\ \,
\Rt_{1,l+1}(x\ns-\ns\frac{L}2\ts\com\ts-d_l)\,\ldots\,
\Rt_{12}(x\ns-\ns\frac{L}2\ts\com\ts-d_1)\,\ts
$$
\begin{equation}\label{5.8}
\textstyle
\times\ \,
R_{12}(x\ns+\ns\frac{L}2\ts\com\ts-d_1)\,\ldots\,
R_{1,l+1}(x\ns+\ns\frac{L}2\ts\com\ts-d_l)\,;
\end{equation}

\vskip4pt\noindent
see Subsection 5.3. Consider the $2\ts l$ factors
displayed in the last two lines of the relation (\ref{5.8}).
Expand the product of these factors as the sum 
$$
\sum_{i,j=1}^L\,E_{ij}\ot C_{ij}(x)
$$
where $C_{ij}(x)\in(\ts\End(\CC^L))^{\ts\ot\ts l}(x)\ts$.
Consider the restrictions of the operator values of the
functions $C_{ij}(x)$ with $1\le i\com j\le N$
to the subspace (\ref{3.41}). By an argument similar to the one
used in the proof of Proposition 5.4,
$$
\sum_{i,j=1}^N\,E_{ij}\ot (\ts C_{ij}(x)\,|\,\CMm_{\,\ts0}\ot\CNn)
$$

\vskip-8pt
$$
\textstyle
=\ \,
\Rt_{1,l+1}^{\,\,\wedge}(x\ns-\ns\frac{L}2\ts\com\ts-d_l)\,\ldots\,
\Rt_{1,m+2}^{\,\,\wedge}(x\ns-\ns\frac{L}2\ts\com\ts-d_1)
$$
$$
\textstyle
\times\ \,
R_{1,m+2}^{\,\,\wedge}(x\ns+\ns\frac{L}2\ts\com\ts-d_l)\,\ldots\,
R_{1,l+1}^{\,\,\wedge}(x\ns+\ns\frac{L}2\ts\com\ts-d_l)\,.
$$

\vskip6pt\noindent
In particular, the operator values of the functions $C_{ij}(x)$
with $1\le i,j\le N$ preserve the subspace (\ref{3.41}).
Now Proposition~5.5 follows from (\ref{5.8}).
\qed
\end{proof}

The $\XN\ts$-module $\Wal$ has been obtained from $W_l$
by pulling back through an automorphism of $\XN\ts$. So
Proposition 5.5 has a corollary.

\begin{Corollary}
The image of the subspace\/ {\rm (\ref{3.41})} under the operator $F_\La$
is an\/ $\XN$-submodule of $\Wal\ts$.
\end{Corollary}


\noindent\textbf{5.6.}
In this subsection we will complete the proof of Theorem 1.8.
Consider the image of the subspace (\ref{3.41}) under the linear operator
$$
J_m\ts F_\La:\CLl\to\CMm_{\,\ts0}\ot\CNn\,.
$$ 
Due to Proposition 3.3, this image coincides with the vector subspace
\begin{equation}\label{5.9}
W_{\ts\Up}\ot\Wom\,\subset\,\CMm_{\,\ts0}\ot\CNn\ts.
\end{equation}
Indeed,
$$
J_m\ts F_\La\ts\,|\,\ts\CMm_{\,\ts0}\ot\CNn\ts=\ts\,
(\ts E_{\ts\Up}\,|\,\ts\CMm_{\,\ts0})\ot\Fom\,;
\hskip-30pt
$$

\vskip4pt
\noindent
see (\ref{3.46}) and (\ref{3.4444}). It now follows from
Corollaries 5.4 and 5.5 that the vector subspace (\ref{5.9})
is a submodule in the $\XN\ts$-module $\Wamn\ts$. 
Let us denote this submodule of $\Wamn$ by $W$.
Note that $W$ is a subquotient of the $\XN\ts$-module $\Wal$
by definition.

The description of the $\XN\ts$-module $\Wamn$ given 
just before stating Proposition 5.4, 
yields the following description of the $\XN\ts$-module $W$. 
Take the $\YS\ts$-module $\Wom$ as defined in Subsection 1.7. 
Then regard $\Wom$ as $\XN\ts$-module by using the homomorphism
$\pi_N\ts$. Pull the $\XN\ts$-module $\Wom$ back through the
automorphism (\ref{autom})
of the algebra $\XN\ts$. 
Now extend the resulting action of $\XN$ on the vector
space $\Wom\ts$, to the vector space $W_{\ts\Up}\ot\Wom\ts$ so that
every element of $\XN$ acts on $W_{\ts\Up}$
as the identity. Then we obtain the $\XN\ts$-module $W$.
Note that the subspace $W_{\ts\Up}$ of 
the vector space (\ref{3.000})
is equivalent to $W_\mu$ as a representation of $G_M\ts$;
see Proposition~3.3.

Consider the vector subspace $W_\La\subset\CLl$ as a
representation of $G_L\ts$,
equivalent to $W_\la$. Then regard $W_\La$ as $\XN\ts$-module
by pulling back through the homomorphism
$\be_{NM}:\XN\to\UgL\ts$, see (\ref{beNM}).

\begin{Proposition}
The\/ $\XN$-module\/ $W$ is a subquotient of the\/ $\XN$-module\/ $W_\La$.
\end{Proposition}

\begin{proof}
By (\ref{beNM}), we have
$\be_{NM}=\,\be_L\circ\ts\eta_L\circ\ts\psi_M\circ\ts\eta_N\ts$.
Consider $W_\La$ as a submodule in the restriction of
the tensor product of
the evaluation $\YL\ts$-modules with the parameters
$d_1\lcd d_l$ to $\YSL\subset\YL\ts$. Then 
regard $W_\La$ as a $\XL\ts$-module, using the homomorphism $\pi_L\ts$.
We have already shown that the resulting action of
$\XL$ in $W_\La$ factors through the homomorphism $\be_L:\XL\to\UgL\ts$.
Hence the $\XN\ts$-module $W_\La$ as defined above
can also be obtained by pulling the 
just determined action of $\XL$ on $W_\La\ts$,
back through the injective homomorphism
$$
\eta_L\circ\ts\psi_M\circ\ts\eta_N:\XN\to\XL\,.
$$
Thus $W_\La$ is a submodule in the $\XN\ts$-module $\Wal$.
But by definition, $W$ is a quotient of a certain $\XN\ts$-submodule of
$\Wal$. The latter submodule of $\Wal$ is contained in $W_\La\,$.
\qed
\end{proof}

We can now complete our proof of Theorem 1.8.
Let us consider the restriction of the representation $W_\La$ of 
the group $G_L$
to the subgroup $G_M\ts$. Realize the vector space (\ref{1.4}) as
\begin{equation}\label{5.10}
{\rm Hom}_{\,G_M}(\ts W_{\ts\Up}\ts\com W_\La\ts)\ts.
\end{equation}
Since the image of the homomorphism $\be_{NM}$ is contained in
the subalgebra of $G_M\ts$-invariants $\BMN\subset\UgL$, the action
of the algebra $\XN$ on $W_\La$ induces an action of $\XN$ on (\ref{5.10}).
If $G_L=Sp_L\ts$, then this action of $\XN$ on (\ref{5.10})
is irreducible, see \cite[Section 4]{MO}. If $G_L=O_L$,
then (\ref{5.10}) is irreducible under the joint action of
the algebra $\XN$ and the subgroup $G_N\subset G_L\ts$.

The operator (\ref{5.53}) is $G_M\times G_N\ts$-equivariant,
and the vector space $W_\Up\ot\Wom$ of the $\XN\ts$-module $W$
comes with the natural action 
of the group $G_M\times G_N\ts$. The action of $G_M$ on $W$
commutes with the action of the algebra $\XN\ts$.
By Proposition 5.6, the $\XN\ts$-module 
\begin{equation}\label{5.11}
{\rm Hom}_{\,G_M}(\ts W_{\ts\Up}\ts\com W\ts)
\end{equation}
is a subquotient of the $\XN\ts$-module (\ref{5.10}).
It is also a subquotient of (\ref{5.10}) as a representation
of the group $G_N\ts$. Since (\ref{5.10}) is irreducible
under the joint action of $\XN$ and $G_N$,
this action must be equivalent to the joint action of $\XN$ and $G_N$
on \text{(\ref{5.11})\ts;} see Proposition~3.5.

The $\XN\ts$-module (\ref{5.11}) can also be obtained
in the following way. Take the 
$\YS\ts$-module $\Wom$ as defined in Subsection 1.7.
Regard $\Wom$ as $\XN\ts$-module by using the homomorphism $\pi_N\ts$.
Finally, pull this $\XN\ts$-module back through the automorphism
(\ref{1.62}) of $\XN\ts$, where
$$
\textstyle
g(x)=g_\mu(\ts x-\frac{M}2\pm\frac12\ts)^{-1}\ts.
$$

\vskip4pt\noindent
Here we use Lemma 4.4.
This explicit description of the $\XN\ts$-module (\ref{5.11})
completes the proof of Theorem 1.8. We also obtain
Proposition 1.4. 


\smallskip\medskip\noindent\textbf{5.7.}
In this subsection, we
prove analogues of the results of Subsection~4.3 for
the twisted Yangian $\YS\ts$.
Due to (\ref{1.751}) and (\ref{1.752}), for
the elements $S_{ij}^{(1)}\in\XN$ we have
$\eta_N(S_{ij}^{(1)})=S_{ij}^{(1)}$
and
$\beta_N(S_{ij}^{(1)})=-E_{ji}-\si(E_{ji})\ts$;
here the matrix unit $E_{ji}$ is regarded as a generator of
the algebra $\US\ts$. 

Therefore for any non-negative integer $M$,
by the definition (\ref{1.769})
of the homomorphism $\be_{NM}:\XN\to\UgMN$ we obtain the equality
\begin{equation}\label{rhs}
\be_{NM}(S_{ij}^{(1)})=-E_{ji}-\si(E_{ji})\,;
\end{equation}
at the right-hand side of this equality we have
an element of the subalgebra $\UgN\subset\UgMN\ts$.
This equality shows that the element $S_{ij}^{(1)}\in\XN$
acts as $-E_{ji}-\si(E_{ji})$ on the $\XN\ts$-module $\Wlm\ts$.
Here we use the fact that  
the coefficient of the series (\ref{gmM}) at $x^{-1}$ is zero;
see also (\ref{1.62}).

Consider the surjective homomorphism $\pi_N:\XN\to\YS\ts$,
and the embedding $\US\to\YS$ defined by (\ref{1.844}).
Using the first statement of Theorem 1.8, we can regard $\Wlm$
as $\YS\ts$-module.
The right-hand side of the equality (\ref{rhs}),
as an element of the subalgebra $\US\subset\YS\ts$,
coincides with $\pi_N(S_{ij}^{(1)})\,$; see (\ref{1.844}).
Independently of Theorem 1.8, this coincidence shows
that the action of the elements $S_{ij}^{(1)}$
on the $\XN\ts$-module $\Wlm$ 
factors through the homomorphism $\pi_N$. 
Moreover, this coincidence shows that the natural action
of the algebra $\US$ on $\Wlm$ coincides with its action
as a subalgebra in $\YS\ts$.

Let us now consider the $\YS\ts$-module $\Wom\ts$. It is 
a submodule in the restriction of the
tensor product (\ref{1.7777777}) of evaluation $\YN\ts$-modules
to the subalgebra $\YS\subset\YN\ts$, see also (\ref{1.777}).
Observe that the embedding $\US\to\YS$ as defined by (\ref{1.844})
can also be obtained by restricting the embedding 
$\UN\to\YN$ to $\YS\subset\YN\,$; see (\ref{4.4}) and (\ref{piN}).
Here we use the equality (\ref{sisi})
in $\glN\ns\ot\glN\ts$.
But the action of the subalgebra
$\UN\subset\YN\ts$ on the $\YN\ts$-module (\ref{1.7777777})
coincides with the natural action of $\UN$
on the vector space $\CNn\ts$; see Subsection 4.4.
Therefore the natural action of $\UgN$ on $\Wom$ coincides with its action
as a subalgebra in $\YS\ts$.


\begin{acknowledgement}
I am grateful to Ivan Cherednik,
Grigori Olshanski, Evgeny Sklyanin and Vitaly Tarasov
for numerous conversations.
The present work was supported by the European Commission under the
grant ERB-FMRX-CT97-0100. 
\end{acknowledgement}




\begin{thebibliography}{MNO}

\bibitem[B]{B}
{A.\,Berele},
\textit{Construction of Sp-modules using tableaux},
{Linear and Multilinear Algebra}
\textbf{19}
(1986),
299--307.

\bibitem[C1]{C1}
{I.\,Cherednik},
\textit{Factorized particles on the half-line and root systems},
{Theor.\ Math. Phys.}
\textbf{61}
(1984),
977--983.

\bibitem[C2]{C2}
{I.\,Cherednik},
\textit{On special bases of irreducible finite-dimensional representations
of the degenerate affine Hecke algebra},
{Funct.\ Analysis Appl.}
\textbf{20}
(1986),
87--89.

\bibitem[C3]{C3}
{I.\,Cherednik},
\textit{A new interpretation of Gelfand--Zetlin bases},
{Duke Math.\ J.}
\textbf{54}
(1987),
563--577.

\bibitem[D]{D}
{V.\,Drinfeld},
\textit{Hopf algebras and the quantum Yang--Baxter equation},
{Soviet Math.\ Dokl.}
\textbf{32}
(1985),
254--258.

\bibitem[GZ1]{GZ1}
{I.\,Gelfand and M.\,Zetlin},
\textit{Finite-dimensional representations of the group
of uni\-modular matrices},
{Dokl.\ Akad.\ Nauk SSSR} 
\textbf{71}
(1950),
825--828.

\bibitem[GZ2]{GZ2}
{I.\,Gelfand and M.\,Zetlin},
\textit{Finite-dimensional representations of groups
of orthogonal matrices},
{Dokl.\ Akad.\ Nauk SSSR}
\textbf{71}
(1950),
1017--1020.

\bibitem[KS]{KS}
{R.\,King and N.\,El-Sharkaway},
\textit{Standard Young tableaux and weight multiplicities
of the classical Lie groups},
{J.\ Phys.}
\textbf{A16}
(1983),
3153--3177.  

\bibitem[KW]{KW}
{R.\,King and T.\,Welsh},
\textit{Construction of orthogonal group modules using tableaux},
{Linear and Multilinear Algebra}
\textbf{33}
(1993),
251--283.

\bibitem[KT]{KT}
{K.\,Koike and I.\,Terada},
\textit{Young-diagrammatic methods for the representation theory
of the classical groups of type $B_n\com C_n\com D_n$},
{J.\ Algebra}
\textbf{107}
(1987),
466--511.

\bibitem[M]{M}
{I.\,Macdonald},
\textit{Symmetric Functions and Hall Polynomials},
Clarendon Press, Oxford, 1995.

\bibitem[M1]{M1}
{A.\,Molev},
\textit{A basis for representations of symplectic Lie algebras},
{Commun.\ Math. Phys.}
\textbf{201} 
(1999),
591--618.

\bibitem[M2]{M2}
{A.\,Molev},
\textit{A weight basis for representations of even orthogonal Lie algebras},
{Adv. Stud.\ Pure\ Math.}
\textbf{28}
(2000),
221--240.

\bibitem[M3]{M3}
{A.\,Molev},
\textit{Weight bases of Gelfand\,--Tsetlin type for representations
of classical Lie algebras},
{J.\ Phys.}
\textbf{A\,33}
(2000),
4143--4158.  

\bibitem[MN]{MN}
{A.\,Molev and M.\,Nazarov},
\textit{Capelli identities for classical Lie algebras},
{Math. Ann.}
\textbf{313}
(1999),
315--357.

\bibitem[MNO]{MNO}
{\hskip2pt\hskip-4pt A.\,Molev, M.\,Nazarov and G.\,Olshanski},
\textit{Yangians and classical Lie algebras},
{Russian Math.\ Surveys}
\textbf{51} 
(1996),
205--282.

\bibitem[MO]{MO}
{A.\,Molev and G.\,Olshanski},
\textit{Centralizer construction for twisted Yangians},
{Selecta Math.}
\textbf{6}
(2000),
269--317.

\bibitem[N1]{N1}
{M.\,Nazarov},
\textit{Young's symmetrizers for projective representations of the symmetric
group},
{Adv.\ Math.}
\textbf{127}
(1997),
190--257.

\bibitem[N2]{N2}
{M.\,Nazarov},
\textit{Yangians and Capelli identities},
{Amer.\ Math.\ Soc.\ Translations}
\textbf{181}
(1998),
139--163.

\bibitem[N3]{N3}
{M.\,Nazarov},
\textit{Capelli elements in the classical universal enveloping algebras},
{Adv. Stud. Pure\ Math.}
\textbf{28}
(2000),
261--285.

\bibitem[N4]{N4}
{M.\,Nazarov},
\textit{Rational representations of Yangians associated 
with skew Young diagrams},
{Math.\,Z.}
\textbf{247}
(2004),
000--000.

\bibitem[NO]{NO}
{M.\,Nazarov and G.\,Olshanski},
\textit{Bethe subalgebras in twisted Yangians},
{Commun. Math.\ Phys.}
\textbf{178}
(1996),
483--506.

\bibitem[NT1]{NT1}
{M.\,Nazarov and V.\,Tarasov},
\textit{Representations of Yangians with Gelfand--Zetlin bases},
{J.\ Reine Angew.\ Math.}
\textbf{496}
(1998),
181--212.

\bibitem[NT2]{NT2}
{M.\,Nazarov and V.\,Tarasov},
\textit{On irreducibility of tensor products of Yangian modules
associated with skew Young diagrams},
{Duke Math. J.}
\textbf{112}
(2002),
342--378.

\bibitem[O1]{O1}
{G.\,Olshanski},
\textit{Extension of the algebra $U(g)$ for infinite-dimensional classical
Lie algebras $g$, and the Yangians $Y(gl(m))$},
{Soviet Math.\ Dokl.}
\textbf{36}
(1988),
569--573.

\bibitem[O2]{O2}
{G.\,Olshanski},
\textit{Twisted Yangians and infinite-dimensional classical Lie algebras},
{Lecture Notes in Math.}
\textbf{1510}
(1992),
103--120.

\bibitem[P]{P}
{R.\,Proctor},
\textit{Young tableaux, Gelfand patterns, and branching rules for
classical groups},
{J.\ Algebra}
\textbf{164}
(1994),
299--360.

\bibitem[S]{S}
{E.\,Sklyanin},
\textit{Boundary conditions for integrable quantum systems},
{J.\ Phys.}
\textbf{A21}
(1988),
2375--2389.  

\bibitem[W]{W}
{H.\,Weyl},
\textit{Classical Groups, their Invariants and Representations},
Princeton University Press, Princeton, 1946.

\bibitem[Y1]{Y1}
{A.\,Young},
\textit{On quantitative substitutional analysis I\ts} and \textit{II\/},
{Proc.\ London Math. Soc.}
\textbf{33}
(1901),
97--146  
and
\textbf{34}
(1902),
361--397.

\bibitem[Y2]{Y2}
{A.\ Young,}
\textit{On quantitative substitutional analysis VI\ts},
{Proc. London Math. Soc.}
\textbf{34}
(1932),
196--230.

\end{thebibliography}
\enddocument